\documentclass[11pt,draft,leqno]{article}
\usepackage{amsmath,amssymb,theorem}
\usepackage{latexsym}
\usepackage[all]{xy}
\usepackage{amsfonts}

\newtheorem{thm}{Theorem}[subsection]

\newtheorem{prop}[thm]{Proposition}
\newtheorem{lem}[thm]{Lemma}
\newtheorem{df}[thm]{Definition}
\newtheorem{cor}[thm]{Corollary}

{\theorembodyfont{\rmfamily} \newtheorem{rem}[thm]{Remark}}

{\theorembodyfont{\rmfamily} }

\newtheorem{thmi}{Theorem}
\newtheorem{cori}[thmi]{Corollary}

\newcommand{\opi}[1]{\operatorname{{\it #1}}}
\newcommand{\op}[1]{\operatorname{{#1}}}
\newcommand{\No}{\boldsymbol{N}}
\newcommand{\De}{\boldsymbol{D}}
\newcommand{\fp}{\boldsymbol{F}\boldsymbol{P}}
\newcommand{\bs}{{\sf{s}}}
\newcommand{\sch}[2]{\mbox{$({#1}\:\!\otimes\:\! {#2})^{\operatorname{sch}}$}}
\newcommand{\sff}[1]{{\sf{#1}}}

\numberwithin{equation}{subsection}

\begin{document}

\title{Schematic homotopy types and non-abelian Hodge theory}
\author{L. Katzarkov\thanks{Partially supported by NSF Career
Award DMS-9875383 and an A.P.Sloan research fellowship} \and
T. Pantev\thanks{Partially supported by
NSF Grants DMS-0099715, DMS-0403884 and an A.P.Sloan research
fellowship} \and B. To\"en} 
\date{}
\maketitle

\tableofcontents

\setcounter{section}{-1}

\section{Introduction}

Modern Hodge theory is concerned with the interaction between the
geometry and the topology of complex algebraic varieties. The general
idea is to furnish topological invariants of varieties with additional
algebraic structures that capture essential information about the
algebraic geometry of these varieties. The very first example of such
a result is the existence of the Hodge decomposition on the Betti
cohomology of a complex projective manifold (see e.g. \cite{gh}), and
more generally the existence of Deligne's mixed Hodge structures on
the cohomology of arbitrary varieties \cite{d-h2,d-h3}.  A more
sophisticated invariant was developed by P.Deligne, P.Griffiths,
J.Morgan and D.Sullivan \cite{dgms} who combined the Hodge theory of
P.Deligne and the rational homotopy theory of D.Quillen and D.Sullivan
to construct a Hodge decomposition on the rational homotopy type of
a complex projective manifold. Their construction was later extended by
J.Morgan who,  in the simply connected case, proved the existence of a
mixed Hodge structure on the complexified homotopy groups \cite{mo},
and by R.Hain who constructed a mixed Hodge structure on the Mal'cev
completion of the fundamental group \cite{ha}. More recently C.Simpson
constructed an action of the discrete group $\mathbb{C}^{\times
\delta}$ on the full pro-algebraic fundamental group of a smooth
projective complex variety that can be thought of as a Hodge
decomposition \cite{s2}.

In \cite{t1} a new homotopy invariant of a space $X$ has been
constructed, which is called the \emph{schematization} of $X$ and is
denoted by $\sch{X}{\mathbb{C}}$.  The main goal of the present paper
is to define a Hodge decomposition on $\sch{X}{\mathbb{C}}$ when $X$
is a projective manifold, and to show that this structure recovers all
of the Hodge structures on cohomology, rational homotopy groups and
completions of the fundamental groups mentioned above.  For this, the
various Hodge decompositions will be viewed as actions of the group
$\mathbb{C}^{\times \delta}$. For example, the Hodge decomposition on
the cohomology of a projective variety, $H^{n}(X,\mathbb{C})\simeq
\oplus H^{n-p}(X,\Omega_{X}^{p})$, can be understood via the action of
$\mathbb{C}^{\times \delta}$ which has weight $p$ on the direct
summand $H^{n-p}(X,\Omega_{X}^{p})$. As a manifestation of this
principle we construct an action of the group $\mathbb{C}^{\times
\delta}$ on the object $\sch{X}{\mathbb{C}}$.

In order to state this result, let us recall that for
any pointed connected space $X$, the schematic homotopy type
$\sch{X}{\mathbb{C}}$ is a pointed connected simplicial presheaf on
the site of affine complex schemes with the faithfully flat
topology. Furthermore, this simplicial presheaf satisfies the
following conditions (see $\S 1$ for details).

\begin{itemize}
\item The sheaf of groups $\pi_{1}(\sch{X}{\mathbb{C}},x)$ is
represented by the pro-algebraic completion
of the discrete group $\pi_{1}(X,x)$.

\item There exist functorial isomorphisms
  $H^{n}(\sch{X}{\mathbb{C}},\mathbb{G}_{a})\simeq 
H^{n}(X,\mathbb{C})$.

\item If $X$ is a simply connected finite $CW$-complex, there exist
functorial isomorphisms
$$\pi_{i}(\sch{X}{\mathbb{C}},x)\simeq \pi_{i}(X,x)\otimes
\mathbb{G}_{a}.$$
\end{itemize}

The object $\sch{X}{\mathbb{C}}$ is therefore a good
candidate for an unificator of the various homotopy invariants
previously studied by Hodge theory. Our main result is
the following:

\begin{thmi}
For any pointed smooth projective complex variety $X$, there exist an
action of $\mathbb{C}^{\times \delta}$ on $\sch{X}{\mathbb{C}}$,
functorial in $X$ and called the Hodge decomposition.  This action
recovers the usual Hodge decompositions on cohomology, on completed
fundamental groups, and (in the simply connected case) on complexified
homotopy groups.
\end{thmi}

We will now give an overview of the content of the present work.

\

\smallskip

\begin{center} \textit{Rational homotopy theory and schematic homotopy
types} \end{center}

\smallskip

It is well known that the rational homotopy type $X_{\mathbb{Q}}$
of an arbitrary space $X$ is a pro-nilpotent homotopy type. As a
consequence, homotopy invariants
of the space $X$ which are not of nilpotent nature are not accessible
through the space $X_{\mathbb{Q}}$. For example, the fundamental group of
$X_{\mathbb{Q}}$ is isomorphic to the Mal'cev completion of
$\pi_{1}(X,x)$, and in particular the full pro-algebraic
fundamental group of $X$ is beyond the scope of the techniques of
rational homotopy theory.

In order to develop a substitute of rational homotopy theory for non
nilpotent spaces, the third author introduced the notion of a
\textit{pointed schematic homotopy type} over a field $k$, as well as
a \textit{schematization functor} (see \cite{t1}).  When the base
field is of characteristic zero, a pointed schematic homotopy type $F$
is essentially a pointed and connected simplicial presheaf on the site
of affine $k$-schemes with the flat topology, whose fundamental group
sheaf $\pi_{1}(F,*)$ is represented by an affine group scheme, and
whose homotopy sheaves $\pi_{i}(F,*)$ are products (possibly infinite)
of copies of the additive group $\mathbb{G}_{a}$. The fundamental
group $\pi_{1}(F,*)$ of a schematic homotopy type can be any affine
group scheme, and its action on a higher homotopy group $\pi_{i}(F,*)$
can be an arbitrary algebraic representation. Furthermore, the
homotopy category of augmented and connected commutative differential
graded algebras (concentrated in positive degrees) is equivalent to
the full sub-category of pointed schematic homotopy types $F$ such
that $\pi_{1}(F,*)$ is a unipotent affine group scheme (see Theorems
\ref{t0} and \ref{t1}). In view of this,  pointed
schematic homotopy types are reasonable models for a generalization of
rational homotopy theory.

For any pointed connected homotopy type $X$ and any field $k$, it was
proved in \cite[Theorem $3.3.4$]{t1} that there exist a universal
pointed schematic homotopy type $\sch{X}{k}$, called the
\textit{schematization of $X$ over $k$}. The schematization is
functorial in $X$. Its universality is explicitly spelled out in
following two properties.

\begin{itemize}
\item The affine group scheme $\pi_{1}(\sch{X}{k},x)$ is the
pro-algebraic completion of
the group \linebreak $\pi_{1}(X,x)$. In particular, there is a one to one
correspondence between
finite dimensional linear representations $\pi_{1}(X,x)$ and
finite dimensional linear representations of the affine group
scheme $\pi_{1}(\sch{X}{k},x)$.

\item For any finite dimensional linear representation $V$ of
$\pi_{1}(X,x)$, also considered as a representation of the
group scheme $\pi_{1}(\sch{X}{k},x)$, there
is a natural isomorphism on cohomology groups with local coefficients
\[
H^{\bullet}(X,V) \simeq H^{\bullet}(\sch{X}{k},V).
\]
\end{itemize}

Moreover, for a simply connected finite $CW$-complex $X$, one can show
\cite[Theorem $2.5.1$]{t1} that the homotopy
sheaves $\pi_{i}(\sch{X}{\mathbb{Q}},x)$ are
isomorphic to $\pi_{i}(X,x)\otimes \mathbb{G}_{a}$, and that the
simplicial set of global sections of the simplicial
presheaf $\sch{X}{\mathbb{Q}}$ is a model for the rational
homotopy type of $X$. This fact justifies the use of the  schematization
functor as a generalization of the
rationalization functor to non-nilpotent spaces.

\

\smallskip

\begin{center} \textit{The Hodge decomposition} \end{center}

\smallskip

As shown in \cite[Section~4.1]{small}, when $X$ is the underlying topological 
space of a compact smooth 
manifold, its schematization $\sch{X}{\mathbb{C}}$ 
can be explicitly described using complexes of differential
forms with coefficients in flat connections on $X$. If $X$ is 
furthermore a complex projective manifold, this description
together with Simpson's non-abelian 
Hodge correspondence gives a model of
$\sch{X}{\mathbb{C}}$ in terms of Dolbeault complexes with coefficients 
in Higgs bundles. Since the group
$\mathbb{C}^{\times \delta}$ acts naturally on the category 
of Higgs bundles and on their Dolbeault complexes, one
gets a natural action of 
$\mathbb{C}^{\times \delta}$ on $\sch{X}{\mathbb{C}}$. The main properties of this $\mathbb{C}^{\times \delta}$ action are summarized in the following theorem.

\begin{thmi}(Theorem \ref{tdec})\label{ti1}
Let $X$ be a pointed smooth projective variety over $\mathbb{C}$, and
let
$$\sch{X^{\op{top}}}{\mathbb{C}}$$
denote the schematization of the underlying topological space of
$X$ (for the classical topology). Then, there exists an action of
$\mathbb{C}^{\times \delta}$ on
$\sch{X^{\op{top}}}{\mathbb{C}}$, called the Hodge decomposition, such
that the following conditions are satisfied.

\begin{enumerate}
\item The induced action of $\mathbb{C}^{\times \delta}$ on the
cohomology groups
$H^{n}(\sch{X^{\op{top}}}{\mathbb{C}},\mathbb{G}_{a})\simeq
H^{n}(X^{\op{top}},\mathbb{C})$ is compatible with the Hodge decomposition
in the following sense.  For any
$\lambda \in \mathbb{C}^{\times \delta}$,
and $y \in H^{n-p}(X,\Omega_{X}^{p}) \subset
H^{n}(X^{\op{top}},\mathbb{C})$ one has
$\lambda(y)=\lambda^{p}\cdot y$.
\item Let $\pi_{1}(\sch{X^{\op{top}}}{\mathbb{C}},x)^{\op{red}}$ be the
maximal reductive quotient of
$\pi_{1}(X^{\op{top}},x)^{\op{alg}}$. Then,
the induced action of $\mathbb{C}^{\times \delta}$ on
$\pi_{1}(\sch{X^{\op{top}}}{\mathbb{C}},x)^{\op{red}}\simeq
\pi_{1}(X^{\op{top}},x)^{\op{red}}$
is the one defined in \cite{s2}.
\item If $X^{\op{top}}$ is simply connected, then the induced action of
$\mathbb{C}^{\times \delta}$ on
\[
\pi_{i}(\sch{X^{\op{top}}}{\mathbb{C}},x)(\mathbb{C})\simeq
\pi_{i}(X^{\op{top}},x)\otimes \mathbb{C}
\]
is compatible with the Hodge decomposition defined in \cite{mo}
in the following
sense. Let $F^{\bullet}\pi_{i}(X^{\op{top}})\otimes \mathbb{C}$ be the
Hodge filtration 
defined in \cite{dgms,mo}, then
\[
F^{p}\pi_{i}(X^{\op{top}})\otimes \mathbb{C}=\{x \in
\pi_{i}(X^{\op{top}}\otimes \mathbb{C}) |
\exists \; q\geq p \text{ so that } \lambda(x)=\lambda^{q}\cdot x,\;  \forall
\lambda \in \mathbb{C}^{\times }\}.
\]
\item Let $R_{n}$ be the set of isomorphism classes of simple
$n$-dimensional linear
representations of $\pi_{1}(\sch{X^{\op{top}}}{\mathbb{C}},x)$.
Then, the induced action of $\mathbb{C}^{\times }$ on the set
\[
R_{n}\simeq
\op{Hom}(\pi_{1}(X^{\op{top}},x),Gl_{n}(\mathbb{C}))/Gl_{n}(\mathbb{C}) 
\]
defines a continuous action of the topological group $\mathbb{C}^{\times }$
(for the analytic topology).
\end{enumerate}
\end{thmi}

\

\smallskip

\begin{center} \textit{The weight spectral sequence} \end{center}

\smallskip

For a general space $X$, the schematization $\sch{X}{\mathbb{C}}$
turns out to be extremely difficult to compute. For instance,
essentially nothing is known about the higher homotopy groups
\linebreak $\pi_{i}(\sch{X}{\mathbb{C}},x)$.

To facilitate computations of this type we propose an approach based
on Curtis' spectral sequence in conventional topology
\cite{curtis,curtis2}. 
For a general pointed schematic homotopy type $F$, we will construct a
\emph{weight tower} $W^{(*)}F^{0}$, which is an algebraic analog of
Curtis' construction from \cite{curtis,curtis2}. To this tower we
associate a 
spectral sequence, called the \emph{weight spectral sequence}, going
from certain homology invariants of $F$ to its homotopy groups.  In
general, this spectral sequence will not converge. However, when $X$
is a smooth projective complex manifold, we can use the existence of
the Hodge decomposition in order to prove that the weight spectral
sequence of $\sch{X^{\op{top}}}{\mathbb{C}}$ degenerates at $E_{2}$
and that its $E_{\infty}$ term computes the higher homotopy groups of
$\sch{X^{\op{top}}}{\mathbb{C}}$.  More precisely, suppose that
$(X,x)$ is a pointed smooth and projective complex manifold with
schematization $F :=(X^{\op{top}}\otimes\mathbb{C},x)^{\op{sch}}$. Then for
any $q \geq 1$ the weight tower of $F$ induces a filtration 
\[
\dots \subset F^{(p)}_{W}\pi_{q}(F,*)\subset
F^{(p-1)}_{W}\pi_{q}(F,*)\subset \dots\subset 
F^{(0)}_{W}\pi_{q}(F,*)=\pi_{q}(F,*)
\]
and a spectral sequence $\{E_{r}^{p,q+p}(W^{(*)}F^{0})\}$, so that
\[
E_{\infty}^{p,p+q}\simeq
F^{(p-1)}_{W}\pi_{q}(F,*)/F^{(p)}_{W}\pi_{q}(F,*), 
\]
and $E^{p,q}_{\bullet}$ degenerates at $E_{2}$ (see
Theorem~\ref{tpurity}).

 As a corollary we get the following
description of the homotopy groups of the schematization. Let $G$
denote the complex pro-reductive completion of $\pi_{1}(X,x)$ and let
$\mathcal{O}(G)$ be the algebra of regular funcions on $G$. By
definition $\pi_{1}(X,x)$ comes with a Zariski dense representation
into $G$, which can be combined with the left regular representation
of $G$ on $\mathcal{O}(G)$ to define a local system of algebras
on $X$ which by abuse of notation will be denoted again by
$\mathcal{O}(G)$. We will call 
${\mathcal{O}}(G)$ the {\em universal reductive local
  system} on $X$.  With this notation we have:

\begin{cori}(see Corollary \ref{c30})\label{ciweight} If $F
  :=(X^{\op{top}}\otimes\mathbb{C},x)^{\op{sch}}$ is the schematization of
  a complex projective manifold, then   
\begin{enumerate}
\item $\pi_{q}(F,*)\simeq \lim\, \pi_{q}(F,*)/F^{(p)}_{W}\pi_{q}(F,*)$.
\item The vector spaces
  $F^{(p-1)}_{W}\pi_{q}(F,*)/F^{(p)}_{W}\pi_{q}(F,*)$ only depend 
on the graded algebra \linebreak 
$H^{*}(X,{\mathcal{O}}(G))$, where
  $H^{*}(X,{\mathcal{O}}(G))$ is  the 
  cohomology algebra of $X$ with coefficients 
in the universal reductive local system.
\end{enumerate}
\end{cori}

\

\noindent
Corollary \ref{ciweight} provides strong eveidence that the
schematization of a smooth and projective complex manifold is much
more simple than the schematization of a generic topological
manifold. Note that even though the weight spectral sequence
$\{E_{r}^{p,q+p}(W^{(*)}F^{0})\}$ is purely topological and exists for
any schematic homotopy type $F$, the proof of the degeneration of this
sequence uses some weight properties of the action of
$\mathbb{C}^{\times \delta}$ induced by our Hodge decomposition.  It
is a striking fact that even if this action is \emph{not} an action of
the multiplicative group $\mathbb{G}_{m}$, the notion of weight
retains some of its algebraic character and ultimately forces the
spectral sequence to degenerate.

\smallskip

\begin{center} \textit{Restrictions on homotopy types} \end{center}

\smallskip

For any pointed connected homotopy type $X$, the schematization
$\sch{X}{\mathbb{C}}$ can be used to define a new homotopy invariant
of the space $X$, which captures information about the action of
$\pi_{1}(X)$ on the higher homotopy groups $\pi_{i}(X)$. More
precisely, we will define $\op{Supp}(\pi_{i}(\sch{X}{\mathbb{C}},x))$,
the \textit{support} of the sheaf $\pi_{i}(\sch{X}{\mathbb{C}},x))$,
as the subset of isomorphism classes of all simple representations of
$\pi_{1}(X,x)$ which appear in a finite dimensional sub-quotient of
$\pi_{i}(\sch{X}{\mathbb{C}},x))$ (see section \ref{dsupp}).

The existence of the Hodge decomposition described in theorem
\ref{ti1} imposes some restrictions on the
homotopy invariants $\op{Supp}(\pi_{i}(\sch{X}{\mathbb{C}},x))$. Our first
observation is the following corollary.

\begin{cori}(Corollary \ref{vhs}) \label{ci1}
Let $X$ be a pointed complex smooth projective algebraic variety, and
let $R(\pi_{1}(X,x))$ be the coarse moduli space of
simple finite dimensional representations of
$\pi_{1}(X,x)$.
\begin{enumerate}
\item The subset $\op{Supp}(\pi_{i}(\sch{X^{\op{top}}}{\mathbb{C}},x))
  \subset 
R(\pi_{1}(X,x))$ is stable under 
the action of $\mathbb{C}^{\times }$ on $R(\pi_{1}(X,x))$.
\item If $\rho \in Supp(\pi_{i}(\sch{X^{\op{top}}}{\mathbb{C}},x)$
is an isolated point
(for the topology induced from the analytic topology of
$R(\pi_{1}(X,x))$), then the local system on $X$ corresponding to
$\rho$ underlies
a polarizable complex variation of Hodge structures.
\item If $\pi_{i}((X^{\op{top}}\otimes \mathbb{C}),x)$ is an affine group
scheme of finite type, then, each
simple factor of the semi-simplification of the representation
of $\pi_{1}(X,x)$ on the vector space $\pi_{i}((X^{\op{top}}\otimes
\mathbb{C}),x)$ underlies a
polarizable complex variation of Hodge structures on $X$.
\item Suppose that $\pi_{1}(X,x)$ is abelian. Then, each isolated
character $\chi \in Supp(\pi_{i}((X^{\op{top}}\otimes
\mathbb{C})^{sch},x))$ is unitary.
\end{enumerate}
\end{cori}

The previous corollary suggests that one can study the action of the
fundamental group of a projective variety on its higher homotopy
groups by means of the support invariants.  However the reader should keep in
mind that the invariant $\op{Supp}(\pi_{i}(\sch{X}{\mathbb{C}},x))$ is
related to the 
action of $\pi_{1}(X,x)$ on $\pi_{i}(X,x)\otimes \mathbb{C}$ in a
highly non-trivial way, which at the moment can be understood only in some
very special cases. Nonetheless, one can use Corollary \ref{ci1} to
produce explicit families of new examples of homotopy types which are
not realizable by smooth and projective algebraic varieties (see
Theorem \ref{cex1}).

In the same vein we discuss two other applications. First, in Theorem
\ref{tfo}, we present a \emph{formality result}, asserting that the
schematization $\sch{X^{\op{top}}}{\mathbb{C}}$ of a smooth projective
complex variety $X$ is completely determined by the pro-reductive
fundamental group $\pi_{1}(X,x)^{\op{red}}$ and the cohomology algebra
$H^{*}(X,\mathcal{O}(\pi_{1}(X,x)^{\op{red}}))$. This theorem 
generalizes and extends the formality result of \cite{dgms}. Finally we
provide topological conditions on a smooth projective manifold $X$
under which the image of the Hurewitz morphism $\pi_{n}(X)
\longrightarrow H_{n}(X)$ is a sub-Hodge structure \ref{philippe}.

\

\smallskip

\begin{center} \textit{Organization of the paper} \end{center}

\smallskip

\

The paper is organized in four chapters. In the first one we recall
the definitions and main results concerning affine stacks and
schematic homotopy types from \cite{t1}. In particular, we recall
Theorem \ref{t3}, which shows how equivariant co-simplicial algebras
are related to the schematization functor. The proofs have not been
included, and can be found in \cite{t1,small}.

In the second chapter we construct the Hodge decomposition on
$\sch{X}{\mathbb{C}}$. For this, we will first review the non-abelian
Hodge correspondence between local systems and Higgs bundles of
\cite{s2}. We will then explain how to describe the schematization of
a space underlying a smooth manifold in terms of differential forms.
Finally after introducing the notion of \textit{fixed-point model
category} that will be used to define the Hodge decomposition, we
conclude the chapter with a proof of Theorem \ref{ti1}.

In the next chapter we define the weight tower of any
schematic homotopy type, and define the associated spectral
sequence in homotopy. Then, the Hodge decomposition
constructed in the previous chapter is used in order to
prove the degeneration statement in Corollary~\ref{ciweight}.

In the last chapter we show how the existence of the Hodge
decomposition imposes restrictions on homotopy types.  In particular,
we construct a whole family of examples of homotopy types which are
not realizable by smooth projective complex varieties. We also prove
the formality statement generalizing the formality theorem of
\cite{dgms} to the schematization, as well as an application of the
existence of the Hodge decomposition to the study of the Hurewitz
map. 

To keep the exposition more focused we concentrate here on the
Hodge theoretic aspects of our consturction. This necessitates
delegating several technical details (in particular several proofs)
concerning background material in schematic homotopy theory to the
companion but independent paper \cite{small}.

\

\smallskip

\

\begin{center} \textit{Related and future works} \end{center}

The construction used in order to define
our Hodge decomposition is similar to the one
used by R. Hain in \cite{ha2}. However, the two approaches differ
as the construction of \cite{ha2} is done relatively to
some variations of Hodge structures, whereas we are
taking into account all local systems. We do not think
that Corollary \ref{ci1} can be obtained by the techniques of \cite{ha2}, as
it uses in a non trivial way the $\mathbb{C}^{\times}$-equivariant
geometry of the whole moduli space of local systems. Note also that
the results of \cite{ha2} do not use at all the non-abelian Hodge
correspondence whereas our constructions highly depend on it.

Recently, L.Katzarkov, T.Pantev and C.Simpson, defined a notion of an
\textit{abstract non-abelian mixed Hodge structure} \cite{kps} and
discussed the existence of such structures on the non-abelian
cohomology of a smooth projective variety. The comparison between the
\cite{kps} approach and the construction of the present work seems
very difficult, essentially because both theories still need to be
developed before one can even state any conjectures comparing the two
points of view. In fact, it seems that a direct comparison of the two
approaches is not feasible within our current understanding of the
Hodge decomposition on the schematic homotopy type. Indeed, the object
$\sch{X^{\op{top}}}{\mathbb{C}}$ is insensitive to the topology on the
space of local systems on $X$, and therefore an important part of the
geometry (which is captured in the construction of \cite{kps}) is
lost.  This problem is reflected concretely in the fact that the
${\mathbb C}^{\times}$-action on $\sch{X^{\op{top}}}{\mathbb{C}}$
given by our Hodge decomposition is only an action of the discrete
group $\mathbb{C}^{\times\delta}$, and therefore does not give rise to
a \textit{filtration} in the sense of \cite{kps}.  The fact that the
object $\sch{X}{\mathbb{C}}$ does not vary well in families is another
consequence of the same problem. This is similar to the fact that the
pro-algebraic fundamental group \textit{does not know that local
systems can vary in algebraic famillies} (see \cite[$1.3$]{d3}). One
way to resolve this problem would be to consider the schematization of
a space $\sch{X^{\op{top}}}{\mathbb{C}}$ as the
\textit{$\mathbb{C}$-points of a schematic shape}, which is a bigger
object involving \textit{schematizations $\sch{X^{\op{top}}}{A}$ over
various $\mathbb{C}$-algebras $A$}. We believe that such an object
does exist and that it can be endowed with a mixed Hodge structure,
extending our construction on $\sch{X^{\op{top}}}{\mathbb{C}}$.

As explained before, our construction of the Hodge decomposition uses
equivariant co-simplicial algebras
as algebraic models for schematic homotopy types. These algebraic
models are very close to the
equivariant differential graded algebras of \cite{bs,ght}.
The main difference between the two approaches is that we use
co-simplicial algebras equipped
with an action of an affine group scheme, whereas in \cite{bs,ght} the
authors use algebras equivariant
for a discrete group action. In a sense, our approach is an
\textit{algebraization} of their approach, adapted for the purpose of
Hodge theory. Our formality theorem \ref{tfo} can also be considered
as a possible answer to \S 7 Problem $2$ of \cite{ght}.

Originally, a conjectural construction of the Hodge decomposition was
proposed in \cite{to2}, where
the notion  of a \textit{simplicial Tannakian category} was used. The
reader may notice the
Tannakian nature of the construction given in $\S 2.3$.

Finally, a crystalline version of
our main result has been recently worked out by M. Olsson in
\cite{ol}, who has deduced from it some new results
on homotopy types of algebraic varieties over fields of positive
characteristics (already on the level of the pro-nilpotent fundamental group).
M. Olsson is also working out a p-adic version of non-abelian Hodge theory,
based in the same way on the theory of schematic homotopy types.

\

\bigskip

\noindent
\textit{Acknowledgements:} We are very grateful to C.Simpson for
introducing us to the subject and for his constant
attention to this work. Special thanks are due to P.Deligne for his
letter \cite{del}
which inspired our construction of the
Hodge decomposition. We also thank A.Beilinson for his interest
and support. Finally, we thank M.Olsson who spotted
few gaps and mistakes in a previous version of this work.

\

\bigskip

\begin{center} \textit{Notations and conventions} \end{center}

\

\bigskip

Let $\Delta$ denote the standard simplicial category. Recall
that its objects
are the ordered finite sets $[n]:=\{0,\dots,n\}$. For any category $C$
and any functor $F : \Delta \longrightarrow C$
or $F : \Delta^{op} \longrightarrow C$, we write $F_{n}$ for the
object $F([n])$.

We fix an universe $\mathbb{U}$ with $\Delta \in \mathbb{U}$ and let
$\sff{Aff}/\mathbb{C}$ be the category of those affine schemes over
$\mathbb{C}$ which belong to $\mathbb{U}$. In this paper any affine
scheme is assumed to be an object in $\sff{Aff}/\mathbb{C}$ (i.e. always
belongs to $\mathbb{U}$), unless it is explicitly stated otherwise.
Throughout the paper we will tacitly identify
$\sff{Aff}/\mathbb{C}$  with the opposite
category of the category of $\mathbb{C}$-algebras belonging to
$\mathbb{U}$,  and the
objects of $\sff{Aff}/\mathbb{C}$ will be sometimes considered  as
algebras via this
equivalence.

We will use the Grothendieck site
$(\sff{Aff}/\mathbb{C})_{{\op{ffqc}}}$, whose underlying category is
$\sff{Aff}/\mathbb{C}$ equipped with the faithfully flat and
quasi-compact topology.

We fix a second universe $\mathbb{V}$ such that $\mathbb{U}\in
\mathbb{V}$, and let $\sff{SSet}$ be the category of simplicial sets
in $\mathbb{V}$.
We denote by  $\sff{SPr}(\mathbb{C})$ the category of $\sff{SSet}$-valued
presheaves on the site $(\sff{Aff}/\mathbb{C})_{{\op{ffqc}}}$. We will always
consider the category $\sff{SPr}(\mathbb{C})$ together
with its local projective model category structure described in
\cite{bl} (see also \cite[Definition $1.1.1$]{t1}), and the words
\textit{equivalence}, \textit{fibration} and \textit{cofibration} will
always refer to this model structure. The homotopy category of
$\sff{SPr}(\mathbb{C})$ will be denoted
by $\op{Ho}(\sff{SPr}(\mathbb{C}))$, and its objects will be called
\textit{stacks}. In the same way \textit{morphism of stacks}
will always refer
to a morphism in $\op{Ho}(\sff{SPr}(\mathbb{C}))$. When we need to consider
objects or morphisms in $\sff{SPr}(\mathbb{C})$ we will use the expressions
\textit{simplicial presheaves} and \textit{morphism of simplicial
presheaves} instead.

Any presheaf of sets on $(\sff{Aff}/\mathbb{C})_{{\op{ffqc}}}$ can be
considered as a presheaf of constant simplicial sets, and so as an
object in $\sff{SPr}(\mathbb{C})$. In particular, the Yoneda embedding
gives rise to a functor from $\sff{Aff}/\mathbb{C}$ to
$\sff{SPr}(\mathbb{C})$, which induces a fully faithful functor from
$\sff{Aff}/\mathbb{C}$ into $\op{Ho}(\sff{SPr}(\mathbb{C}))$. Via this
embedding, all affine schemes will be considered both as objects in
$\sff{SPr}(\mathbb{C})$ and as objects in
$\op{Ho}(\sff{SPr}(\mathbb{C}))$.

We denote by $\sff{SPr}_{*}(\mathbb{C})$ the model category of pointed
objects in $\sff{SPr}(\mathbb{C})$.  For any pointed simplicial
presheaf $F \in \sff{SPr}_{*}(\mathbb{C})$ we write $\pi_{i}(F,*)$ for
the homotopy sheaves of $F$ (see for example \cite[$1.1.1$]{t1}).
Objects in $\op{Ho}(\sff{SPr}_{*}(\mathbb{C}))$ will be called
\textit{pointed stacks}, and those with $\pi_{0}(F)=*$ \textit{pointed
connected stacks}.

Our references for model categories are \cite{hi,ho}, and we will
always suppose that model categories are $\mathbb{V}$-categories
(i.e. the set of morphism between two objects belongs to $\mathbb{V}$).

For any simplicial model category $M$, $\underline{\op{Hom}}(a,b)$ will be
the simplicial set of morphisms between $a$ and $b$
in $M$, and $\mathbb{R}\underline{\op{Hom}}(a,b)$ will be its derived version.
The objects $\mathbb{R}\underline{\op{Hom}}(a,b)$ are
well defined and functorial in the homotopy category $\op{Ho}(\sff{SSet})$.
The simplicial $Hom$ in the model category of pointed objects $M_{*}$
will be denoted by $\underline{\op{Hom}}_{*}$, and
its derived version by $\mathbb{R}\underline{\op{Hom}}_{*}$.

All complexes considered in this paper are
co-chain complexes (i.e. the differential increases degrees). We will
use freely (and often implicitly) the dual Dold-Kan correspondence
between positively graded complexes of vector spaces and co-simplicial
vector spaces (see \cite[$1.4$]{kr}).

Finally, for a group $\Gamma$ which belongs to the universe
$\mathbb{U}$, we will denote by $\Gamma^{\op{alg}}$ (respectively
$\Gamma^{\op{red}}$) its pro-algebraic completion (respectively its
pro-reductive completion).  By definition,
$\Gamma^{\op{alg}}$ is the universal affine group scheme (respectively
affine and reductive group scheme) which admits
a morphism (resp. a morphism with Zariski dense image)
from the constant sheaf of groups $\Gamma$.

\section{Review of the schematization functor}

In this first chapter, we review the theory of affine stacks and
schematic homotopy types introduced in \cite{t1}. The main goal
is to recall the theory and fix the notations and the
terminology. For further details and proofs the reader may
wish to consult \cite{t1,small}.

\subsection{Affine stacks and schematic homotopy types}

To begin with, let us recall some basic facts about co-simplicial
algebras. For us, an algebra will always
mean a commutative unital $\mathbb{C}$-algebra which belongs to
$\mathbb{V}$.

We will denote by $\sff{Alg}^{\Delta}$ the category
of co-simplicial algebras belonging to the universe $\mathbb{V}$. By
definition $\sff{Alg}^{\Delta}$ is the category of functors from
the standard simplicial category $\Delta$, to the category of
$\mathbb{C}$-algebras in $\mathbb{V}$.
For any co-simplicial algebra $A \in \sff{Alg}^{\Delta}$, one can
consider its underlying co-simplicial $\mathbb{C}$-vector
space, and associate to it its normalized co-chain complex $\No(A)$
(see \cite[$1.4$]{kr}). Any
morphism $f : A \longrightarrow B$ in $\sff{Alg}^{\Delta}$, induces a
natural morphism of complexes 
\[
\No(f) : \No(A) \longrightarrow \No(B).
\]
We will say that $f$ is an equivalence if $\No(f) : \No(A)
\longrightarrow 
\No(B)$ is a quasi-isomorphism.

The category
$\sff{Alg}^{\Delta}$ is endowed with a simplicial closed model
category structure for which the fibrations are epimorphisms,  the
equivalences are defined above, and the cofibrations are defined by
the usual lifting property
with respect to trivial fibrations. This model category is known to be
cofibrantly generated, and even finitely
generated (see \cite[$\S 2.1$]{ho}).

The category we are really interested in is the full subcategory of
$\sff{Alg}^{\Delta}$ of objects belonging to $\mathbb{U}$.  This
category is also a finitely generated simplicial closed model
category, which is a sub-model category of $\sff{Alg}^{\Delta}$ in a
very strong sense. For example, it satisfies the following stability
properties, which will be tacitly used in the rest of this work. Both 
properties can be deduced easily from the fact that the model
categories in question are finitely generated, and the fact that the
sets of generating cofibrations and trivial cofibrations (as well as
their domains and codomains) of $\sff{Alg}^{\Delta}$ all belong to
$\mathbb{U}$.

\begin{enumerate}
\item If $A$ is a co-simplicial algebra which is equivalent to a
co-simplicial algebra in $\mathbb{U}$, then
there exist a co-simplicial algebra in $\mathbb{U}$, which is
cofibrant in 
$\sff{Alg}^{\Delta}$, and which is 
equivalent to $A$.
\item Let $\op{Ho}(\sff{Alg}^{\Delta}_{\mathbb{U}})$ denote
the homotopy category of co-simplicial algebras belongings to
$\mathbb{U}$. Then
the natural functor $\op{Ho}(\sff{Alg}^{\Delta}_{\mathbb{U}})
\longrightarrow \op{Ho}(\sff{Alg}^{\Delta})$ is fully faithful. Its
essential image 
consists
of co-simplicial algebras which are isomorphic in
$\op{Ho}(\sff{Alg}^{\Delta})$ to co-simplicial
algebras in $\mathbb{U}$.
\end{enumerate}

\

\bigskip

Let $\sff{CDGA}$ denote the model category of commutative differential
(positively) graded $\mathbb{C}$-algebras belonging to
$\mathbb{V}$ (see \cite{bg}). Let
\[
\opi{Th} : \op{Ho}(\sff{Alg}^{\Delta}) \longrightarrow \op{Ho}(\sff{CDGA}),
\]
denote the functor of Thom-Sullivan cochains introduced in
\cite[Theorem $4.1$]{his}. In modern language, the functor $Th$ is
equivalent to the
functor of homotopy limits along the category $\Delta$ in the model category
$\sff{CDGA}$ (see \cite[$\S 20$]{hi}).
The functor $\opi{Th}$ has an inverse
\[
D : \op{Ho}(\sff{CDGA}) \longrightarrow
\op{Ho}(\sff{Alg}^{\Delta})
\]
which is defined as follows. Let
$\De : C^{+} \longrightarrow Vect^{\Delta}$ denote the {\em
denormalization functor} (described e.g. in
\cite[$1.4$]{kr}) from the category $C^{+}$ of co-chain complexes
of ${\mathbb C}$-vector spaces concentrated in positive degrees to
the category $\opi{Vect}^{\Delta}$ of co-simplicial vector
spaces. Recall from \cite[$1.4$]{kr} that $\De$ admits
functorial morphisms (the \textit{shuffle product})
$$\De(E)\otimes \De(F) \longrightarrow \De(E\otimes F)$$
which are unital associative and commutative. In particular, if
$A \in
\sff{CDGA}$ is a commutative
differential graded algebra, one can use the shuffle product to
define a natural commutative
algebra structure on $\De(A)$ by the
composition $\De(A)\otimes \De(A) \longrightarrow \De(A\otimes A)
\longrightarrow \De(A)$. This defines a functor
$$\De : \sff{CDGA} \longrightarrow \sff{Alg}^{\Delta},$$
which preserves equivalences, and thus induces a functor on the level
of homotopy categories
$$D : \op{Ho}(\sff{CDGA}) \longrightarrow \op{Ho}(\sff{Alg}^{\Delta}).$$
Now Theorem \cite[$4.1$]{his} implies that $\opi{Th}$ and $D$ are
inverse equivalences.

Since the functors $D$ and $\opi{Th}$ induce equivalences of homotopy
categories, readers who are not comfortable with co-simplicial objects
can replace co-simplicial algebras by commutative differential graded
algebras in the discussion below. However, in order to emphasize the
homotopy-theoretic context of our construction we choose to work
systematically with co-simplicial algebras rather than commutative
differential graded algebras.

Next we define the geometric {\em spectrum} of a co-simplicial algebra. More
precisely we define a functor 
\[
\op{Spec} : (\sff{Alg}^{\Delta})^{op} \longrightarrow
\sff{SPr}(\mathbb{C}),
\]
by the following formula
\[
\xymatrix@R=1pt@C=9pt{
\op{Spec} A \; : & (\sff{Aff}/\mathbb{C})^{op} \ar[r] &
\opi{\sff{SSet}} \\
&  \op{Spec} B   \ar@{|->}[r] & \underline{\op{Hom}}(A,B),
}
\]
where as usual $\underline{\op{Hom}}(A,B)$ denotes the simplicial set of
morphisms from the co-simplicial algebra $A$ to the co-simplicial algebra $B$.
In other words, if $A$ is given by a co-simplicial object $[n] \mapsto
A_{n}$, then the presheaf of $n$-simplices
of $\op{Spec} A$ is given by $(\op{Spec} B) \mapsto \op{Hom}(A_{n},B)$.

The functor $\op{Spec}$ is a right Quillen functor. Its left adjoint
functor  $\mathcal{O}$ associates to each simplicial  presheaf $F$ the
co-simplicial algebra of
\textit{functions} on $F$.  Explicitly
\[
\mathcal{O}(F)_{n}:=\op{Hom}(F,\mathbb{G}_{a}),
\]
where $\mathbb{G}_{a}:=\op{Spec} \mathbb{C}[T]\in \sff{Aff}/\mathbb{C}$
is the additive group scheme.

The right derived functor  of $\op{Spec}$ induces a functor on the
level of homotopy categorie
\[
\mathbb{R}\op{Spec} : \op{Ho}(\sff{Alg}^{\Delta})^{op} \longrightarrow
\op{Ho}(\sff{SPr}(\mathbb{C})),
\]
whose restriction to the full sub-category of $\op{Ho}(\sff{Alg}^{\Delta})$
consisting of objects isomorphic to a co-simplicial
algebra in $\mathbb{U}$ is fully faithful. Furthermore, for any
object $A \in \op{Ho}(\sff{Alg}^{\Delta})$, isomorphic to
some co-simplicial algebra belonging of $\mathbb{U}$,
the adjunction
morphism 
\[
A \longrightarrow \mathbb{L}\mathcal{O}(\mathbb{R}\op{Spec} A),
\]
in $\op{Ho}(\sff{Alg}^{\Delta})$, 
is an isomorphism. We are now ready to define affine stacks:

\begin{df}{(\cite[Definition $2.2.4$]{t1})}\label{d1}
An \emph{affine stack} is a stack $F \in \op{Ho}(\sff{SPr}(\mathbb{C}))$
isomorphic to $\mathbb{R}\op{Spec} A$, for some
co-simplicial algebra $A$ belonging to $\mathbb{U}$.
\end{df}

\

\medskip

\

The following important  result characterizes pointed connected
affine stacks, and relates it to homotopy theory over the
complex numbers.

\begin{thm}{(\cite[Theorem $2.4.1$, $2.4.5$]{t1})}\label{t0}
Let $F\in \op{Ho}(\sff{SPr}_{*}(\mathbb{C}))$ be a pointed stack. The
following three conditions are equivalent.
\begin{enumerate}
\item The pointed stack $F$ is affine and connected.
\item The pointed stack $F$ is connected and for all $i>0$ the sheaf
$\pi_{i}(F,*)$ is represented by
an affine unipotent group scheme (see \cite[IV \S 2]{dg}).
\item There exist a cohomologically connected co-simplicial algebra
$A$ (i.e. $H^{0}(A)\simeq \mathbb{C}$), which
belongs to $\mathbb{U}$, and such that
$F\simeq \mathbb{R}\op{Spec} A$.
\end{enumerate}
\end{thm}

\

\noindent
Recall that a commutative and unipotent affine group scheme
over $\mathbb{C}$ is the same thing as
a linearly compact vector space (see \cite[IV, \S 2 Proposition $4.2$
$(b)$]{dg} and \cite[II.1.4]{sa}). This implies that
for $F$ a pointed connected affine stack, and $i>1$,
the sheaves $\pi_{i}(F,*)$
are associated to well defined $\mathbb{C}$-vector spaces
$\pi^{i}(F,*)\in \mathbb{U}$, by the formula
\[
\xymatrix@R=1pt@C=9pt{
\pi_{i}(F,*) \; : & (\sff{Aff}/\mathbb{C})_{{\op{ffqc}}} \ar[r] & \sff{Ab} \\
& \op{Spec} B \ar@{|->}[r] &
\op{Hom}_{\mathbb{C}-\sff{Vect}}(\pi^{i}(F,*),B).
}
\]
If $F \simeq \mathbb{R}\op{Spec} A$, the vector spaces
$\pi^{i}(F,*)$ are isomorphic to the co-homotopy groups
$\pi^{i}(\opi{Th}(A))$ of the corresponding commutative differential
graded algebra as defined in \cite[$6.12$]{bg}.

Let us recall now that for a pointed simplicial presheaf $F$, one can
define its simplicial presheaf of loops
$\Omega_{*}F$. The functor $\Omega_{*} : \sff{SPr}_{*}(\mathbb{C})
\longrightarrow \sff{SPr}(\mathbb{C})$ is right Quillen,
and can be derived to a functor defined on the level of homotopy categories
$$\mathbb{R}\Omega_{*}F : \op{Ho}(\sff{SPr}_{*}(\mathbb{C})) \longrightarrow
\op{Ho}(\sff{SPr}(\mathbb{C})).$$

\begin{df}\label{d2}
A pointed and connected stack $F \in
\op{Ho}(\sff{SPr}_{*}(\mathbb{C}))$ is called a {\em pointed affine
$\infty$-gerbe} if the loop stack $\mathbb{R}\Omega_{*}F \in
\op{Ho}(\sff{SPr}(\mathbb{C}))$ is affine.
\end{df}

By definition, a pointed affine $\infty$-gerbe is a pointed stack, and
will always be considered in the category
of pointed stacks $\op{Ho}(\sff{SPr}_{*}(\mathbb{C}))$. In the same way, a
morphism between pointed affine $\infty$-gerbes
will always mean a morphism in $\op{Ho}(\sff{SPr}_{*}(\mathbb{C}))$.

A pointed schematic homotopy type will be a pointed affine
$\infty$-gerbe which in addition satisfies
a cohomological condition. Before we state this condition, let us recall that
for any algebraic group $G$, for any finite dimensional
linear representation $V$ of $G$ and any integer $n > 1$  one can
define a stack $K(G,V,n)$
(see \cite[$\S 1.2$]{t1}), which is the unique
pointed and connected stack having
\[
\pi_{1}(K(G,V,n),*)\simeq G \qquad \pi_{n}(K(G,V,n),*)\simeq V
\qquad \pi_{i}(K(G,V,n),*)\simeq 0 \; for \; i\neq 1,n,
\]
together with the given action of $G$ on $V$, and with trivial
associated Postnikov invariant in $H^{n+1}(G,V)$. With this notation we have:

\begin{df}
\begin{itemize}
\item A morphism of pointed stacks $f : F \longrightarrow F'$ is a
{\em $\boldsymbol{P}$-equivalence} if for any algebraic group $G$, any
linear representation of finite dimension $V$ and any integer $n>1$,
the induced morphism
\[
f^{*} : \mathbb{R}\underline{\op{Hom}}_{*}(F',K(G,V,n))
\longrightarrow \mathbb{R}\underline{\op{Hom}}_{*}(F,K(G,V,n))
\]
is an isomorphism.
\item A pointed stack $H$ is {\em $\boldsymbol{P}$-local}, if for any
$\boldsymbol{P}$-equivalence $f  : F \longrightarrow F'$ the induced
  morphism 
\[
f^{*} : \mathbb{R}\underline{\op{Hom}}_{*}(F',H) \longrightarrow
\mathbb{R}\underline{\op{Hom}}_{*}(F,H)
\]
is an isomorphism.
\item A {\em pointed schematic homotopy type} is a pointed affine
$\infty$-gerbe which is $\boldsymbol{P}$-local.
\end{itemize}
\end{df}

\

By definition, a pointed schematic homotopy type is a pointed stack,
and will always be considered in the category
of pointed stacks $\op{Ho}(\sff{SPr}_{*}(\mathbb{C}))$. In the same way, a
morphism between pointed schematic homotopy types
will always mean a morphism in $\op{Ho}(\sff{SPr}_{*}(\mathbb{C}))$.

The following theorem is  a partial analogue for
pointed schematic homotopy types of Theorem \ref{t0}.

\begin{thm}{(\cite[Theorem $3.2.4$]{t1}, \cite[Cor. 3.6]{small})}\label{t1}
Let $F$ be a pointed and connected stack. Then, $F$ is a pointed schematic
homotopy type if and only if it satisfies the following
two conditions.
\begin{enumerate}
\item The sheaf $\pi_{1}(F,*)$ is represented by an affine group scheme.
\item For any $i>1$, the sheaf $\pi_{i}(F,*)$ is represented by a
unipotent group scheme.
In other words, the sheaf $\pi_{i}(F,*)$ is represented by a linearly
compact vector space in
$\mathbb{U}$ (see \cite[II.1.4]{sa}).
\end{enumerate}
\end{thm}

\

\bigskip

We finish this section by recalling the main existence theorem of
schematic homotopy theory. The category $\opi{\sff{SSet}}$ can be
embedded into the category $\sff{SPr}(\mathbb{C})$ by vieweing a
simplicial set $X$ as a constant simplicial presheaf on
$(\sff{Aff}/\mathbb{C})_{{\op{ffqc}}}$.  With this convention we have
the following important definition:

\begin{df}{(\cite[Definition $3.3.1$]{t1})}\label{d6}
Let $X$ be a pointed and connected simplicial set in $\mathbb{U}$. A
\emph{schematization of $X$ over $\mathbb{C}$} is
a pointed schematic homotopy type $\sch{X}{\mathbb{C}}$,
together with a morphism in $\op{Ho}(\sff{SPr}_{*}(\mathbb{C}))$
\[
u : X \longrightarrow \sch{X}{\mathbb{C}}
\]
which is a universal for morphisms from $X$ to pointed schematic
homotopy types.
\end{df}

\

\smallskip

\noindent
\quad We have stated  Definition \ref{d6} only for
simplicial sets in order the simplify the exposition. However, by
using the {\em singular functor} $\opi{Sing}$ (see for example
\cite{ho}), from the
category of topological spaces to the category of simplicial sets, one
can define  the schematization of a pointed connected
topological space.  In what follows we will always assume implicitly
that the functor $\opi{Sing}$ has been applied when necessary  and we
will generally not distinguish between topological
spaces and simplicial sets when considering the schematization functor.

\

\begin{thm}{(\cite[Theorem $3.3.4$]{t1})}\label{t2}
Any pointed and connected simplicial set $(X,x)$ in $\mathbb{U}$
possesses a schematization over $\mathbb{C}$.
\end{thm}

\

As already mentioned, specifying a commutative and
unipotent affine group scheme over $\mathbb{C}$ is  equivalent to
specifying
a linearly compact vector space. This implies that the
sheaves \linebreak $\pi_{i}(\sch{X}{\mathbb{C}},x)$
are associated to well defined $\mathbb{C}$-vector space
$\pi^{i}(\sch{X}{\mathbb{C}},x)\in \mathbb{U}$, by the formula
\[
\xymatrix@R=1pt@C=9pt{
\pi_{i}(\sch{X}{\mathbb{C}},x) \; : &
(\sff{Aff}/\mathbb{C})_{{\op{ffqc}}} \ar[r] & \opi{Ab} \\
& \op{Spec} B \ar@{|->}[r]
&\op{Hom}_{\mathbb{C}-\sff{Vect}}(\pi^{i}((X\otimes
\mathbb{C})^{sch},x),B).
}
\]
Furthermore, the natural action of the affine group scheme
$\pi_{1}(\sch{X}{\mathbb{C}},x)$ on
$\pi_{i}(\sch{X}{\mathbb{C}},x)$ is continuous in the sense
that it is induced by an algebraic action
on the vector space \linebreak $\pi^{i}(\sch{X}{\mathbb{C}},x)$.

In general the homotopy sheaves $\pi_{i}((X\otimes
\mathbb{C})^{sch},x)$ are relatively big (they are not of finite type over
$\mathbb{C}$, even when
$X$ is a finite homotopy type) and are hard  to compute. The only two
general cases where one knows something
are the following.

\begin{prop}{(\cite[Corollaries $3.3.7$, $3.3.8$,
$3.3.9$]{t1})}\label{p2}
Let $X$ be a pointed connected simplicial set in $\mathbb{U}$, and let
$\sch{X}{\mathbb{C}}$ be its schematization.
\begin{enumerate}

\item The affine group scheme $\pi_{1}(\sch{X}{\mathbb{C}},x)$
is naturally isomorphic to the pro-algebraic completion
of the discrete group $\pi_{1}(X,x)$ over $\mathbb{C}$.

\item There is a natural isomorphism
\[
H^{\bullet}(X,\mathbb{C})\simeq
H^{\bullet}(\sch{X}{\mathbb{C}},\mathbb{G}_{a}). 
\]
\item If $X$ is simply connected and of finite type (i.e. the homotopy
type of a simply connected and
finite $CW$ complex), then for any $i>1$, the group scheme
$\pi_{i}(\sch{X}{\mathbb{C}},x)$ is naturally isomorphic to
the pro-unipotent
completion of the discrete groups $\pi_{i}(X,x)$. In other words, for
any $i>1$  
\[
\pi_{i}(\sch{X}{\mathbb{C}},x)\simeq
\pi_{i}(X,x)\otimes_{\mathbb{Z}}\mathbb{G}_{a}.
\]
\end{enumerate}
\end{prop}

\subsection{Equivariant stacks} \label{ss:equivariant}

For the duration of this section we fix a presheaf of groups $G$ on
$(\sff{Aff}/\mathbb{C})_{{\op{ffqc}}}$, which will be considered as a
group object 
in $\sff{SPr}(\mathbb{C})$. We will make the assumption that $G$ is
cofibrant as an object of $\sff{SPr}(\mathbb{C})$.
For example, $G$ could be representable (i.e. an affine group scheme),
or a constant presheaf associated to
a group in $\mathbb{U}$.

Let $\opi{G-\sff{SPr}}(\mathbb{C})$ be the category of
simplicial presheaves equipped with a left action of $G$, which is
again a closed model category
(see \cite{ss}). Recall that the fibrations
(respectively equivalences) in $\opi{G-\sff{SPr}}(\mathbb{C})$ are defined
to be the morphisms inducing fibrations (respectively equivalences)
between the underlying
simplicial presheaves. The model category $\opi{G-\sff{SPr}}(\mathbb{C})$
will be called the {\em model category of $G$-equivariant simplicial
presheaves}, and the objects in $\op{Ho}(\opi{G-\sff{SPr}}(\mathbb{C}))$
will be called {\em $G$-equivariant stacks}. For any $G$-equivariant
stacks $F$ and $F'$ we will denote by
$\underline{\op{Hom}}_{G}(F,F')$ the simplicial set of morphisms in
$\opi{G-\sff{SPr}}(\mathbb{C})$, and  by
$\mathbb{R}\underline{\op{Hom}}_{G}(F,F')$ its derived version.

Next recall that to any group $G$ one can associate  its classifying
simplicial presheaf $\opi{BG} \in \sff{SPr}_{*}(\mathbb{C})$
(see \cite[$\S 1.3$]{t1}). The
object $\opi{BG} \in \op{Ho}(\sff{SPr}_{*}(\mathbb{C}))$ is well defined up to
a unique isomorphism by the following properties
$$\pi_{0}(\opi{BG})\simeq * \qquad \pi_{1}(\opi{BG},*)\simeq G \qquad
\pi_{i}(\opi{BG},*)=0 \; for \; i>1.$$
Consider the coma category $\sff{SPr}(\mathbb{C})/\opi{BG}$, of
objects over the classifying simplicial presheaf $\opi{BG}$, endowed
with its natural simplicial closed model structure (see 
\cite{ho}). Recall that fibrations, equivalences
and cofibrations in $\sff{SPr}(\mathbb{C})/\opi{BG}$ are defined on
the underlying objects in $\sff{SPr}(\mathbb{C})$. We set
$\opi{BG}:=\opi{EG}/G$, where $\opi{EG}$ is a cofibrant model  (fixed
once and for all) of $*$ in
$\opi{G-\sff{SPr}}(\mathbb{C})$.

Next we define a pair of adjoint functors
\[
\xymatrix{
\opi{G-\sff{SPr}}(\mathbb{C}) \ar@<1ex>[r]^{\opi{De}} &
\sff{SPr}(\mathbb{C})/\opi{BG},
\ar@<1ex>[l]^{\opi{Mo}}
}
\]
where $\opi{De}$ stands for \textit{descent} and $\opi{Mo}$ for
\textit{monodromy}.
If $F$ is a $G$-equivariant simplicial presheaf, then
$\opi{De}(F)$ is defined to be $(\opi{EG}\times F)/G$, where $G$ acts
diagonally on $\opi{EG}\times F$.
Note that there is a natural projection
$\opi{De}(F) \longrightarrow \opi{EG}/G=\opi{BG}$, and so
$\opi{De}(F)$ is naturally
an object in $\sff{SPr}(\mathbb{C})/\opi{BG}$.

The functor $\opi{Mo}$ which is right adjoint to $\opi{De}$ can be
defined in the following way. For an object $F
\longrightarrow \opi{BG}$
in $\sff{SPr}(\mathbb{C})/\opi{BG}$, the simplicial presheaf
underlying $\opi{Mo}(F)$ is defined by
\[
\xymatrix@R=1pt@C=9pt{
\opi{Mo}(F) \; : & (\sff{Aff}/\mathbb{C})^{op} \ar[r] & \sff{SSet} \\
& Y \ar@{|->}[r] & \underline{\op{Hom}}_{\opi{BG}}(\opi{EG}\times Y,F),
}
\]
where $\underline{\op{Hom}}_{\opi{BG}}$ denotes the simplicial set of
morphisms in the
coma category $\sff{SPr}(\mathbb{C})/\opi{BG}$, and $\opi{EG}$ is
considered as an object in $\sff{SPr}(\mathbb{C})/\opi{BG}$ by
the natural projection $\opi{EG}\longrightarrow
\opi{EG}/G=\opi{BG}$. The action of $G$ on $\opi{Mo}(F)$ is then defined by
making $G$ act on $\opi{EG}$. We have the following useful

\begin{lem}{(\cite[Lemma 3.10]{small})}\label{l1}
The Quillen adjunction $(\opi{De},\opi{Mo})$ is a Quillen equivalence.
\end{lem}

\medskip

\noindent
The previous lemma implies that the derived Quillen adjunction induces
an equivalence of categories
\[
\op{Ho}(\opi{G-\sff{SPr}}(\mathbb{C})) \simeq
\op{Ho}(\sff{SPr}(\mathbb{C})/\opi{BG}). 
\]

\begin{df}
For any $G$-equivariant stack $F \in \op{Ho}(\opi{G-\sff{SPr}}(\mathbb{C}))$
define the {\em quotient stack $[F/G]$} of $F$ by $G$ as the
object $\mathbb{L}\opi{De}(F) \in
\op{Ho}(\sff{SPr}(\mathbb{C})/\opi{BG})$ corresponding to $F \in
\op{Ho}(\opi{G-\sff{SPr}}(\mathbb{C}))$.
\end{df}

The construction also implies that the homotopy fiber of the natural
projection
$$
p : [F/G] \longrightarrow \opi{BG}
$$
is naturally isomorphic to the
underlying stack of the $G$-equivariant stack $F$.

An important example of a quotient stack to keep in mind is the
following. Suppose that $G$ acts on a sheaf of groups
$V$. Then, $G$ acts also naturally on the simplicial presheaf
$K(V,n)$. The quotient $[K(V,n)/G]$ of $K(V,n)$ by
$G$ is naturally isomorphic to $K(G,V,n)$ described in the previous
section.

\subsection{Equivariant co-simplicial algebras and equivariant affine
  stacks} 

Suppose that $G$ is an affine group scheme and consider the category
of linear representations of $G$. By definition this is the category
of quasi-coherent sheaves of $\mathcal{O}$-modules in $\mathbb{V}$, on
the big site $(\sff{Aff}/\mathbb{C})_{\op{ffqc}}$, which are equipped
with a linear action of the presheaf of groups $G$.  Equivalently, it
is the category of co-modules in $\mathbb{V}$ over the co-algebra 
$\mathcal{O}(G)$ of regular functions on $G$.  This category will be denoted by
$\sff{Rep}(G)$. Note that it is an abelian $\mathbb{C}$-linear tensor
category, which admits all $\mathbb{V}$-limits and
$\mathbb{V}$-colimits.  The category of co-simplicial $G$-modules is
defined to be the category $\sff{Rep}(G)^{\Delta}$, of co-simplicial
objects in $\sff{Rep}(G)$.

Recall from \cite[\S 3.2]{small} that there
exists a simplicial finitely generated closed model structure on
the category $\sff{Rep}(G)^{\Delta}$, such that
the following properties are satisfied
\begin{itemize}
\item A morphism $f : E \longrightarrow E'$ is an equivalence if and
  only if, 
for any $i$, the induced morphism $H^{i}(f) : H^{i}(E) \longrightarrow
H^{i}(E')$ is an isomorphism.
\item A morphism $f : E \longrightarrow E'$ is a cofibration if and
only if, for any $n>0$,
the induced morphism $f_{n} : E_{n} \longrightarrow E'_{n}$ is a monomorphism.
\item A morphism $f : E \longrightarrow E'$ is a fibration if and only
if it is an epimorphism whose
kernel $K$ is such that for any $n\geq 0$, $K_{n}$ is an injective
object in $\sff{Rep}(G)$.
\end{itemize}

\bigskip

The category $\sff{Rep}(G)$ is endowed with a symmetric monoidal
structure, given by the tensor product of co-simplicial
$G$-modules (defined levelwise). In particular we can consider the
category $\opi{G-\sff{Alg}}^{\Delta}$ of commutative unital monoids in
$\sff{Rep}(G)^{\Delta}$. It is reasonable to view the objects in
$\opi{G-\sff{Alg}}^{\Delta}$ as co-simplicial
algebras equipped with an action of the group scheme $G$. Motivated by
this remark we will refer to the category $\opi{G-\sff{Alg}}^{\Delta}$ as
the {\em category of $G$-equivariant co-simplicial algebras}. From
another point of view, the category $\opi{G-\sff{Alg}}^{\Delta}$
is also the category of simplicial affine schemes of $\mathbb{V}$
equipped with an action of $G$.

Every $G$-equivariant co-simplicial algebra $A$ has an underlying
co-simplicial $G$-module again denoted by $A \in
\sff{Rep}(G)^{\Delta}$. This defines a forgetful functor
$$\opi{G-\sff{Alg}}^{\Delta} \longrightarrow \sff{Rep}(G)^{\Delta}$$
which has a left adjoint $L$, given by the free commutative
monoid construction. 

As proved in \cite[\S 3.2]{small}, there
exists a simplicial cofibrantly generated closed model structure
on the category $\opi{G-\sff{Alg}}^{\Delta}$, such that
the following two properties are satisfied
\begin{itemize}
\item A morphism $f : A \longrightarrow A'$ is an equivalence if and
only if the induced morphism in $\sff{Rep}(G)^{\Delta}$
is an equivalence.
\item A morphism $f : A \longrightarrow A'$ is a fibration if and only
if the induced morphism in
$\sff{Rep}(G)^{\Delta}$ is a fibration.
\end{itemize}

\

\bigskip

\

\noindent
Similarly to the non-equivariant case, we could  have
defined a model structure of $G$-equivariant
commutative differential graded algebras. This is the category
$\opi{G-\sff{CDGA}}$, of commutative monoids in $C^{+}(\sff{Rep}(G))$ - the
symmetric monoidal model category of positively graded co-chain complexes in
$\sff{Rep}(G)$.

There exists a cofibrantly generated closed model structure
on the category $\opi{G-\sff{CDGA}}$, such that
the following two properties are satisfied
\begin{itemize}
\item A morphism $f : A \longrightarrow A'$ in $\opi{G-\sff{CDGA}}$
is an equivalence if and
only if the induced morphism in $C^{+}(\sff{Rep}(G))$
is a quasi-isomorphism.
\item A morphism $f : A \longrightarrow A'$ in $\opi{G-\sff{CDGA}}$
is a fibration if and only
if the induced morphism in
$C^{+}(\sff{Rep}(G))$ is a fibration (i.e. $A_{n} \longrightarrow
A_{n}'$ is surjective
for any $n\geq 0$).
\end{itemize}

\

\bigskip

\noindent
As in the non-equivariant case, there exist a denormalization functor
\begin{equation} \label{eq:equivariantD}
D : \op{Ho}(\opi{G-\sff{CDGA}}) \longrightarrow
\op{Ho}(\opi{G-\sff{Alg}}^{\Delta}).
\end{equation}
Indeed, the denormalization functor and the shuffle products exist
over any $\mathbb{C}$-linear tensor
abelian base category, and in particular over the tensor category
$\sff{Rep}(G)$. Therefore one can repeat the definition of the
functor $\opi{Th}$ in \cite[$4.1$]{his} and produce a functor of
$G$-equivariant Thom-Sullivan co-chains:
\begin{equation} \label{eq:equivariantTh}
\opi{Th} : \op{Ho}(\opi{G-\sff{Alg}}^{\Delta}) \longrightarrow
\op{Ho}(\opi{G-\sff{CDGA}})
\end{equation}
which is an inverse of $D$. This allows us to view any $G$-equivariant
commutative differential graded algebra as a well defined object
in $\op{Ho}(\opi{G-\sff{Alg}}^{\Delta})$ and vise-versa. We will make
a frequent 
use of this point of view in what follows.

\

\bigskip

For any $G$-equivariant co-simplicial algebra $A$, one can define its
(geometric) spectrum \linebreak
$\op{Spec}_{G}\, A \in \opi{G-\sff{SPr}}(\mathbb{C})$, by
taking the usual spectrum of its underlying co-simplicial algebra and
keeping track of the $G$-action. Explicitly,
if $A$ is given by a morphism of co-simplicial algebras $A
\longrightarrow A\otimes \mathcal{O}(G)$, one finds
a morphism of simplicial schemes
$$G\times \op{Spec}\, A\simeq \op{Spec}\, (A\otimes \mathcal{O}(G))
\longrightarrow \op{Spec}\, A,$$
which induces a well defined $G$-action on the simplicial scheme
$\op{Spec}\, A$. Hence, by passing to the simplicial presheaves
represented by $G$ and $\op{Spec}\, A$, one gets the $G$-equivariant
simplicial presheaf $\op{Spec}_{G}(A)$.

This procedure defines
a functor
$$\op{Spec}_{G}\, : (\opi{G-\sff{Alg}}^{\Delta})^{op} \longrightarrow
\opi{G-\sff{SPr}}(\mathbb{C}),$$
which by \cite[\S 3.2]{small} is a right Quillen functor.
The left adjoint of $\op{Spec}_{G}$ will be denoted by
$\mathcal{O}_{G} : \opi{G-\sff{SPr}}(\mathbb{C}) \longrightarrow
\opi{G-\sff{Alg}}^{\Delta}$. \\

One can form the right derived functor of
$\op{Spec}_{G}$:
$$\mathbb{R}\op{Spec}_{G} : \op{Ho}(\opi{G-\sff{Alg}}^{\Delta})^{op}
\longrightarrow \op{Ho}(\opi{G-\sff{SPr}}(\mathbb{C})),$$
which possesses a left adjoint $\mathbb{L}\mathcal{O}_{G}$.
One can then compose this functor with the quotient stack functor
$[-/G]$, and obtain a functor
$$[\mathbb{R}\op{Spec}_{G}(-)/G] : \op{Ho}(\opi{G-\sff{Alg}}^{\Delta})^{op}
\longrightarrow \op{Ho}(\sff{SPr}(\mathbb{C})/\opi{BG}),$$
which still possesses a left adjoint due to the fact that
$[-/G]$ is an equivalence
of categories. We will denote this left adjoint again by
$$\mathbb{L}\mathcal{O}_{G} : \op{Ho}(\sff{SPr}(\mathbb{C})/\opi{BG})
\longrightarrow \op{Ho}(\opi{G-\sff{Alg}}^{\Delta})^{op}.$$

\

As proved in \cite[Proposition~3.14]{small},
If $A \in \op{Ho}(\opi{G-\sff{Alg}}^{\Delta})^{op}$ is isomorphic to some
$G$-equivariant co-simplicial algebra in
$\mathbb{U}$, then the adjunction morphism
$$A \longrightarrow
\mathbb{L}\mathcal{O}_{G}(\mathbb{R}\op{Spec}_{G}\, A)$$
is an isomorphism.
In particular, the functors $\mathbb{R}\op{Spec}_{G}$ and
$[\mathbb{R}\op{Spec}_{G}(-)/G]$ become fully faithful when
restricted to the full sub-category of
$\op{Ho}(\opi{G-\sff{Alg}}^{\Delta})$ consisting of $G$-equivariant co-simplicial
algebras isomorphic to some object in $\mathbb{U}$. This justifies the
following definition.

\begin{df}{\cite[\S 3.2]{small}}\label{d7}
An equivariant stack $F \in \op{Ho}(\opi{G-\sff{SPr}}(\mathbb{C}))$ is a
\emph{$G$-equivariant affine stack} if it is isomorphic to some
$\mathbb{R}\op{Spec}_{G}(A)$, with $A$ being a $G$-equivariant
co-simplicial algebra in $\mathbb{U}$.
\end{df}

\bigskip

We conclude this section with a proposition showing
that stacks of the form $[\mathbb{R}\op{Spec}_{G}\, (A)/G]$
are often pointed schematic homotopy types. Therefore the theory of
equivariant differential graded algebras gives a way to produce
schematic homotopy types.

\begin{prop}{\cite[\S 3.2]{small}}\label{p5}
Let $A \in \opi{G-\sff{Alg}}^{\Delta}$ be a $G$-equivariant co-simpicial
algebra in $\mathbb{U}$, such that the underlying algebra
of $A$ has an  augmentation $x : A \longrightarrow \mathbb{C}$ and is
connected (i.e.  $H^{0}(A)\simeq \mathbb{C}$).
Then, the
quotient stack $[\mathbb{R}\op{Spec}_{G}\, (A)/G]$ is a pointed
schematic homotopy type. Furthermore, one has
$$
\pi_{i}([\mathbb{R}\op{Spec}_{G}\, (A)/G],x)\simeq
\pi_{i}(\mathbb{R}\op{Spec}\, A,x) \text{ for }i>1,$$
and the fundamental group $\pi_{1}([\mathbb{R}\op{Spec}_{G}\, (A)/G],x)$
is an extension of $G$ by the pro-unipotent group
$\pi_{1}(\mathbb{R}Spec\, A,x)$. 
\end{prop}

\subsection{An explicit model for $\sch{X}{\mathbb{C}}$}
\label{s-explicit}

Let $X$ be a pointed and connected simplicial set in $\mathbb{U}$. In
this section we describr an explicit
model for $\sch{X}{\mathbb{C}}$ which is based on the notion of
equivariant affine stacks.

Let $G$ be the complex pro-reductive completion of the discrete group
$\pi_{1}(X,x)$. By definition $G$ comes with a universal homomorphism 
$\pi_{1}(X,x) \longrightarrow G$ with a Zariski dense image. The
universal homomorphism 
induces a morphism  $X \longrightarrow B(G(\mathbb{C}))$ be the
of simplicial sets.  This latter morphism is
well defined up to homotopy, and we choose a representative once and
for all.  Let $p : P \longrightarrow X$ be the corresponding
$G$-torsor in $\sff{SPr}(\mathbb{C})$. More precisely, $P$ is the
simplicial presheaf sending an affine scheme $\op{Spec} A \in
\sff{Aff}/\mathbb{C}$ to the simplicial set
$P(A):=(\opi{EG}(A)\times_{\opi{BG}(A)} X)$.  The morphism \linebreak 
$p : P \longrightarrow X$ is then a well defined morphism in
$\op{Ho}(\opi{G-\sff{SPr}}(\mathbb{C}))$, the group $G$ acting on
$P=(\opi{EG}\times_{\opi{BG}}X)$ by its action on $\opi{EG}$, and
trivially on $X$.  Alternatively we can describe $P$ by the formula
\[
P\simeq (\widetilde{X}\times G)/\pi_{1}(X,x),
\]
where $\widetilde{X}$ is the universal covering of $X$, and
$\pi_{1}(X,x)$ acts on $\widetilde{X}\times G$ by the diagonal
action (our convention here is that $\pi_{1}(X,x)$ acts on $G$ by left
translation). We assume a this point that $\widetilde{X}$ is chosen to
be  cofibrant in the model
category of $\pi_{1}(X,x)$-equivariant simplicial sets. For example,
we may assume that $\widetilde{X}$ is
a $\pi_{1}(X,x)$-equivariant cell complex.

We consider now the $G$-equivariant affine stack
$\mathbb{R}\op{Spec}_{G}\, \mathcal{O}_{G}(P) \in
\op{Ho}(\opi{G-\sff{SPr}}(\mathbb{C}))$,
which comes naturally equipped with its adjunction morphism $P
\longrightarrow \mathbb{R}\op{Spec}_{G}\,
\mathbb{L}\mathcal{O}_{G}(P)$.
This induces a well defined morphism in $\op{Ho}(\sff{SPr}(\mathbb{C}))$:
\[
X \simeq [P/G] \longrightarrow [\mathbb{R}\op{Spec}_{G}\,
\mathbb{L}\mathcal{O}_{G}(P)/G].
\]
Furthermore, as $X$ is pointed, this morphism induces a natural morphism in
$\op{Ho}(\sff{SPr}_{*}(\mathbb{C}))$
\[
u : X \longrightarrow [\mathbb{R}\op{Spec}_{G}\,
\mathbb{L}\mathcal{O}_{G}(P)/G].
\]
With this notation we now have the following important

\begin{thm}{(\cite[Theorem 3.20]{small})}\label{t3}
The natural morphism
\[
u : X \longrightarrow [\mathbb{R}\op{Spec}_{G}\,
\mathbb{L}\mathcal{O}_{G}(P)/G]
\]
is a model for the schematization of $X$.
\end{thm}

\

\bigskip

Let $(X,x)$ be a pointed connected simplicial set in $\mathbb{U}$,
let $G$ be the pro-reductive completion
of the group $\pi_{1}(X,x)$, and let $\widetilde{X}$ be the universal
covering of $X$. Again, we assume that $\widetilde{X}$ is
chosen to be cofibrant
as a $\pi_{1}(X,x)$-simplicial set.
Consider $\mathcal{O}(G)$ as a locally constant sheaf of algebras
on $X$ via the natural action of $\pi_{1}(X,x)$. Let
\[
C^{\bullet}(X,\mathcal{O}(G)):=
(\mathcal{O}(G)^{\widetilde{X}})^{\pi_{1}(X,x)}
\]
be the co-simplicial algebra of co-chains on $X$ with coefficients
in $\mathcal{O}(G)$ (see \cite[Section~3.2]{small} for
details). This 
co-simplicial algebra is equipped with a natural
$G$-action, induced by the
regular representation of $G$. One can thus consider
$C^{\bullet}(X,\mathcal{O}(G))$ as an object
in $\opi{G-\sff{Alg}}^{\Delta}$.  We have the following immediate

\begin{cor}{(\cite[Corollary 3.21]{small})}\label{c1}
With the previous notations, one has
$$\sch{X}{\mathbb{C}}\simeq [\mathbb{R}\op{Spec}_{G}
C^{\bullet}(X,\mathcal{O}(G))/G].$$
\end{cor}

\

\bigskip

\noindent
\textit{Remarks:} The previous corollary shows that the
schematic homotopy type $\sch{X}{\mathbb{C}}$
can be explicitly described by the $G$-equivariant co-simplicial
algebra of co-chains on $X$ with coefficients
in $\mathcal{O}(G)$. For example, it is possible to describe the
$\mathbb{C}$-vector spaces
$\pi_{i}(\sch{X}{\mathbb{C}},x)(\mathbb{C})$ for $i>1$ by
using the minimal model for the corresponding commutative differential
graded algebra
$\opi{Th}(C^{\bullet}(X,\mathcal{O}(G)))$. This was the original description
given by
P. Deligne in his letter \cite{del}.

The previous corollary can be restated in terms of the pro-algebraic
completion of $\pi_{1}(X,x)$ rather than the pro-reductive one. In
fact if in the statement in Corollary \ref{c1} we take $G$ to be
$\pi_{1}(X,x)^{\op{alg}}$, the statement remains valid and the proof
is exactly the same. In summary:  if $G$ is the pro-algebraic
completion of $\pi_{1}(X,x)$, and $C^{\bullet}(X,\mathcal{O}(G))$ is the
$G$-equivariant co-simplicial algebra
of co-chains  on $X$ with coefficients in $G$, one again has
$\sch{X}{\mathbb{C}}\simeq [\mathbb{R}\op{Spec}_{G}
C^{\bullet}(X,\mathcal{O}(G))/G]$.

\section{The Hodge decomposition}

In this section, we will construct the Hodge decomposition on
$\sch{X^{\op{top}}}{\mathbb{C}}$, when $X^{\op{top}}$
is the underlying homotopy type of a smooth projective complex
algebraic variety. This Hodge decomposition
is a higher analogue of the \textit{Hodge filtration} on the
pro-algebraic fundamental group defined by C.Simpson
in \cite{s2}. More precisely, it is an action of the discrete group
$\mathbb{C}^{\times \delta}$  on the pointed stack
$\sch{X^{\op{top}}}{\mathbb{C}}$, such that the induced actions on its
fundamental group, cohomology and homotopy groups
coincide with the various previously defined Hodge filtrations of
\cite{mo,ha,s2}.

The construction we propose here is based on two results. The first is
the explicit description (reviewed in section \ref{s-explicit}) of
$\sch{X^{\op{top}}}{\mathbb{C}}$ in terms of the
equivariant co-simplicial algebra of co-chains on $X^{\op{top}}$
 with coefficients in the universal reductive local system. The second
is the non-abelian Hodge theorem of \cite{s2} establishing a
correspondence between local systems and
Higgs bundles as well as their cohomology.

Here is a short outline of the construction. Let $G$ be the
pro-reductive completion of $\pi_{1}(X^{\op{top}},x)$, and let
$C^{\bullet}(X^{\op{top}},\mathcal{O}(G))$ be the $G$-equivariant
co-simplicial algebra of co-chains on $X^{\op{top}}$ with coefficients
in the local system of algebras $\mathcal{O}(G)$ (see
section~\ref{s-explicit}). First we use the non-abelian Hodge
correspondence to show that this $G$-equivariant co-simplicial algebra
can be constructed from the Dolbeault cohomology complexes of certain
Higgs bundles.  Since the category of Higgs bundles is equipped with a
natural $\mathbb{C}^{\times }$-action, this will induce a well defined
action of $\mathbb{C}^{\times \delta}$ on the pointed stack
$[\mathbb{R}\op{Spec}_{G} C^{\bullet}(X^{\op{top}},\mathcal{O}(G))/G]$
which, as we have seen, is isomorphic to
$\sch{X^{\op{top}}}{\mathbb{C}}$.

In order to implement this program we first recall the non-abelian
Hodge correspondence between local systems and Higgs bundles in the
form presented by C.Simpson in \cite{s2}. In particular we explain
briefly how this correspondence extends to the relevant
categories of $\opi{Ind}$-objects. Next we make precise the relation
between the explicit model for the schematization presented in
Corollary \ref{c1} and certain algebras of differential forms. We also
define the notion of a fixed-point model category, which generalizes
the notion of a model category of objects together with a group
action.  Finally we use these results to endow
$\sch{X^{\op{top}}}{\mathbb{C}}$ with an action of the group
$\mathbb{C}^{\times \delta}$. This action is what we call the {\em
Hodge decomposition of $\sch{X^{\op{top}}}{\mathbb{C}}$}.

\subsection{Review of the non-abelian Hodge correspondence} \label{ss-review}

In this section we recall some results concerning the non-abelian
Hodge correspondence of \cite{s2}, that will be needed for the proof
of the formality theorem, and for the construction of the Hodge
decomposition.

Let $X$ be a smooth projective algebraic variety over $\mathbb{C}$,
and let $X^{\op{top}}$ be the underlying topological space (in the
complex topology). As we explained in the first chapter the functor
$\opi{Sing}$ allows us to also view $X^{\op{top}}$ as a simplicial
set. Fix a base point $x \in X^{\op{top}}$. Let $L_{B}(X)$ be the
category of semi-simple local systems of finite dimensional
$\mathbb{C}$-vector spaces on $X^{\op{top}}$. It is a rigid
$\mathbb{C}$-linear tensor category which is naturally equivalent to
the category of finite dimensional semi-simple representations of the
fundamental group $\pi_{1}(X^{\op{top}},x)$.

The category of semi-simple $C^{\infty}$-bundles with flat
connections on $X$ will be denoted by $L_{DR}(X)$.  The objects in
$L_{DR}(X)$ are pairs
$(V,\nabla)$, where $V$ is a $C^{\infty}$-bundle, $\nabla : V
\longrightarrow V\otimes A^{1}$
is an integrable connection, and $A^{1}$ is the sheaf of
$C^{\infty}$-differential forms on $X$.
$L_{DR}(X)$ is also a rigid $\mathbb{C}$-linear tensor
category with monoidal structure given by
\[
(V,\nabla_{V})\otimes (W,\nabla_{W})
:= (V\otimes W,\nabla_{V}\otimes \op{Id}_{W} + \op{Id}_{V}\otimes
\nabla_{W}).
\]
The functor which maps a flat bundle to its monodromy
representations at $x$, induces an equivalence (the
Riemann-Hilbert correspodence)
of tensor categories $L_{B}(X)\simeq L_{DR}(X)$.

Recall next that a Higgs bundle on $X$ is a
$C^{\infty}$-vector bundle $V$, together with
an operator $D'' : V \longrightarrow V\otimes A^{1}$,
satisfying $(D'')^{2}=0$, and
the Leibniz's rule $D''(a\cdot s)=\overline{\partial}(a)\cdot s+a\cdot
D''(s)$, for any section $s$ of $V$ and any
function $a$. We will denote by $L_{Dol}(X)$ the
category of poly-stable
Higgs bundles $(V,D'')$ with vanishing
rational Chern classes (see \cite[$\S 1$]{s2}). Again, $L_{Dol}(X)$ is a
rigid $\mathbb{C}$-linear tensor
category, with monoidal structure given by
$$
(V,D''_{V})\otimes (W,D''_{W})
:= (V\otimes W,D''_{V}\otimes \op{Id}_{W} + \op{Id}_{V}\otimes
D''_{W}).$$

A harmonic bundle on $X$, is a triple $(V,\nabla,D'')$, where $V$ is a
$C^{\infty}$-bundle such that $(V,\nabla)$ is a flat bundle, $(V,D'')$
is a Higgs bundle, and the operators $\nabla$ and $D''$ are related by
a harmonic metric (see \cite[$\S 1$]{s2}). A morphism of harmonic
bundles is a morphism of $C^{\infty}$-bundles which preserves $\nabla$
and $D''$ (actually preserving $\nabla$ or $D''$ is enough, see
\cite[Lemma~1.2]{s2}). The category of such harmonic bundles will be
denoted by $L_{D'}(X)$, and is again a $\mathbb{C}$-linear tensor
category.  The existence of the harmonic metric implies that
$(V,\nabla)$ is semi-simple, and that $(V,D'')$ is poly-stable with
vanishing chern classes. Therefore, there exist natural projections
$$
\xymatrix{
& L_{D'}(X) \ar[dl] \ar[dr] & \\
L_{DR}(X) & & L_{Dol}(X).}
$$
The essence of the non-abelian Hodge correspondence is captured in
the following theorem.

\begin{thm}\label{nhc}{{\rm (\cite[Theorem~1]{s2})}}
The natural projections
$$\xymatrix{
& L_{D'}(X) \ar[dl] \ar[dr] & \\
L_{DR}(X) & & L_{Dol}(X).}
$$
are equivalence of $\mathbb{C}$-linear tensor categories.
\end{thm}
\textit{Remark:} The two projection  functors of the previous theorem
are functorial in $X$.
As a consequence of the functoriality of the equivalence
$L_{DR}(X)\simeq L_{D'}(X)\simeq L_{Dol}(X)$, one
obtains that it is compatible with the fiber functors at $x \in X$:
\[
\begin{split}
L_{DR}(X) & \longrightarrow L_{DR}(\{x\})\simeq \sff{Vect}, \\
L_{D'}(X) & \longrightarrow L_{D'}(\{x\})\simeq \sff{Vect}, \\
L_{Dol}(X) & \longrightarrow L_{Dol}(\{x\})\simeq \sff{Vect}.
\end{split}
\]
For an object
$(V,\nabla) \in L_{DR}(X)$, one can form its de Rham complex of
$C^{\infty}$-differential forms
\[
(A_{DR}^{\bullet}(V),\nabla) := \xymatrix@1{A^{0}(V)
\ar[r]^-{\nabla} & A^{1}(V)
\ar[r]^-{\nabla} & \ldots  \ar[r]^-{\nabla} &
A^{n}(V) \ar[r]^-{\nabla} &  \ldots},
\]
where $A^{n}(V)$ is the space of global sections of the
$C^{\infty}$-bundle
$V\otimes A^{n}$, and $A^{n}$ is the sheaf of (complex valued) smooth
differential forms of degree $n$ on $X$. These complexes are functorial
in $(V,\nabla)$, and
compatible with the tensor products, in the sense that there exists
morphims of complexes
\[
(A^{\bullet}_{DR}(V),\nabla_{V})\otimes
(A^{\bullet}_{DR}(W),\nabla_{W}) \longrightarrow
(A^{\bullet}_{DR}(V\otimes W),\nabla_{V}\otimes \op{Id}_{W} +
\op{Id}_{V}\otimes \nabla_{W}),
\]
which are functorial, associative, commutative and unital in the arguments
$(V,\nabla_{V})$ and $(W,\nabla_{W})$. In other words,
the functor $(V,\nabla) \mapsto (A^{\bullet}_{DR}(V),\nabla)$, from
$L_{DR}(X)$ to the category of complexes, is
a symmetric pseudo-monoidal functor. We denote by
\[
H^{\bullet}_{DR}(V):=H^{\bullet}(A^{\bullet}_{DR}(V),\nabla),
\]
the cohomology of the de Rham complex of $(V,\nabla)$.

For an object $(V,D'')$ of $L_{Dol}(X)$, one can define its Dolbeault
complex of $C^{\infty}$-differential
forms as
$$(A_{Dol}^{\bullet}(V),D'') := \xymatrix@1{A^{0}(V) \ar[r]^-{D''} &
A^{1}(V) \ar[r]^-{D''} & \ldots  \ar[r]^-{D''} &
A^{n}(V) \ar[r]^-{D''} &  \ldots,}$$
Again, these complexes are functorial and compatible with the monoidal
structure, in the sense explained above.
We denote by
$$H^{\bullet}_{Dol}(V):=H^{\bullet}(A^{\bullet}_{Dol}(V),D''),$$
the cohomology of the Dolbeault complex of $(V,D'')$.

Finally, for a harmonic bundle $(V,\nabla,D'') \in L_{D'}(X)$, one can consider the
operator $D':=\nabla-D''$ on $V$, and
the sub-complex $(\op{Ker}(D'),D'')$ of the Dolbeault complex
$(A_{Dol}^{\bullet}(V),D'')$ consisting of
differential forms $\alpha$ with $D'(\alpha)=0$. Note that
$(\op{Ker}(D'),D'')=(\op{Ker}(D'),\nabla)$ is also the sub-complex of
the de Rham complex $(A_{DR}^{\bullet}(V),\nabla)$ of forms $\alpha$
with $D'(\alpha)=0$.
This complex will be denoted
by
\[
(A_{D'}^{\bullet}(V),\nabla)=(A^{\bullet}_{D'}(V),D'') :=
(\op{Ker}(D'),D'')=(\op{Ker}(D'),\nabla).
\]
Again, the functor $(V,\nabla,D'') \mapsto (A_{D'}(V),\nabla)$ is
compatible with the tensor product in the sense
explained above.
In this way we obtain a natural
diagram of complexes
\[
\xymatrix{
& (A_{D'}^{\bullet}(V),\nabla) =(A_{D'}^{\bullet}(V),D'')
\ar[dr] \ar[dl] & \\
(A_{DR}^{\bullet}(V),\nabla) &
 & (A_{Dol}^{\bullet}(V),D'').}
\]
An important property of the correspondence \ref{nhc} is the following
compatibility
and formality results with respect to de Rham and Dolbeault
complexes. In the following
theorem, $(H_{DR}^{\bullet}(V),0)$ (respectively
$(H^{\bullet}_{Dol}(V),0)$) denotes the complex which underlying
graded vector space
is $H^{\bullet}_{DR}(V)$ (respectively  $H^{\bullet}_{Dol}(V)$) and
with zero differential.

\begin{thm}\label{fo}{{\rm (\cite[Lemma~2.2]{s2})}}
Let $(V,\nabla) \in L_{DR}(X)$, and $(V,D'')\in L_{Dol}(X)$ be the
corresponding Higgs bundle
via Theorem \ref{nhc}. Then, there exists a functorial diagram of
quasi-isomorphisms of complexes
\begin{equation*}
({\bf formality}) \; \xymatrix{
(A^{\bullet}_{DR}(V),\nabla) & & \ar[ll] (A^{\bullet}_{D'},D)=(A^{\bullet}_{D'}(V),D'') \ar[rr]
\ar[dll] \ar[drr] & & (A^{\bullet}_{Dol}(V),D'') \\
 (H^{\bullet}_{Dol}(V),0) & & & & (H^{\bullet}_{DR}(V),0) }
\end{equation*}
The diagram {\em ({\bf formality})} is moreover functorial
for pull-backs along morphisms of pointed smooth
and projective manifolds $Y \rightarrow X$.
\end{thm}

\

\medskip

\noindent
In the next section we will use theorem \ref{fo} extended to the
case of $\opi{Ind}$-objects.
The reader may find
definitions and general properties of categories of $\opi{Ind}$-objects in
\cite{sa,sga}. The categories of $\opi{Ind}$-objects
(belonging to the universe $\mathbb{U}$) in $L_{B}(X)$, $L_{DR}(X)$,
$L_{Dol}(X)$ and $L_{D'}(X)$ will be denoted respectively by
$T_{B}(X)$, $T_{DR}(X)$, $T_{Dol}(X)$, and $T_{D'}(X)$.
The equivalences of theorem \ref{nhc} extend verbatim to $\mathbb{C}$-linear
tensor
equivalences
\[
T_{B}(X)\simeq T_{DR} \simeq T_{D'}(X) \simeq T_{Dol}(X).
\]
Similarly the pseudo-monoidal functors
\[
\xymatrix@R=1pt@C=9pt{
L_{DR}(X) \ar[r] & \sff{C}^{+} \\
(V,\nabla) \ar@{|->}[r] & (A_{DR}^{\bullet}(V),\nabla), \\
 & \\
L_{Dol}(X) \ar[r] & \sff{C}^{+} \\
(V,D'')  \ar@{|->}[r] &  (A_{Dol}^{\bullet}(V),D''), \\
 & \\
L_{D'}(X) \ar[r] & \sff{C}^{+} \\
(V,\nabla,D')  \ar@{|->}[r] &  (A_{D'}^{\bullet}(V),D''),
}
\]
to the category $\sff{C}^{+}$ of complexes of $\mathbb{C}$-vector
spaces concentrated in non-negative degrees, 
extend naturally to pseudo-monoidal
functors from the categories of $\opi{Ind}$-objects
\[
\xymatrix@R=1pt@C=10pt{
T_{DR}(X) \ar[r] & \sff{C}^{+} \\
\{(V_{i},\nabla_{i})\}_{i \in I} \ar@{|->}[r] &  \op{colim}_{i \in I}
(A_{DR}^{\bullet}(V_{i}),\nabla_{i}) \\
& \\
T_{Dol}(X) \ar[r] & \sff{C}^{+} \\
\{(V_{i},D''_{i})\}_{i \in I} \ar@{|->}[r] &
\op{colim}_{i \in I} (A_{Dol}^{\bullet}(V_{i}),D''_{i}), \\
& \\
T_{D'}(X) \ar[r] & \sff{C}^{+} \\
\{(V_{i},\nabla_{i},D''_{i})\}_{i \in I} \ar@{|->}[r] &
\op{colim}_{i \in I} (A_{D'}^{\bullet}(V_{i}),D''_{i}).
}
\]
By convention, objects in $T_{DR}$ (respectively $T_{Dol}$,
respectively $T_{D'}$)
will again be denoted by
$(V,\nabla)$ (respectively $(V,D'')$, respectively
$(V,\nabla,D'')$). If $(V,\nabla)$
(respectively $(V,D'')$, respectively $(V,\nabla,D'')$) is the
$\opi{Ind}$-object
$\{(V_{i},\nabla_{i})\}_{i \in I} \in T_{DR}(X)$ (respectively
$\{(V_{i},D''_{i})\}_{i \in I}
\in T_{Dol}(X)$, respectively $\{(V_{i},\nabla_{i},D''_{i})\}_{i \in
I} \in T_{D'}(X)$),
we will also put
\[
\begin{split}
(A^{\bullet}_{DR}(V),\nabla) & := \op{colim}_{i \in I}
(A_{DR}^{\bullet}(V_{i}),\nabla_{i}) \\
(A^{\bullet}_{Dol}(V),D'') & := \op{colim}_{i \in I}
(A_{Dol}^{\bullet}(V_{i}),D''_{i}) \\
(A^{\bullet}_{D'}(V),D'') & := \op{colim}_{i \in I}
(A_{D'}^{\bullet}(V_{i}),D''_{i}).
\end{split}
\]
Theorem \ref{fo} then extends to the following formality result for
$\opi{Ind}$-objects.
\begin{cor}\label{cfo}
Let $(V,\nabla) \in T_{DR}(X)$, and $(V,D'')\in T_{Dol}(X)$ be the
corresponding Higgs bundle
via Theorem \ref{nhc} for $\opi{Ind}$-objects.
Then, there exists a functorial diagram of quasi-isomorphisms of complexes
\begin{equation*}
\tag{{\it Ind}-{\bf formality}} \xymatrix{
(A^{\bullet}_{DR}(V),\nabla) & & \ar[ll]
(A^{\bullet}_{D'},D)=(A^{\bullet}_{D'}(V),D'') \ar[rr]
\ar[dll] \ar[drr] & & (A^{\bullet}_{Dol}(V),D'') \\
 (H^{\bullet}_{Dol}(V),0) & & & & (H^{\bullet}_{DR}(V),0) }
\end{equation*}
The diagram {\em ({\it Ind}-{\bf formality})} is moreover functorial
for pullbacks along morphisms of pointed smooth
projective manifolds $Y \rightarrow X$.
\end{cor}
The categories $T_{B}(X)$, $T_{DR}(X)$,
$T_{Dol}(X)$ and $T_{D'}(X)$ are naturally $\mathbb{C}$-linear tensor abelian
categories. Furthermore, they
possess all $\mathbb{U}$-limits and $\mathbb{U}$-colimits.
This last property allows one to
define an action of the monoidal category $\sff{Vect}$, of vector spaces
in $\mathbb{U}$, on the categories
$T_{B}(X)$, $T_{DR}(X)$,  $T_{Dol}(X)$ and $T_{D'}(X)$, making them
into closed module categories
over $\sff{Vect}$ (see \cite[$\S 4.2$]{ho}). Precisely, this means that there
exists bilinear functors (external products)
\[
\otimes : \sff{Vect}\times T_{?}(X) \longrightarrow T_{?}(X),
\]
where $T_{?}(X)$ is one of the categories $T_{B}(X)$, $T_{DR}(X)$,
$T_{Dol}(X)$  and $T_{D'}(X)$. Moreover these functors
satisfy the usual adjunction property
\[
\op{Hom}(A\otimes V,W)\simeq
\underline{\op{Hom}}_{\sff{Vect}}(A,\underline{\op{Hom}}(V,W)),
\]
where $\underline{\op{Hom}}(V,W) \in \sff{Vect}$ is the vector space
of morphisms coming from the linear
structure on $T_{?}(X)$.

Using these external products, one can define the notion
of an action of an affine group scheme $G$ on
an object $V \in T_{?}(X)$. Indeed, let $\mathcal{O}(G)$
be the Hopf algebra of functions on $G$. Then,
a $G$-action on $V \in T_{?}(X)$ is a morphism
\[
V \longrightarrow \mathcal{O}(G)\otimes V,
\]
which turns $V$ into a
co-module over the co-algebra $\mathcal{O}(G)$.
Dually an action of $G$ on $V$ is the
data of a factorization
of the functor
\[
\xymatrix@R=1pt@C=9pt{
h_{V} : \;  T_{?}(X) \ar[r] & \sff{Vect} \\
V' \ar@{|->}[r] & \underline{\op{Hom}}(V',V)
}
\]
through the category of linear representations of $G$.
These definitions will be used in the next section, and we will
talk freely about
action of an affine group scheme $G$ on an object of
$T_{B}(X)$, $T_{DR}(X)$, $T_{Dol}(X)$ and $T_{D'}(X)$.

\

\begin{rem} \label{rem-actions} Let  ${\mathcal
C}$ be a Tannakian category with ${\mathcal C} = \sff{Rep}(G)$ for some
pro-algebraic group $G$. There are two natural actions (by algebra
automorphisms) of the group $G$ on the algebra of functions ${\mathcal
O}(G)$. These actions are
induced respectively by the rigtht and left  translation action of $G$
on itself. If we now view $\mathcal{O}(G)$ as an algebra object in
$\mathcal{C}$ via say the left action, then the right action of $G$
turns $\mathcal{O}(G)$ into a $G$-equivariant algebra object in
$\mathcal{C}$.
\end{rem}

\subsection{Schematization and differential forms} \label{ss-difforms}

In this section $(X,x)$ will denote
a pointed connected compact smooth manifold and $X^{\op{top}}$ will
denote the underlying
topological space of $X$. Let $L_{B}(X)$ be the
category of semi-simple local systems of
finite dimensional $\mathbb{C}$-vector spaces on $X^{\op{top}}$. It is
a rigid $\mathbb{C}$-linear tensor category which is naturally equivalent
to the category of finite dimensional semi-simple
representations of the fundamental group $\pi_{1}(X^{\op{top}},x)$.

The category of semi-simple $C^{\infty}$ complex vector bundles with flat
connections on $X$ will be denoted as before by $L_{DR}(X)$. Recall
that the category $L_{DR}(X)$ is a rigid $\mathbb{C}$-linear tensor
category, and the functor which maps a flat bundle to its monodromy
representations at $x$, induces an equivalence
of tensor categories $L_{B}(X)\simeq L_{DR}(X)$ (this again is the
Riemann-Hilbert correspondence).

Let $G_{X}:=\pi_{1}(X^{\op{top}},x)^{\op{red}}$ be the pro-reductive
completion of the group $\pi_{1}(X^{\op{top}},x)$. Note that it is the
Tannaka dual of the category $L_{B}(X)$.  The algebra 
$\mathcal{O}(G_{X})$ of regular functions on $G_{X}$ 
can be viewed as the left regular representation
of $G_{X}$. Through the universal morphism $\pi_{1}(X^{\op{top}},x)
\longrightarrow G_{X}$, we can also consider ${\mathcal O}(G_{X})$ as
a linear representation of $\pi_{1}(X^{\op{top}},x)$. This linear
representation is not finite dimensional, but it is admissible in the
sense that it equals the union of its finite dimensional
sub-representations. Therefore, the algebra $\mathcal{O}(G_{X})$
corresponds to an object in the $\mathbb{C}$-linear tensor category
$T_{B}(X)$, of $\opi{Ind}$-local systems on $X^{\op{top}}$. By
convention all of our $\opi{Ind}$-objects are labeled by
$\mathbb{U}$-small index categories.

Furthermore, the algebra structure on $\mathcal{O}(G_{X})$, gives rise
to a map
\[
\mu : \mathcal{O}(G_{X}) \otimes \mathcal{O}(G_{X}) \longrightarrow
\mathcal{O}(G_{X}),
\]
which is easily checked to be a morphism in
$T_{B}(X)$. This means that if we write
$\mathcal{O}(G_{X})$ as
the colimit of finite dimensional local systems $\{V_{i}\}_{i \in
I}$, then the product $\mu$ will be given
by a compatible system of morphisms in $L_{B}(X)$
\[
\mu_{i,k} : V_{i} \otimes V_{i} \longrightarrow V_{k},
\]
for some index $k \in I$ with $i\leq k \in I$.  The morphism $\mu$ (or
equivalentely the collection of morphisms $\mu_{i,k}$), endows the
object $\mathcal{O}(G_{X}) \in T_{B}(X)$ with a structure of a
commutative unital monoid.  Through the Riemann-Hilbert correspondence
$T_{B}(X)\simeq T_{DR}(X)$, the algebra $\mathcal{O}(G_{X})$ can also
be considered as a commutative monoid in the tensor category
$T_{DR}(X)$ of $\opi{Ind}$-objects in $L_{DR}(X)$.

Let $\{(V_{i},\nabla_{i})\}_{i \in I} \in T_{DR}(X)$ be the
object corresponding to $\mathcal{O}(G_{X})$.
For any $i \in I$, one can form the de Rham complex of
$C^{\infty}$-differential forms
$$(A^{\bullet}_{DR}(V_{i}),\nabla_{i}):=\xymatrix@1{A^{0}(V_{i})
\ar[r]^-{\nabla_{i}} & A^{1}(V_{i})
\ar[r]^-{\nabla_{i}} & \ldots  \ar[r]^-{\nabla_{i}}&
A^{n}(V_{i}) \ar[r]^-{D_{i}} &  \ldots}.$$
In this way we obtain an inductive system of complexes
$\{(A^{\bullet}_{DR}(V_{i}),D_{i})\}_{i \in I}$ whose colimit
complex was defined to be the de Rham complex of $\mathcal{O}(G_{X})$
$$(A^{\bullet}_{DR}(\mathcal{O}(G_{X})),\nabla) :=
\op{colim}_{i \in I}(A^{\bullet}_{DR}(V_{i}),\nabla_{i}).$$
The complex $(A^{\bullet}_{DR}(\mathcal{O}(G_{X})),\nabla)$ has a
natural structure of a commutative differential graded algebra,
coming from the commutative monoid structure on
$\{(V_{i},\nabla_{i})\}_{i \in I} \in T_{DR}(X)$. Using wedge products
of differential forms,
these morphisms induce in the usual fashion morphisms of complexes
$$(A^{\bullet}_{DR}(V_{i}),\nabla_{i}) \otimes
(A^{\bullet}_{DR}(V_{j}),\nabla_{j}) \longrightarrow
(A^{\bullet}_{DR}(V_{k}),\nabla_{k})$$
which, after passing to the colimit along $I$, turn
$(A^{\bullet}_{DR}(\mathcal{O}(G_{X})),\nabla)$ into a commutative
differential graded algebra.

The affine group scheme $G_{X}$ acts via the right regular
representation on the $\opi{Ind}$-local system
$\mathcal{O}(G_{X})$. As explained in Remark~\ref{rem-actions} this
action is compatible with the algebra structure. By functoriality, we
get an action of $G_{X}$ on the corresponding object in
$T_{DR}(X)$. Furthermore, if $G_{X}$ acts on an inductive system of
flat bundles $(V_{i},\nabla_{i})$, then it acts naturally on its de
Rham complex $\op{colim}_{i \in
I}(A^{\bullet}_{DR}(V_{i}),\nabla_{i})$, by acting on the spaces of
differential forms with coefficients in the various $V_{i}$. Indeed,
if the action of $G_{X}$ is given by a co-module structure
\[
\{(V_{i},\nabla_{i})\}_{i \in I}\longrightarrow \{\mathcal{O}(G_{X})
\otimes (V_{i},\nabla_{i})\}_{i \in I},
\]
then one obtains a morphism of $\opi{Ind-C}^{\infty}$-bundles by tensoring
with the sheaf $A^{n}$ of differential
forms on $X$
\[
\{V_{i}\otimes A^{n}\}_{i \in I}\longrightarrow
\{\mathcal{O}(G_{X})\otimes (V_{i}\otimes A^{n})\}_{i \in I}.
\]
Taking global sections on $X$, one has a morphism
\[
\op{colim}_{i \in I}A^{n}(V_{i}) \longrightarrow \op{colim}_{i \in I}
A^{n}(V_{i})\otimes \mathcal{O}(G_{X}),
\]
which defines an action of $G_{X}$ on the space of differential forms
with values in the $\opi{Ind-C}^{\infty}$-bundle $\{V_{i}\}_{i \in
I}$. Since this action is compatible with the differentials
$\nabla_{i}$, one obtains an action of $G_{X}$ on the de Rham complex
$\op{colim}_{i \in I
}(A^{\bullet}_{DR}(V_{i}),\nabla_{i})$. Furthermore since the action
is compatible with the algebra structure on $\mathcal{O}(G_{X})$ it
follows that $G_{X}$ acts on $\op{colim}_{i \in
I}(A^{\bullet}_{DR}(V_{i}),\nabla_{i})$ by algebra automorphisms.
Thus, the group scheme $G_{X}$ acts in a natural way on the complex
$(A^{\bullet}_{DR}(\mathcal{O}(G_{X})),\nabla)$, turning it into a
well defined $G_{X}$-equivariant commutative differential graded
algebra.

Using the denormalization functor (see \eqref{eq:equivariantD})
$D : \op{Ho}(\opi{G_{X}-\sff{CDGA}}) \longrightarrow
\op{Ho}(\opi{G_{X}-\sff{Alg}}^{\Delta})$,
we obtain a well defined $G_{X}$-equivariant co-simplicial algebra denoted by
\[
C^{\bullet}_{DR}(X,\mathcal{O}(G_{X})):=
D(A^{\bullet}_{DR}(\mathcal{O}(G_{X})),\nabla)
\in \op{Ho}(\opi{G_{X}-\sff{Alg}}^{\Delta}).
\]
To summarize:
\begin{df}\label{ddiff}
Let $(X,x)$ be a pointed connected smooth manifold, and let
$G_{X}:=\pi_{1}(X,x)^{\op{red}}$ be the pro-reductive completion of
its fundamental group. The \emph{$G_{X}$-equivariant commutative
differential graded algebra of de Rham cochains of $X$ with coefficients in
$\mathcal{O}(G_{X})$} will be denoted by
\[
(A^{\bullet}_{DR}(\mathcal{O}(G_{X})),\nabla) \in
\op{Ho}(\opi{G_{X}-\sff{CDGA}}).
\]
Its denormalization will be denoted by
\[
C^{\bullet}_{DR}(X,\mathcal{O}(G_{X})):=
D(A^{\bullet}_{DR}(\mathcal{O}(G_{X})),\nabla) \in
\op{Ho}(\opi{G_{X}-\sff{Alg}}^{\Delta}).
\]
\end{df}
\

\smallskip

\noindent
Any smooth map $f : (Y,y) \longrightarrow
(X,x)$ of pointed connected
smooth manifolds induces a morphism
$G_{Y}:=\pi_{1}(Y,y)^{\op{red}} \longrightarrow
G_{X}:=\pi_{1}(X,x)^{\op{red}}$, and 
therefore a well defined functor
\[
f^{*} : \op{Ho}(\opi{G_{X}-\sff{Alg}}^{\Delta}) \longrightarrow
\op{Ho}(\opi{G_{Y}-\sff{Alg}}^{\Delta}).
\]
Clearly the pull-back of differential forms
via $f$ induces
a well defined morphism in $\op{Ho}(\opi{G_{Y}-\sff{Alg}}^{\Delta})$
\[
f^{*} : f^{*}C^{\bullet}_{DR}(X,\mathcal{O}(G_{X}))
\longrightarrow C^{\bullet}_{DR}(Y,\mathcal{O}(G_{Y})).
\]
Since this morphism depends functorially (in an obvious
fashion) on the
morphism $f$, we get a well defined functor
\[
(X,x) \mapsto [\mathbb{R}\op{Spec}_{G_{X}}\,
C^{\bullet}_{DR}(X,\mathcal{O}(G_{X}))/G_{X}],
\]
from the category of pointed connected smooth manifolds to the
category of pointed schematic
homotopy types. To simplify notation, we will denote this functor by
$(X,x) \mapsto (X\otimes \mathbb{C})^{\op{diff}}$.

\begin{prop}{(\cite[Proposition~4.8]{small})}\label{pdiff}
For a pointed connected smooth and compact manifold $(X,x)$, there
exists a functorial isomorphism
in $\op{Ho}(\sff{SPr}_{*}(\mathbb{C}))$
\[
\sch{X^{\op{top}}}{\mathbb{C}}\simeq
(X\otimes \mathbb{C})^{\op{diff}}.
\]
\end{prop}

\subsection{Fixed-point model categories}

Before we define the Hodge decomposition on the schematization
$(\sch{X^{\op{top}}}{\mathbb{C}},x)$,  we will need some preliminary
results about fixed-point model categories.

Let $M$ be a cofibrantly generated model category with
$\mathbb{U}$-limits and colimits, and let $\Gamma$ be a group in
$\mathbb{U}$. Suppose that the group $\Gamma$ acts on $M$ by
auto-equivalences. This means that the action is given by a monoidal
functor $\Gamma \longrightarrow \underline{\op{End}}(M)$ (see \cite[$\S
6$]{s2} for a precise definition).  For $\gamma\in \Gamma$ and $x \in
\operatorname{ob}(M)$, let $\gamma\cdot x$ denote the image of $x$
under the auto-equivalence $\gamma$. Similarly for a morphism $f$ in
$M$ we will use the notation $\gamma\cdot f$ for the image of $f$ by
the auto-equivalence $\gamma$.

\medskip

\noindent
Define a new category $M^{\Gamma}$ of
$\Gamma$-fixed points in $M$ as follows.
\begin{itemize}
\item An object in $M^{\Gamma}$ is the data consisting
of an object $x \in M$, together with isomorphisms
$u_{\gamma} : \gamma\cdot x \longrightarrow x$ specified for each
$\gamma \in 
\Gamma$, satisfying the relation
$u_{\gamma_{2}}\circ (\gamma_{2}\cdot u_{\gamma_{1}}) =
u_{\gamma_{2}\cdot \gamma_{1}}$. Such an object will be
denoted by $(x,u)$.
\item A morphism $(x,u) \longrightarrow (y,v)$, between two objects in
$M^{\Gamma}$, is a morphism $f : x \longrightarrow y$
in $M$, such that for any $\gamma \in \Gamma$ one has $(\gamma\cdot
f)\circ u_{\gamma}=v_{\gamma}\circ f$.
\end{itemize}

\

The following proposition is a particular case of the existence of a
model structure on \textit{the category of sections of a left Quillen
presheaf} introduced in \cite[Theorem $17.1$]{hs}. The proof is a
straitgforward generalization of \cite[Theorem $13.8.1$]{hi}, and is
left to the reader.

\begin{prop}
There exist a unique structure of a cofibrantly generated model
category on $M^{\Gamma}$ such that
a morphism $f : (x,u) \longrightarrow (y,v)$ is an equivalence
(respectively a fibration) if and only if
the underlying morphism in $M$, $f : x \longrightarrow y$ is an
equivalence (respectively a fibration).
\end{prop}

\

\medskip

From the definition of a fixed-point model category it follows
immediately that any $\Gamma$-equivariant left Quillen
functor $F : M \longrightarrow N$ induces a
left Quillen functor $F^{\Gamma}$
on the model category of fixed points $F^{\Gamma}
: M^{\Gamma} \longrightarrow N^{\Gamma}$.
Furthermore, if $G^{\Gamma}$ is the right adjoint to $F^{\Gamma}$, then
the functors
\[
\mathbb{L}F^{\Gamma} : \op{Ho}(M^{\Gamma}) \longrightarrow
\op{Ho}(N^{\Gamma}), 
\qquad
\mathbb{R}G^{\Gamma} : \op{Ho}(N^{\Gamma}) \longrightarrow
\op{Ho}(M^{\Gamma}) 
\]
commute with the forgetful functors
\[
\op{Ho}(M^{\Gamma}) \longrightarrow \op{Ho}(M), \qquad \op{Ho}(N^{\Gamma})
\longrightarrow \op{Ho}(N).
\]
\

Our main example of a fixed-point model category will be the
following. Fix an affine group scheme $G$, and let $\Gamma$ be a group
acting on $G$, by group automorphisms. We will suppose that $\Gamma$
and $G$ both belong to $\mathbb{U}$.  The group $\Gamma$ acts
naturally (by auto-equivalences) on the category
$\opi{G-\sff{Alg}}^{\Delta}$, of $G$-equivariant co-simplicial
algebras. To be more precise, the group $\Gamma$ acts on the Hopf
algebra $B$ corresponding to $G$. Then, for any co-module $u :
M\rightarrow M\otimes B$ over $B$, and $\gamma\in \Gamma$, one defines
$\gamma\cdot M$ to be the $B$-co-module with co-action
\[
\xymatrix@1{M \ar[r]^-{u} & M\otimes B \ar[r]^-{\op{Id}\otimes
\gamma} & M\otimes B.}
\]
This defines an action of the group $\Gamma$ on the category
$\sff{Rep}(G)$ of linear representations of $G$. By functoriality,
this action extends to an action on the category of co-simplicial
objects, $\sff{Rep}(G)^{\Delta}$. Finally, as this action is by
monoidal auto-equivalences, it induces a natural action on the
category of $G$-equivariant co-simplicial algebras,
$\opi{G-\sff{Alg}}^{\Delta}$.

In the same way, the group $\Gamma$ acts on the category
$\opi{G-\sff{SPr}}(\mathbb{C})$, of $G$-equivariant simplicial
presheaves. For any $F \in \opi{G-\sff{SPr}}(\mathbb{C})$ and $\gamma
\in \Gamma$, we write $\gamma\cdot F$ for the $G$-equivariant
simplicial presheaf whose underlying simplicial presheaf is the same
as the one for $F$ and whose $G$-action is defined by composing the
$G$-action on $F$ with the automorphism $\gamma : G \longrightarrow
G$.

The right Quillen functor $\op{Spec}_{G} :
(\opi{G-\sff{Alg}}^{\Delta})^{op} \longrightarrow
\opi{G-\sff{SPr}}(\mathbb{C})$ commutes with the action of $\Gamma$,
and therefore induces a right Quillen functor on the model categories
of fixed points:
\[
\op{Spec}_{G}^{\Gamma} : ((\opi{G-\sff{Alg}}^{\Delta})^{\Gamma})^{op}
\longrightarrow (\opi{G-\sff{SPr}}(\mathbb{C}))^{\Gamma}.
\]
The group $\Gamma$ also acts on the comma category
$\sff{SPr}(\mathbb{C})/\opi{BG}$. For $\gamma \in \Gamma$, and
$u : F \rightarrow \opi{BG}$ an object in
$\sff{SPr}(\mathbb{C})/\opi{BG}$, $\gamma\cdot u : F \rightarrow \opi{BG}$
is defined by the
composition
\[
\xymatrix@1{F \ar[r]^-{u} & \opi{BG} \ar[r]^-{\gamma} & \opi{BG}.}
\]
\

One can then consider the left Quillen functor $\opi{De} :
\opi{G-\sff{SPr}}(\mathbb{C}) \longrightarrow
\sff{SPr}(\mathbb{C})/\opi{BG}$ defined
in section~\ref{ss:equivariant}. Since this functor is
$\Gamma$-equivariant, it induces a 
functor of the categories of fixed points
\[
\opi{De}^{\Gamma} : (\opi{G-\sff{SPr}}(\mathbb{C}))^{\Gamma}
\longrightarrow (\sff{SPr}(\mathbb{C})/\opi{BG})^{\Gamma},
\]
which is a Quillen equivalence (see Lemma \ref{l1}). For $F \in
\op{Ho}((\opi{G-\sff{SPr}}(\mathbb{C}))^{\Gamma})$, the corresponding object
in $\op{Ho}(\sff{SPr}(\mathbb{C})/\opi{BG})^{\Gamma})$ will be denoted by
$[F/G]^{\Gamma}$. By composing with the right derived functor
of $\op{Spec}_{G}^{\Gamma}$ one obtains a functor
\[
[\mathbb{R}\op{Spec}_{G}^{\Gamma}(-)/G]^{\Gamma} :
\op{Ho}((\opi{G-\sff{Alg}}^{\Delta})^{\Gamma})^{\op{op}} \longrightarrow
\op{Ho}((\sff{SPr}(\mathbb{C})/\opi{BG})^{\Gamma}).
\]
Finally, the forgetful functor $\sff{SPr}(\mathbb{C})/\opi{BG}
\longrightarrow \sff{SPr}(\mathbb{C})$ becomes $\Gamma$-equivariant if we
endow $\sff{SPr}(\mathbb{C})$ with the trivial $\Gamma$-action. The forgetful
functor induces a well defined functor
\[
\op{Ho}((\sff{SPr}(\mathbb{C})/\opi{BG})^{\Gamma}) \longrightarrow
\op{Ho}(\sff{SPr}(\mathbb{C})^{\Gamma}).
\]
By composing with $[\mathbb{R}\op{Spec}_{G}^{\Gamma}(-)/G]^{\Gamma}$
one gets a new
functor, which will is the one we are ultimately interested in
\[
[\mathbb{R}\op{Spec}_{G}^{\Gamma}(-)/G] :
\op{Ho}((\opi{G-\sff{Alg}}^{\Delta})^{\Gamma})^{op} \longrightarrow
\op{Ho}(\sff{SPr}(\mathbb{C})^{\Gamma}).
\]
This last functor sends a $G$-equivariant co-simplicial algebra fixed
under the $\Gamma$ action to a $\Gamma$-equivariant
stack.

\subsection{Construction of the Hodge decomposition}

First recall the notion of an action of a group $\Gamma$ on a pointed
stack $F \in \op{Ho}(\sff{SPr}_{*}(\mathbb{C}))$.

\begin{df}
An \emph{action of a group $\Gamma$ on a pointed stack} $F \in
\op{Ho}(\sff{SPr}_{*}(\mathbb{C}))$, is
a pair $(F_{0},u)$, consisting of a $\Gamma$-equivariant pointed stack
$F_{0} \in \op{Ho}(\opi{\Gamma-\sff{SPr}}_{*}(\mathbb{C}))$,
and an isomorphism $u$ of pointed stacks (i.e. an isomorphism in
$\op{Ho}(\sff{SPr}_{*}(\mathbb{C}))$) between $F$ and
the pointed stack underlying  $F_{0}$.

A \emph{morphism $(F,F_{0},u) \longrightarrow (F',F'_{0},u')$ between two
stacks $F$ and $F'$ equipped with $\Gamma$-actions}, is the data of two
morphisms $a : F \rightarrow F'$ and $b : F_{0} \rightarrow F_{0}'$,
in $\op{Ho}(\sff{SPr}_{*}(\mathbb{C}))$ and
$\op{Ho}(\opi{\Gamma-\sff{SPr}}_{*}(\mathbb{C}))$,
such that the
following diagram
\[
\xymatrix{
F \ar[r]^-{a} \ar[d]_-{u} & F' \ar[d]^-{u'} \\
F_{0} \ar[r]_-{b} & F_{0}'}
\]
commutes in $\op{Ho}(\sff{SPr}_{*}(\mathbb{C}))$.
\end{df}

\

With these definitions, the pointed stacks equipped with a $\Gamma$
action form a category. By definition, this category is also the
$2$-fiber product of the categories
$\op{Ho}(\sff{SPr}_{*}(\mathbb{C}))$ and
$\op{Ho}(\opi{\Gamma-\sff{SPr}}_{*}(\mathbb{C}))$ over
$\op{Ho}(\sff{SPr}_{*}(\mathbb{C}))$, and is therefore equivalent to
$\op{Ho}(\opi{\Gamma-\sff{SPr}}_{*}(\mathbb{C}))$ via the functor
$(F,F_{0},u) \mapsto F_{0}$. Using this equivalence, we will view
the pointed stacks equipped with a $\Gamma$ action as
objects in $\op{Ho}(\opi{\Gamma-\sff{SPr}}_{*}(\mathbb{C}))$.

\

\medskip

Now, fix a complex smooth projective algebraic variety $X$, a point
$x \in X$ and
let $X^{\op{top}}$ be the underlying topological space for the
classical topology. 
Consider the fiber-at-$x$ functor
\[
\omega_{x} : L_{Dol}(X) \longrightarrow \sff{Vect}.
\]
With this choice $(L_{Dol}(X),\otimes,\omega_{x})$ becomes a neutral
Tannakian category, whose affine group scheme of tensor automorphisms
will be denoted by $G_{X}$. In \cite{s2} C.Simpson introduced and
studied and action of $\mathbb{C}^{\times \delta}$ on $G_{X}$.  The
discrete group $\mathbb{C}^{\times\delta}$ acts on the category
$L_{Dol}(X)$ as follows. For $\lambda \in \mathbb{C}^{\times\delta}$,
and $(V,D'')$,write $D''=\overline{\partial}+\theta$, where $\theta$
is the Higgs field of $(V,D'')$.  Then, we define $\lambda \cdot
(V,D'')$ to be $(V,\overline{\partial}+\lambda\cdot \theta)$. This
defines an action of $\mathbb{C}^{\times\delta}$ by tensor
auto-equivalences of $L_{Dol}(X)$.  Since $\omega_{x}$ is invariant
under the action of the group $\mathbb{C}^{\times\delta}$ on
$L_{Dol}(X)$ it follows that group scheme $G_{X}$ is endowed with a
natural $\mathbb{C}^{\times\delta}$-action.

Let $T_{Dol}$ be the category of those $\opi{Ind}$-objects in
$L_{Dol}(X)$, which
belong to
$\mathbb{U}$. Recall that $T_{Dol}$ has a natural structure of a
$\mathbb{C}$-linear tensor category.
Furthermore, the action of $\mathbb{C}^{\times}$ extends to an action
by tensor auto-equivalences of $T_{Dol}$.
The resulting fiber functor $\omega_{x} : T_{Dol} \longrightarrow
\sff{Vect}$ has a right adjoint $p : \sff{Vect} \longrightarrow
T_{Dol}$, which is still $\mathbb{C}^{\times}$-invariant. The image
of the trivial algebra $\boldsymbol{1}$ under $p$
is therefore a monoid in $T_{Dol}$ (commutative and unital, as
usual). Moreover, as $p$ is fixed under
the action of $\mathbb{C}^{\times}$, $p(\boldsymbol{1})$ is naturally
fixed
under $\mathbb{C}^{\times}$ as a commutative monoid
in $T_{Dol}$. In other words
\[
p(\boldsymbol{1})\in (\sff{Comm}(T_{Dol}))^{\mathbb{C}^{\times}},
\]
where $\sff{Comm}(T_{Dol})$ denotes the category of commutative monoids in
$T_{Dol}$. Through the equivalence of $T_{Dol}$ with the category
of linear representations
of $G_{X}$, the commutative monoid $p(\boldsymbol{1})$ corresponds to the
left regular representation of $G_{X}$ on the algebra of functions
$\mathcal{O}(G_{X})$.

Consider now the adjunction morphism
\[
c : p(\boldsymbol{1}) \longrightarrow
p(\omega_{x}(p(\boldsymbol{1})))\simeq
p(\mathcal{O}(G)) \simeq \mathcal{O}(G_{X})\otimes p(\boldsymbol{1}),
\]
where the tensor product on the right is the external tensor product of
$p(\boldsymbol{1})$ with the vector space $\mathcal{O}(G_{X})$. The
morphism
$c$ defines a structure of $\mathcal{O}(G_{X})$-co-module on the
commutative monoid
$p(\boldsymbol{1})$. The fact that $p$ is fixed under
$\mathbb{C}^{\times\delta}$ implies that
$p(\boldsymbol{1})$ is naturally a fixed point of the category of
commutative monoids
in $T_{Dol}$, which are equipped with an action of the affine group
scheme $G_{X}$. Indeed, this is equivalent to the fact that the
multiplication
morphism $G_{X}\times G_{X} \longrightarrow G_{X}$ is
$\mathbb{C}^{\times\delta}$-equivariant.

The object $p(\boldsymbol{1})\in T_{Dol}$ is an inductive system of
objects in $L_{Dol}(X)$, i.e. an inductive system of Higgs bundles
$\{(V_{i},D''_{i})\}_{i \in I}$. One can consider for any $i \in I$,
the Dolbeault complex $(A^{\bullet}_{Dol}(V_{i}),D''_{i})$ (see $\S$
\ref{ss-review}), and thus define an inductive system of complexes
$(A^{\bullet}_{Dol}(V_{i}),D''_{i})_{i\in I}$. The colimit along $I$
of this inductive system was denoted by
$(A^{\bullet}_{Dol}(\boldsymbol{1}),D'')$. The commutative monoid
structure on $p(\boldsymbol{1})$ induces a well defined commutative
differential graded algebra structure on
$(A^{\bullet}_{Dol}(p(\boldsymbol{1})),D'')$.

Now, as $p(\boldsymbol{1})$ is a commutative monoid fixed
under the $\mathbb{C}^{\times\delta}$-action,
$(A^{\bullet}_{Dol}(p(\boldsymbol{1})),D'')$ possesses a natural
action of $\mathbb{C}^{\times\delta}$ defined in the following way. For
$\lambda \in
\mathbb{C}^{\times\delta}$, let
$\lambda \cdot p(\boldsymbol{1}) \in T_{Dol}(X)$ be the image of
$p(\boldsymbol{1})$ under the action of
$\mathbb{C}^{\times\delta}$ on $T_{Dol}(X)$, and
\[
u_{\lambda} : \lambda\cdot p(\boldsymbol{1}) \simeq p(\boldsymbol{1})
\]
be the isomorphism coming from the structure of being a fixed object.
We define
an isomorphism $\phi_{\lambda}$ of
$(A^{\bullet}_{Dol}(p(\boldsymbol{1})),D'')$ by the following
commutative
diagram
\[
\xymatrix{
(A^{\bullet}_{Dol}(p(\boldsymbol{1})),D'') \ar[r]^-{u_{\lambda}}
\ar[d]_-{\phi_{\lambda}} &
(A^{\bullet}_{Dol}(\lambda\cdot p(\boldsymbol{1})),D'')
\ar[d]^-{[\lambda]} \\
(A^{\bullet}_{Dol}(p(\boldsymbol{1})),D'') \ar[r]_-{u_{\lambda}} &
(A^{\bullet}_{Dol}(\lambda\cdot p(\boldsymbol{1})),D''),}
\]
where $[\lambda]$ is the automorphism of
$A^{\bullet}_{Dol}(\lambda\cdot p(\boldsymbol{1})),D'')$ which is
multiplication by $\lambda^{p}$ on the differential forms of type
$(p,q)$. The assignment $\lambda \mapsto \phi_{\lambda}$ defines an
action of $\mathbb{C}^{\times\delta}$ on the complex
$(A^{\bullet}_{Dol}(p(\boldsymbol{1})),D'')$. Moreover, this action is
compatible with the wegde product of differential forms, and is
therefore an action of $\mathbb{C}^{\times\delta}$ on the commutative
differential graded algebra
$(A^{\bullet}_{Dol}(p(\boldsymbol{1})),D'')$.

Finally, the action of $G_{X}$ on $p(\boldsymbol{1})$ induces a well
defined action of $G_{X}$ on the commutative differential graded
algebra $(A^{\bullet}_{Dol}(p(\boldsymbol{1})),D'')$. Here the
argument is the same as in the definition of the $G_{X}$-equivariant
commutative differential graded algebra
$(A_{DR}^{\bullet}(\mathcal{O}(G_{X})),D)$ given in $\S$
\ref{ss-difforms}.  The group scheme $G_{X}$ acts on the underlying
$\opi{Ind-C}^{\infty}$-bundle of $p(\boldsymbol(1))$, and therefore on
the spaces of differential forms with coefficients in this
bundle. This action preserves the operator $D''$ coming from the
structure of Higgs bundle, and therefore induces an action of $G_{X}$
on the Dolbeault complex $(A^{\bullet}_{Dol}(p(\boldsymbol{1})),D'')$.

This action is compatible with the action of
$\mathbb{C}^{\times\delta}$ on
$(A^{\bullet}_{Dol}(p(\boldsymbol{1})),D'')$ and on $G_{X}$, and gives
rise to a fixed point
\[
(A^{\bullet}_{Dol}(p(\boldsymbol{1})),D'') \in
(\opi{G_{X}-CDGA})^{\mathbb{C}^{\times}}.
\]

In summary we have defined the following object. Let
$\sff{CRep}(G_{X})$ be
the category of (positively graded) complexes of the
linear representations of $G_{X}$ which belong to $\mathbb{U}$. The group
$\mathbb{C}^{\times\delta}$ acts on $G_{X}$, and therefore
on the category $\opi{G_{X}-CDGA}$ of commutative and unital
algebras in $\opi{CRep}(G_{X})$.
The algebra $(A^{\bullet}_{Dol}(\boldsymbol{1}),D'')$ is then an
object in the fixed
category $(\opi{G_{X}-CDGA})^{\mathbb{C}^{\times}}$.
Using the denormalization functor which sends a differential graded
algebra to a
co-simplicial algebra, one can then
consider
\[
C^{\bullet}_{Dol}(X,\mathcal{O}(G_{X}))
:=D(A^{\bullet}_{Dol}(p(\boldsymbol{1})),D'') \in
\op{Ho}((\opi{G_{X}-\sff{Alg}}^{\Delta})^{\mathbb{C}^{\times\delta}}).
\]
This $G_{X}$-equivariant co-simplicial algebra is the Dolbeault
counter-part of $C^{\bullet}_{DR}(X,\mathcal{O}(G_{X}))$ appearing in
definition \ref{ddiff}.

\begin{prop} \label{prop-action}
Let $A:=C^{\bullet}_{Dol}(X,\mathcal{O}(G_{X})) \in
\op{Ho}((\opi{G_{X}-\sff{Alg}}^{\Delta})^{\mathbb{C}^{\times\delta}})$. 
\begin{itemize}
\item[(a)] After forgetting the action of $G_{X}$, there exist a natural
$\mathbb{C}^{\times\delta}$-equivariant
augmentation morphism in
$\op{Ho}((\sff{Alg}^{\Delta})^{\mathbb{C}^{\times\delta}})$ 
\[
e : A \longrightarrow \mathbb{C}.
\]
\item[(b)]
The underlying pointed stack of the
$\mathbb{C}^{\times\delta}$-equivariant pointed stack
\[
\xymatrix@1{
\op{Spec}\, \mathbb{C} \ar[r]^-{e^{*}} &
\mathbb{R}\op{Spec}^{\mathbb{C}^{\times\delta}}\, (A) \ar[r] &
[\mathbb{R}\op{Spec}_{G_{X}}^{\mathbb{C}^{\times\delta}}\, (A)/G_{X}]}
\]
is functorially isomorphic to $\sch{X^{\op{top}}}{\mathbb{C}}$ in
$\op{Ho}(\sff{SPr}_{*}(\mathbb{C}))$.
\end{itemize}
\end{prop}
{\bf Proof:} For the proof of (a) note that there is a natural
morphism of $\mathbb{C}^{\times\delta}$-fixed commutative differential
graded algebras
$A \longrightarrow A^{0}$,
obtained by projection on the $0$-th terms of the co-simplicial object
$A$.
By definition, $A^{0}$ is the
space of $C^{\infty}$-sections of the
$\opi{Ind}$-object $p(\boldsymbol{1}) \in T_{Dol}$.
In particular we can
can compose the above morphism
with the evaluation at $x \in X$, followed by the evaluation at
$e \in G_{X}$ to get a $\mathbb{C}^{\times}$-equivariant morphism
\[
A \longrightarrow A_{0} \longrightarrow p(\boldsymbol{1})_{x}\simeq
\mathcal{O}(G_{X})\longrightarrow \mathbb{C}.
\]
This proves part (a) of the proposition.

\smallskip

For the proof of (b) recall first
that the $G_{X}$-equivariant commutative differential graded algebra
$(A^{\bullet}_{Dol}(p(\boldsymbol{1})),D'')$
admits an interpretation as
the colimit of the inductive system
$\{(A_{Dol}^{\bullet}(V_{i}),D_{i}'')\}_{i \in I}$  of complexes with
$G_{X}$-action. Theorem \ref{nhc} and Theorem \ref{fo} (or
Corollary \ref{cfo}) imply
that $(A^{\bullet}_{Dol}(p(\boldsymbol{1})),D'')$ is naturally isomorphic
in $\op{Ho}(\opi{G_{X}-CDGA})$ to
the colimit of the corresponding inductive system of connections
$\{(A_{DR}^{\bullet}(V_{i}),\nabla_{i})\}_{i \in I}$.
On the other hand this inductive system of connections corresponds in
turn to the
$\opi{Ind}$-local system $\mathcal{O}(\pi_{1}(X,x)^{\op{red}})\simeq
\mathcal{O}(\pi_{1}(X,x)^{\op{red}})$.
Therefore,
$(A^{\bullet}_{Dol}(p(\boldsymbol{1})),D'')$ is naturally
isomorphic in $\op{Ho}(\opi{G_{X}-\sff{CDGA}})$ to
$(A^{\bullet}_{DR}(\mathcal{O}(G_{X}),\nabla)$ defined in the last
paragraph.
Proposition \ref{pdiff} implies the existence of a functorial
isomorphism in $\op{Ho}(\sff{SPr}_{*}(\mathbb{C}))$
\[
\sch{X^{\op{top}}}{\mathbb{C}}\simeq
[\mathbb{R}\op{Spec}_{G_{X}}\, A/G_{X}].
\]
In conclusion, we have constructed a $\mathbb{C}^{\times\delta}$-equivariant
stack $[\mathbb{R}\op{Spec}_{G_{X}}^{\mathbb{C}^{\times\delta}}\, A/G_{X}] \in
\op{Ho}(\opi{\mathbb{C}^{\times\delta}-\sff{SPr}}_{*}(\mathbb{C}))$, and
a functorial isomorphism between its underlying pointed stack and
$\sch{X^{\op{top}}}{\mathbb{C}}$.
By defintion, this is the data of a functorial
$\mathbb{C}^{\times\times}$-action on the pointed stack $(X^{\op{top}}\otimes
\mathbb{C})^{\op{sch}}$.
\hfill $\Box$

\

\bigskip

We are now
ready to define the Hodge decomposition:

\begin{df}\label{ddec}
The {\em Hodge decomposition} on the schematic homotopy type
$\sch{X^{\op{top}}}{\mathbb{C}}$ is the
$\mathbb{C}^{\times\delta}$-action defined by Proposition~\ref{prop-action}.
By
abuse of notation we will write
\[
\sch{X^{\op{top}}}{\mathbb{C}} \in
\op{Ho}(\sff{SPr}_{*}(\mathbb{C})^{\mathbb{C}^{\times}})
\]
for the schematic homotopy type of $X^{\op{top}}$ together with its Hodge
decomposition.
\end{df}

The functoriality statement in Theorem \ref{nhc} and \ref{fo}
 immediately implies that the Hodge
decomposition is functorial. Note however that the
construction depends on the choice
of the point $x \in X^{\op{top}}$ and so the functoriality is well defined
only for pointed morphisms of algebraic varieties. Finally we have the
following natural compatibility result:

\begin{thm}\label{tdec}
Let $X$ be a pointed smooth and projective algebraic variety over
$\mathbb{C}$, and let
\[
\sch{X^{\op{top}}}{\mathbb{C}}\in
\op{Ho}(\sff{SPr}_{*}(\mathbb{C})^{\mathbb{C}^{\times\delta}})
\]
be the schematization
of $X^{\op{top}}$ together with its
Hodge decomposition.
\begin{itemize}
\item[(i)] The induced action of $\mathbb{C}^{\times\delta}$ on the
cohomology groups
$H^{n}(\sch{X^{\op{top}}}{\mathbb{C}},\mathbb{G}_{a})\simeq
H^{n}(X^{\op{top}},\mathbb{C})$ is compatible with the Hodge decomposition
in the following sense. For any
$\lambda \in \mathbb{C}^{\times}$,
and $y \in H^{n-p}(X,\Omega_{X}^{p}) \subset
H^{n}(X^{\op{top}},\mathbb{C})$ one has
$\lambda(y)=\lambda^{p}\cdot y$.
\item[(ii)] Let $\pi_{1}(\sch{X^{\op{top}}}{\mathbb{C}},x)^{\op{red}}$
be the maximal reductive
quotient of $\pi_{1}(\sch{X^{\op{top}}}{\mathbb{C}},x))$.
The induced action of $\mathbb{C}^{\times\delta}$ on
$\pi_{1}(\sch{X^{\op{top}}}{\mathbb{C}},x)^{\op{red}}\simeq
\pi_{1}(X^{\op{top}},x)^{\op{red}}$
coincides with the one defined in \cite{s2}.
\item[(iii)] If $X^{\op{top}}$ is simply connected, then the induced action
of $\mathbb{C}^{\times}$ on
$$
\pi_{i}((X^{\op{top}}\otimes
\mathbb{C})^{\op{sch}},x)\simeq
\pi_{i}(X^{\op{top}},x)\otimes \mathbb{C}
$$
is compatible with the Hodge
decomposition defined in \cite{mo}. More precisely, suppose 
that $F^{\bullet}\pi_{i}(X^{\op{top}})\otimes \mathbb{C}$ is the
Hodge filtration
defined in \cite{mo}, then
\[
F^{p}\pi_{i}(X^{\op{top}})\otimes \mathbb{C}=\left\{ \left. x \in
\pi_{i}(X^{\op{top}}\otimes \mathbb{C}) \right|
\exists \; q\geq p, \text{ such that  } \lambda(x)=\lambda^{q}\cdot x,\;
\forall \lambda \in \mathbb{C}^{\times}\right\}.
\]
\item[(iv)] Let $R_{n}$ be the set of isomorphism classes of
$n$-dimensional simple
linear representations of
$\pi_{1}((X^{\op{top}}\otimes\mathbb{C})^{\op{sch}},x)$ 
Then, the induced action of $\mathbb{C}^{\times}$ on the set
\[
R_{n}\simeq
\op{Hom}(\pi_{1}(X^{\op{top}},x),Gl_{n}(\mathbb{C}))/Gl_{n}(\mathbb{C})
\]
defines a continuous action of the topological group
$\mathbb{C}^{\times}$
(for the analytic topology).
\end{itemize}
\end{thm}
{\bf Proof:} (i) Since the group scheme
$G_{X}:=\pi_{1}(X^{\op{top}},x)^{\op{red}}$ is reductive, its
Hochschild cohomology with 
coefficients in finite dimensional 
local systems vanishes. Therefore, one has
\[
H^{n}(\sch{X^{\op{top}}}{\mathbb{C}},\mathbb{G}_{a})\simeq
H^{n}(\mathbb{R}\op{Spec}_{G_{X}}\, A, \mathbb{G}_{a})^{G_{X}}\simeq
H^{n}(B)^{G_{X}}\simeq H^{n}(B^{G_{X}}),
\]
where $B=(A^{\bullet}_{Dol}(p(\boldsymbol{1})),D'')$ is the
$G_{X}$-equivariant
co-simplicial algebra defined during the construction of the Hodge
decomposition.
Clearly $B^{G_{X}}$ is isomorphic to
the commutative differential graded algebra of Dolbeault cochains
$A_{Dol}^{\bullet}(\boldsymbol{1},\overline{\partial})$, of the trivial
Higgs
bundle $(\boldsymbol{1},\overline{\partial})$.
Therefore,
\[
B^{G_{X}}\simeq
\bigoplus_{p}(A^{p,\bullet}(X),\overline{\partial})[-p].
\]
Moreover, the action of $\lambda \in \mathbb{C}^{\times}$ on $x \in
A^{p,q}(X)$ is given
by $\lambda(x)=\lambda^{p}\cdot x$. This implies that the induced
action of $\mathbb{C}^{\times}$
on $H^{n}(\sch{X}{\mathbb{C}},\mathbb{G}_{a})$ is the one
required.

\

\smallskip

\noindent
(ii) This is clear by construction.

\

\smallskip

\noindent
(iii) If $X^{\op{top}}$ is simply connected, then
$\pi_{1}(X^{\op{top}},x)^{\op{red}}=*$, and therefore one has
that
$\sch{X^{\op{top}}}{\mathbb{C}}$ is
naturally isomorphic to
\[
\mathbb{R}\op{Spec} (C^{\bullet}_{Dol}(X,\boldsymbol{1})),
\]
where $C^{\bullet}_{Dol}(X,\boldsymbol{1})$ is the dernormalization of
$(A_{Dol}^{\bullet}(\boldsymbol{1}),\overline{\partial}))$.
The action of $\mathbb{C}^{\times}$ has weight $p$ on
$A^{p,q}(X)$. This action also correspond to a decreasing filtration
$F^{\bullet}$ of $(A_{Dol}^{\bullet}(\boldsymbol{1}),\overline{\partial})$
by sub-$\sff{CDGA}$, where $F^{p}$ consists of differential forms of
type $(p',q)$ 
with $p'\geq p$.

The vector space $\pi_{i}(\sch{X}{\mathbb{C}},x)(\mathbb{C})\simeq
\pi_{i}(X,x)\otimes \mathbb{C}$ is the $i$-th homotopy group of the
mapping space
$\mathbb{R}\underline{\op{Hom}}_{cdga}((A_{Dol}^{\bullet}
(\boldsymbol{1}),\overline{\partial}),\mathbb{C})$,
which is in turn naturally isomorphic to the $i$-th homotopy group of
the commuttative dga $(A_{Dol}^{\bullet}(\boldsymbol{1}),\overline{\partial})$,
as defined in \cite{bg}. In other words, it is the dual of the space
of indecomposable elements of degree $i$ in the minimal model of
$(A_{Dol}^{\bullet}(\boldsymbol{1}),\overline{\partial})$. As the
action of $\mathbb{C}^{\times}$ on
$(A_{Dol}^{\bullet}(\boldsymbol{1}),\overline{\partial})$ is
compatible with the Hodge filtration $F^{\bullet}$, the induced action
of $\mathbb{C}^{\times}$ on $\pi_{i}(X,x)\otimes \mathbb{C}$ is
compatible with the filtration induced by $F^{\bullet}$ on the space
of indecomposables. This filtration being the Hodge filtration defined
in \cite{mo}, this proves (iii).

\

\smallskip

\noindent
(iv) This follows from (iii) and \cite[Proposition $1.5$]{s2}. \hfill
$\Box$

\

\section{Weights}

As we have already mentioned, the Hodge decomposition on
$\sch{X^{\op{top}}}{\mathbb{C}}$ constructed in the previous section
is only a part of a richer structure, which includes a \textit{weight
filtration}. In the simply connected case, the weight filtration is
constructed in \cite{mo}, and its existence reflects the fact that the
Hodge structures defined on the complexified homotopy groups are only
mixed Hodge structures. The existence of the weight filtration on the
rational homotopy type of a projective $X$ is related to the fact that
any simply connected homotopy type (or more generally any nilpotent
homotopy type) is a \textit{successive extensions of abelian homotopy
types}.

\subsection{Linear stacks and homology theory}

In this section we present a linear version of the 
schematic homotopy types, and use it to define
homology of schematic homotopy types. 

Let $\mathcal{O}$ be the sheaf on $(\sff{Aff}/\mathbb{C})_{\op{ffqc}}$
represented by the affine line ${\mathbb{A}}^{\! 1}$. By definition
$\mathcal{O}$ is a sheaf of $\mathbb{C}$-algebras and can be viewed as
an object in $\sff{SPr}(\mathbb{C})$. We consider the category
$\sff{SMod}(\mathcal{O})$, of objects in $\sff{SPr}(\mathbb{C})$
equipped with a structure of $\mathcal{O}$-modules (and morphisms
preserving this structure). In other words, $\sff{SMod}(\mathcal{O})$
is the category of simplicial objects in the category
$\sff{Mod}(\mathcal{O})$ of presheaves of $\mathcal{O}$-modules on
$\sff{Aff}/\mathbb{C}$ (note that we do not impose  any
sheaf conditions and consider all presheaves of
$\mathcal{O}$-modules).

Forgetting the $\mathcal{O}$-module structure defines a forgetful
functor
\[
\sff{SMod}(\mathcal{O}) \longrightarrow \sff{SPr}(\mathbb{C}).
\]
This functor has  a left adjoint
\[
-\otimes \mathcal{O} : \sff{SPr}(\mathbb{C}) \longrightarrow
\sff{SMod}(\mathcal{O}). 
\]
The category $\sff{SMod}(\mathcal{O})$ has a natural model structure
for which the fibrations and equivalences are defined in
$\sff{SPr}(\mathbb{C})$ through the forgetful functor. We will call
the model category $\sff{SMod}(\mathcal{O})$ the category of
\emph{linear stacks}. It is Quillen equivalent (via the Dold-Kan
correspondence) to the model category of negatively graded co-chain
complexes of $\mathcal{O}$-modules. In particular its homotopy
category $\op{Ho}(\sff{SMod}(\mathcal{O})$ is equivalent to
$D^{-}((\sff{Aff}/\mathbb{C})_{\op{ffqc}},\mathcal{O})$, the bounded
above derived category of $\mathcal{O}$-modules on the ringed site
$((\sff{Aff}/\mathbb{C})_{\op{ffqc}},\mathcal{O})$.

Let $E$ be a co-simplicial $\mathbb{C}$-vector space. We define
a linear stack $\op{Spel}\, E$ (\emph{Spel} stands for
\emph{\textbf{s}pectre \textbf{l}in\'eaire}) by the formula 
\[
\xymatrix@C-1pc@R-2pc{
\op{Spel}\, E : & \sff{Aff}/\mathbb{C} \ar[r] & \sff{SSet} \\
& A  \ar@{|->}[r] & \underline{\op{Hom}}(E,A),
}
\]
where $\underline{\op{Hom}}(E,A)$ is the simplicial set having
$\op{Hom}(E_{n},A)$ as sets of $n$-simplicies. The
$\mathcal{O}$-module structure on $\op{Spel}\, E$ is given by the
natural $A$-module structure on each simplicial set
$\underline{\op{Hom}}(E,A)$.  Note that the simplicial presheaf
underlying $\op{Spel}\, E$ is isomorphic to $\op{Spec}\,
(\op{Sym}^{\bullet}(E))$, where \linebreak $\op{Sym}^{\bullet}(E)$ is the free
commutative co-simplicial $\mathbb{C}$-algebra generated by the
co-simplicial vector space $E$.

The functor $\op{Spel}$ is a right Quillen functor
\[
\op{Spel} : \sff{Vect}^{\Delta} \longrightarrow
\sff{SMod}(\mathcal{O}),
\]
from the model category of co-simplicial vector spaces (with the usual model
structure for which equivalences and fibrations correspond to
quasi-isomorphisms 
and epimorphisms through the Dold-Kan correspondence) to the model
category of linear stacks.
Therefore the functor $\op{Spel}$ can be derived into a functor 
\[
\mathbb{R}\op{Spel} : \op{Ho}(\sff{Vect}^{\Delta}) \longrightarrow
\sff{SMod}(\mathcal{O}).
\]
(actually $\op{Spel}$ preserves equivalences, so $\op{Spel}\simeq
\mathbb{R}\op{Spel}$ Its is easy to check that this functor is fully
faithful when restricted to the full sub-category of
$\mathbb{U}$-small co-simplicial vector spaces.

\begin{df}\label{d20}
The essential image of the functor $\mathbb{R}\op{Spel}$, restricted
to $\mathbb{U}$-small co-simplicial vector spaces, is called
the category of \emph{linear affine stacks}.
\end{df}

Note that as $\mathbb{R}\op{Spel}\, E\simeq \mathbb{R}\op{Spec}\,
(\op{Sym}^{\bullet}(E))$, 
the forgetful 
functor
\[
\op{Ho}(\sff{SMod}(\mathcal{O})) \longrightarrow
\op{Ho}(\sff{SPr}(\mathbb{C})), 
\]
sends linear affine stacks to affine stacks.

Let $F$ be an affine stack, and let 
$A=\mathbb{L}\mathcal{O}(F)$ be the co-simplicial algebra of cohomology
on $F$, which can be chosen to belong to $\mathbb{U}$. Viewing $A$ as
a co-simplicial vector space by forgetting the multiplicative structure
one gets a linear affine stack
\[
\mathcal{H}(F,\mathcal{O}):=\op{Spel}\, A,
\]
which we will call the \emph{homology type of $F$}.  The functor $F
\mapsto \mathcal{H}(F,\mathcal{O})$ is left adjoint to the forgetful
functor from linear affine stacks to affine stacks. In particular, one
has an adjunction morphism of affine stacks
\[
F \longrightarrow \mathcal{H}(F,\mathcal{O}),
\]
and it is reasonable to interpret this morphism as the
\emph{abelianization} of the schematic homotopy type $F$.

By construction, the $i$-th homotopy sheaf of
$\mathcal{H}(F,\mathcal{O})$ is the linearly compact vector space dual
to $H^{i}(F,\mathcal{O})$ (a linearly compact vector space is an
object in the category of pro-vector spaces of finite dimension). It
is denoted by $H_{i}(F,\mathcal{O})$ and is called the $i$-th homology
sheaf of $F$. The adjunction morphism above induces Hurewitz
maps
\[
\pi_{i}(F) \longrightarrow H_{i}(F,\mathcal{O})
\]
which are well defined for all $i$.

\subsection{The weight tower}

First we introduce a general notion of a weight filtration on an
affine stack.

\

\smallskip

\begin{center} \textit{The weight tower of an affine stack} \end{center}

\

\smallskip

Let $\mathbb{T}$ be the category corresponding to the poset of natural
numbers $\mathbb{N}$.  The category
$\sff{SPr}_{*}(\mathbb{C})^{\mathbb{T}}$, of functors from
$\mathbb{T}$ to $\sff{SPr}_{*}(\mathbb{C})$ is a model category for
which equivalences and fibrations are defined levelwise (see
e.g. \cite{hi}). We will call the model category
$\sff{SPr}_{*}(\mathbb{C})^{\mathbb{T}}$ the \emph{model category of
towers in $\sff{SPr}_{*}(\mathbb{C})$}. We will call the objects in
the homotopy category
$\op{Ho}(\sff{SPr}_{*}(\mathbb{C})^{\mathbb{T}})$ \emph{towers of
pointed stacks}.

We will define a functor
\[
W^{*} : \op{Ho}(\sff{AS}_{*}^{c}) \longrightarrow
\op{Ho}(\sff{SPr}_{*}(\mathbb{C})^{\mathbb{T}}),
\]
 from the homotopy category of
pointed and connected affine stacks to the homotopy category of towers
of pointed stacks. This construction is an algebraic counter-part of a
well known topological construction (see \cite{curtis,curtis2}), and will be
achieved using the homotopy theory of simplicial affine group schemes
discussed in detail in \cite{small}. 

As explained in \cite{small}, there exists a model category of
($\mathbb{U}$-small) simplicial affine group schemes $\sff{sGAff}$, whose
equivalences are defined as the morphisms inducing quasi-isomorphisms
on the corresponding co-simplicial Hopf algebras. Recall also, that
the natural 'classifying stack' functor
\[
\xymatrix@R-2pc{
\sff{sGAff} \ar[r] & \sff{SPr}_{*}(\mathbb{C}) \\
G_{*} \ar@{|->}[r] & BG_{*},
}
\]
can be derived on the right to a functor
\[
\xymatrix@R-2pc{
\op{Ho}(\sff{sGAff}) \ar[r] & \op{Ho}(\sff{SPr}_{*}(\mathbb{C})) \\
G_{*} \ar@{|->}[r] & BR(G_{*})
}
\]
(here, $R$ is a fibrant replacement functor in the model category
$\sff{sGAff}$).  The functor $BR$ is fully faithful and its essential
image consists precisely of all pointed schematic homotopy types (see
\cite[Theorem~3.5]{small}).

The functor of complex points $H_{*} \mapsto H_{*}(\mathbb{C})$ is an
exact conservative (see \cite[III \S 3 Corollary 7.6]{dg}) functor
from the category $\sff{sGAff}$ to the category $\sff{SSet}$ of
simplicial sets in $\mathbb{U}$. Furthermore, this functor possesses a
left adjoint, sending a simplicial set to the free simplicial affine
group scheme it generates\footnote{The free affine group scheme
  generated by a set $I$ is the 
pro-algebraic completion of the free group over $I$}. Using this
adjunction, one can construct (as explained e.g. in \cite{il}) a
standard free resolution $L_{*}G_{*} \longrightarrow G_{*}$. The
object $L_{*}G_{*}$ is a bi-simplicial affine group scheme, such that
each $L_{*}G_{n}$ is the standard free resolution of $G_{n}$. We
denote by
$LG_{*}$ the diagonal of $L_{*}G_{*}$, and consider the induced
morphism of simplicial affine group schemes $LG_{*} \longrightarrow
G_{*}$.  The morphism $LG_{*}(\mathbb{C}) \longrightarrow
G_{*}(\mathbb{C})$ is a homotopy equivalence (see \cite{il}), which
implies that the morphism of pointed simplicial presheaves
\[
BLG_{*} \longrightarrow BG_{*}
\]
induces an isomorphisms on all
homotopy presheaves $\pi_{i}^{\op{pr}}$. In particular it induces an
isomorphism of pointed stacks.  This construction gives a functorial
morphism of simplicial affine group schemes $LG_{*} \longrightarrow
G_{*}$ such that the induced morphism
\[
BLG_{*} \longrightarrow BG_{*}
\]
is an equivalence of pointed stacks, and furthermore
each affine group scheme $LG_{n}$ is free. 

To any $G_{*}\in \sff{sGAff}$, we can now associate a new simplicial
affine group scheme $H \in \sff{sGAff}$ which is defined to be the
maximal unipotent quotient $H:=LG_{*}^{\op{uni}}$ of $LG_{*}$.  The
assignment $G_{*} \mapsto H_{*}$ can be thought of as the
\emph{left derived functor of the pro-unipotent completion
functor}.

We now consider the lower central series
\[
\dots \subset H_{*}^{(i)} \subset H_{*}^{(i-1)} \subset \dots \subset
H_{*}^{(1)} \subset H_{*}^{(0)}=H_{*},
\]
given for any $i\geq 1$ by
\[
H_{*}^{(i)}:=[H_{*}^{(i-1)},H_{*}] \subset H_{*}.
\]
The filtration $H_{*}^{(i)}$ on $H_{*}$ gives rise to a tower of morphisms
of simplicial affine group schemes
\[
\xymatrix@1{
H_{*} \ar[r] & \cdots \ar[r] & H_{*}/H_{*}^{(i)} \ar[r] &
H_{*}/H_{*}^{(i-1)} \ar[r] & \cdots \ar[r] & H_{*}/H_{*}^{(1)}.}
\]
Passing to classifying stacks gives a tower of morphisms
of pointed simplicial presheaves
\[
\xymatrix@1{
BH_{*} \ar[r] & \dots \ar[r] & B\left(H_{*}/H_{*}^{(i)}\right) \ar[r] &
B\left(H_{*}/H_{*}^{(i-1)}\right) \ar[r] & \dots \ar[r] &
B\left(H_{*}/H_{*}^{(1)}\right).} 
\]
Combining all these we get a functor
\begin{equation} \label{eq:towerfunctor}
\xymatrix@R-2pc{
\sff{sGAff}  \ar[r] & \sff{SPr}_{*}(\mathbb{C})^{\mathbb{T}} \\
G_{*} \ar@{|->}[r] & \left\{B\left(H_{*}/H_{*}^{(i)}\right)\right\}_{i}.
}
\end{equation}

\begin{lem}\label{l30}
The functor \eqref{eq:towerfunctor}
preserves equivalences and induces a well defined functor
\[
\op{Ho}(\sff{sGAff}) \longrightarrow
\op{Ho}(\sff{SPr}_{*}(\mathbb{C})^{\mathbb{T}}). 
\]
\end{lem}
{\bf Proof:} First of all, by construction each $H_{n}$ is a
unipotent group scheme. Since sheaves represented by affine and
unipotent group schemes are stable under passage to finite limits and finite
colimits, we conclude that the homotopy sheaves $\pi_{j}(BH_{*},*)$ are
unipotent. Therefore, Theorem \ref{t0} implies that the pointed
stacks $BH_{*}$ and $BH_{*}'$ are affine stacks.

Now suppose $G_{*} \longrightarrow G_{*}'$ is an equivalence in the
model category $\sff{sGAff}$. First we will check that the induced
morphism $BH_{*} \longrightarrow BH_{*}'$ is an isomorphism of pointed
stacks. Indeed, since these stacks are affine stacks, it is enough to
show that the morphism $BH_{*} \longrightarrow BH_{*}'$ induces an
isomorphism on the cohomology with coefficient in $\mathbb{G}_{a}$
(see \cite[Theorem~2.2.9]{t1}).  As $G_{*} \longrightarrow G_{*}'$ is an
equivalence, it induces an isomorphism on the Hochschild cohomologies,
and in particular it induces a quasi-isomorphism $\mathcal{O}(BG_{*})
\longrightarrow \mathcal{O}(BG_{*}')$.  Similarly the natural
morphism 
\[
\mathcal{O}(BLG_{*}) \longrightarrow \mathcal{O}(BLG_{*}')
\]
is a quasi-isomorphism. Due to the fact that 
\[
BLG_{*}\simeq \op{hocolim} BLG_{n},
\]
the problem reduces to showing that for each integer $n$ the
induced morphism on Hochschild cohomology
\begin{equation} \label{eq:hhmorphism}
H^{*}(LG_{n},\mathcal{O}) \longrightarrow H^{*}(H_{n},\mathcal{O})
\end{equation}
is an isomorphism.
But  $LG_{n}$ is a free affine group scheme, and 
$H_{n}$ are free unipotent affine group schemes, and therefore
\[
H^{i}(LG_{n},\mathcal{O})=H^{i}(H_{n},\mathcal{O})=0 \; \quad \forall \;
i>1.
\]
Furthermore $H^{1}(K,\mathcal{O})\simeq
H^{1}(K^{\op{uni}},\mathcal{O})$ for any 
affine group scheme $K$ and hence \eqref{eq:hhmorphism} is an isomorphism. 

This implies that $BH_{*} \longrightarrow BH_{*}'$ and so, 
Lemma \ref{l30} reduces to  the following result:

\begin{lem}\label{l31}
Let $H_{*}$ and $K_{*}$ be objects in $\sff{sGAff}$, such that each
$H_{n}$ and $K_{n}$ is a free unipotent affine group scheme.  Let $f :
H_{*}\longrightarrow K_{*}$ be an equivalence in $\sff{sGAff}$. Then,
for all $i>0$, the induced morphism
\[
BH_{*}/H_{*}^{(i)} \longrightarrow BK_{*}/K_{*}^{(i)}
\]
is an equivalence of pointed stacks.
\end{lem}
{\bf Proof:} By induction it is enough to check that
each morphism
\[
BH_{*}^{(i)}/H_{*}^{(i+1)} \longrightarrow BK_{*}^{(i)}/K_{*}^{(i+1)}
\]
is an equivalence.
To see this, we consider
the induced morphism
\[
\prod_{i}H_{*}^{(i)}/H_{*}^{(i+1)} \longrightarrow
\prod_{i}K_{*}^{(i)}/K_{*}^{(i+1)}.
\]
This morphism can be seen as a
morphism of simplicial filtered linearly compact Lie algebras, where
the filtration is given by the product decomposition,
\[
\dots \subset \prod_{i>i_{0}+1} \subset \prod_{i>i_{0}} \subset \dots
\prod_{i}, 
\]
and the Lie algebra structure is given by the commutator bracket.
Set
$H_{*}^{\op{ab}}:=H_{*}/H_{*}^{(1)}$ and
$K_{*}^{\op{ab}}:=K_{*}/K_{*}^{(1)}$. 
One has a commutative diagram of simplicial linearly compact vector spaces
\[
\xymatrix{ \prod_{i}H_{*}^{(i)}/H_{*}^{(i+1)} \ar[r] &
\prod_{i}K_{*}^{(i)}/K_{*}^{(i+1)} \\ H_{*}^{\op{ab}}\ar[u] \ar[r] &
K_{*}^{\op{ab}}.\ar[u]}
\] 
Furthermore, as each $H_{n}$ and $K_{n}$ is free,
the vertical morphisms induce isomorphisms
\[
L(H_{*}^{\op{ab}}) \simeq \prod_{i}H_{*}^{(i)}/H_{*}^{(i+1)} \qquad
L(K_{*}^{\op{ab}}) \simeq \prod_{i}K_{*}^{(i)}/K_{*}^{(i+1)},
\]
where $L(V)$ denotes the free simplicial linearly compact Lie algebra
generated by a simplicial linearly compact vector space $V$. Moreover,
these isomorphisms are compatible with the natural filtrations on
$L(H_{*}^{\op{ab}})$ and $L(K_{*}^{\op{ab}})$ defined by the iterated
brackets. Hence, it suffices to prove that the induced morphism of
linear affine stacks $BH_{*}^{\op{ab}} \longrightarrow BK_{*}^{\op{ab}}$ is an
isomorphism. But, this follows  immediately by observing  that
$BH_{*}^{\op{ab}}$ is the homology type of $BH_{*}$, and $BK_{*}^{\op{ab}}$ is
the homology type of $BK_{*}$. Thus
\[
BH_{*}^{\op{ab}}\simeq \mathcal{H}(BH_{*},\mathcal{O}) \simeq
\mathcal{H}(BK_{*},\mathcal{O}) \simeq
BK_{*}^{\op{ab}}
\]
and so the lemma is proven.
 \hfill $\Box$

This completes the proof of Lemma~\ref{l30} \ \hfill $\Box$

\

\bigskip

Lemma \ref{l30}, combined with the equivalence between the homotopy
category $\op{Ho}(\sff{SHT}_{*})$ of pointed schematic homotopy types and
$\op{Ho}(\sff{sGAff})$, yield a functor
\[
W^{(*)} : \op{Ho}(\sff{SHT}_{*}) \longrightarrow
    \op{Ho}(\sff{SPr}_{*}(\mathbb{C})^{\mathbb{T}}).
\]
By restricting this functor to the full sub-category of pointed and connected
affine stacks one gets the {\em weight tower functor}
\[
W^{(*)} : \op{Ho}(\sff{AS}^{c}_{*}) \longrightarrow
\op{Ho}(\sff{SPr}_{*}(\mathbb{C})^{\mathbb{T}}),
\]
from the homotopy category of pointed and connected affine stacks to the
homotopy category of towers of pointed stacks.

\

\medskip

\begin{center} \textit{The weight tower of a schematic homotopy type}
\end{center} 

\

\medskip

Let now $F$ be a pointed schematic homotopy type. Consider
the natural projection
\[
F \longrightarrow K(\pi_{1}(F,*)^{\op{red}},1),
\]
where $\pi_{1}(F)^{\op{red}}$ is the maximal reductive quotient of
$\pi_{1}(F,*)$. Let $F^{0}$ denote the homotopy fiber of this
morphism.  By theorems \ref{t0} and \ref{t1} the stack $F^{0}$ is a
pointed and connected affine stack.  We can therefore consider its 
weight tower $W^{(*)}F^{0}$:

\begin{df}\label{d30}
For a schematic homotopy type $F$, the weight tower of $F$ is the
object $W^{(*)}F^{0} \in \op{Ho}(\sff{SPr}_{*}(\mathbb{C})^{\mathbb{T}})$
defined above. For any $i\geq 0$, the $i$-th graded piece of the
weight tower of $F$, denoted by $\op{Gr}_{W}^{(i)}F^{0}$, is defined to be
the homotopy fiber of the morphism $W^{(i+1)}F^{0} \longrightarrow
W^{(i)}F^{0}$.
\end{df}

\

\noindent
{\bf Caution:} The construction $F \mapsto
W^{(*)}F^{0}$ is not fully functorial in $F$, simply because the
assignment $F
\mapsto \pi_{1}(F,*)^{\op{red}}$ is not functorial in $F$. However, $F
\mapsto \pi_{1}(F,*)^{\op{red}}$ and hence $F \mapsto
W^{(*)}F^{0}$ are functorial with respect to \emph{reductive morphisms}.

\

\begin{df}\label{d29}
A morphism $f : F \longrightarrow G$ of pointed schematic homotopy types
is \emph{reductive} if the image of the induced morphism
$\pi_{1}(F,*) \longrightarrow \pi_{1}(G,*)^{\op{red}}$
is a reductive affine group scheme.
\end{df}

\

\noindent
Equivalently, the morphism $f$ is reductive if and only if
the induced functor from the category of linear representations
of $\pi_{1}(G,*)$ to the category of linear representations
of $\pi_{1}(F,*)$ preserves semi-simple objects. In particular, any
isomorphisms of pointed schematic homotopy types is a reductive morphism.

Clearly, the construction $F \mapsto W^{(*)}F^{0}$ induces a functor
\[
W^{(*)} : \op{Ho}(\sff{SHT}_{*})^{\op{red}} \longrightarrow
\op{Ho}(\sff{SPr}_{*}(\mathbb{C})^{\mathbb{T}}), 
\]
from the category $\op{Ho}(\sff{SHT}_{*})^{\op{red}}$ of pointed
schematic homotopy types and reductive morphisms to the homotopy
category of towers of pointed stacks.

\begin{rem}\label{r31}
If $f : (X,x) \longrightarrow (Y,y)$ is a morphism of pointed smooth
projective complex manifolds, then the induced morphism
$(X^{\op{top}}\otimes\mathbb{C},x)^{\op{sch}} \longrightarrow
(Y^{\op{top}}\otimes\mathbb{C},y)^{\op{sch}}$ is a reductive
morphism. Indeed, we can factor every such $f$ as the composition $X \to
\mathbb{P}^{N}\times Y \to Y$ of a closed immersion and a
projection. The projection induces an isomorphism on fundamental
groups and so is obviously a reductive morphism. For closed immersions
the reductivity follows from the pluriharmonicity of the equivariant
harmonic map associated to a reductive local system, see
\cite[Proposition~2.2]{s3}. 
\end{rem}

\

\medskip

\noindent
By construction, there exists a natural morphism of towers of pointed stacks
\[
F^{0} \longrightarrow W^{(*)}F^{0},
\]
corresponding to the projections
\[
LG_{*}\longrightarrow H_{*} \longrightarrow H_{*}/H_{*}^{(i)}
\]
(here $F^{0}=BG_{*}\simeq BLG_{*}$ is considered as the constant tower).
By adjunction, this morphism correspond to a well defined morphism of
pointed stacks
\[
F^{0} \longrightarrow W^{(\infty)}:=\op{holim}_{i \in \mathbb{T}}
\displaylimits W^{(i)}F^{0}.
\]
The essential properties of this morphism are summarized in the
following propodition.

\begin{prop}\label{p21} Let $F$ be a pointed schematic homotopy type.
\begin{enumerate}
\item[(1)] For any $i$, the pointed stack $W^{(i)}F^{0}$ is a pointed and
  connected affine stack. 
\item[(2)] The natural morphism
\[
F^{0} \longrightarrow W^{(\infty)}:=\op{holim}_{i\in
  \mathbb{T}}\displaylimits W^{(i)}F^{0}
\]
is an isomorphism of pointed stacks.
\item[(3)] For any $i\geq 0$,
$\op{Gr}_{W}^{(i)}F^{0}$ is the underlying stack of a
linear schematic homotopy type.
\item[(4)] There is an isomorphism of affine stacks
\[
\op{Gr}_{W}^{(0)}F^{0}\simeq \mathcal{H}(F^{0},\mathcal{O}).
\]
\item[(5)] Let $L$ be  the (linearly compact) free
Lie algebra generated by the (linearly compact) vector space
$H_{\bullet >0}(F^{0},\mathcal{O})$. For any $p\geq 1$ and $q\geq 2$,
let $L_{p,q}$ 
be the (closed) sub-vector space of $L$ consisting of
elements generated by all brackets $[x_{1},[x_{2},\dots,[x_{q}]\dots]$,
with $x_{j}\in H_{d_{j}}(F^{0},\mathcal{O})$ and
$\sum d_{j}=p$.

Then, for any $i\geq 1$ and any $p\geq 1$ there is an isomorphism
\[
\pi_{p}(\op{Gr}_{W}^{(i)}F^{0},*) \simeq L_{p+i,i+1}.
\]
\end{enumerate}
\end{prop}
{\bf Proof:} We already proved items (1), (3), (4) and (5) in the
process of proving lemmas \ref{l30} and \ref{l31}. 
For the proof of (2) we will keep the same notation as for the
construction of $W^{(*)}F^{0}$ (i.e. $F^{0}=BG_{*}$ and
$H_{*}=LG_{*}^{\op{uni}}$). First of all, as we saw in the proof of
lemma \ref{l30} the natural morphism
\[
F^{0}=BG_{*} \longrightarrow BH_{*}
\]
induces an isomorphism on cohomology with coefficients in
$\mathbb{G}_{a}$.  Since $F^{0}$ and $BH_{*}$ are affine stacks this
implies that $F^{0}\simeq BH_{*}$.

Furtermore, each $W^{(i)}F^{0}$ is an affine stack and in
particular a pointed 
schematic homotopy type. Therefore, in order to prove that
\[
F^{0}\simeq BH_{*} \longrightarrow \op{holim}_{i \in
  \mathbb{T}}\displaylimits W^{(i)}F^{0}
\]
is an isomorphism, it is enough to prove that the natural morphism of
co-simplicial Hopf algebras
\begin{equation} \label{eq:even_iso}
\op{colim}_{i \in \mathbb{T}}\displaylimits
\mathcal{O}(H_{*}/H_{*}^{(i)}) \longrightarrow \mathcal{O}(H_{*})
\end{equation}
is a quasi-isomorphism. However $H_{*}$ is unipotent, and so 
$\cap H_{*}^{(i)}=\{ * \}$. Thus \eqref{eq:even_iso} is even an
isomorphism.  The proposition is proven. \hfill $\Box$ 

\

\bigskip

\begin{center} \textit{Equivariant weight filtration} \end{center}

\

\medskip

All of the constructions presented above have enough built-in 
functoriality to be compatible with the extra structure
of a discrete group action. Indeed, suppose $F \in
\op{Ho}(\sff{SPr}_{*}(\mathbb{C})^{\Gamma})$  is a $\Gamma$-equivariant
schematic homotopy type, where $\Gamma$ is a discrete group.
The action of $\Gamma$ can be understood as morphism of simplicial monoids
\[
\Gamma \longrightarrow \underline{\op{End}}(F).
\]
Without a loss of generality we can assume that $F$ is fibrant as a
pointed simplicial presheaf, and therefore the last morphism can be
written as morphism of simplicial monoids
\[
\Gamma \longrightarrow \mathbb{R}\underline{\op{End}}(F),
\]
which is well defined in the homotopy category of simplicial monoids.

Let $G_{*}$ be a fibrant object in $\sff{sGAff}$ such that $F^{0}\simeq
BG_{*}$. According to
\cite[Theorem 3.5]{small}, we have a natural isomorphism of
simplicial monoids
\[
\Gamma \longrightarrow \mathbb{R}\underline{\op{End}}_{\sff{sGAff}}(G_{*}).
\]
Let $\widetilde{\Gamma}$ be a cofibrant replacement of $\Gamma$ in the model
category of simplicial monoids. The above morphism can then be represented as
an actual morphism of simplicial monoids
\[
\widetilde{\Gamma} \longrightarrow
\underline{\op{End}}_{\sff{sGAff}}(G_{*}),
\]
giving rise to a $\widetilde{\Gamma}$-action on $G_{*}$. Applying the
constructions 
$G_{*} \mapsto LG_{*} \mapsto LG_{*}^{\op{uni}} \mapsto W^{(*)}F$
preserves the 
$\widetilde{\Gamma}$-action, and gives rise to a 
$\widetilde{\Gamma}$-object in the model category
$\sff{SPr}_{*}(\mathbb{C})^{\mathbb{T}}$. 
This yields a natural functor
\[
W^{(*)} : \op{Ho}(\sff{SHT}_{*}^{\Gamma})^{\op{red}} \longrightarrow
\op{Ho}((\sff{SPr}_{*}(\mathbb{C})^{\mathbb{T}})^{\widetilde{\Gamma}}),
\]
from the homotopy category of $\Gamma$-equivariant pointed schematic
homotopy types and reductive morphisms to the homotopy category of
$\Gamma$-equivariant towers. This functor is a $\Gamma$-equivariant
refinement of the functor $W^{(*)}$ constructed before. Finally, it is
well known that the natural functor
\[
\op{Ho}((\sff{SPr}_{*}(\mathbb{C})^{\mathbb{T}})^{\Gamma})\longrightarrow
\op{Ho}((\sff{SPr}_{*}(\mathbb{C})^{\mathbb{T}})^{\widetilde{\Gamma}})
\]
is an equivalence of categories (see \cite{dk}). Therefore, we get a functor
\[
W^{(*)} : \op{Ho}(\sff{SHT}_{*}^{\Gamma})^{\op{red}} \longrightarrow
\op{Ho}((\sff{SPr}_{*}(\mathbb{C})^{\mathbb{T}})^{\Gamma}),
\]
which is a $\Gamma$-equivariant version of our previous construction. 

\begin{rem}\label{r30}
In the construction of the weight tower of a pointed schematic
homotopy type $F$, we ignored the natural action of the group
$\pi_{1}(F,*)^{\op{red}}$ on $F^{0}$. Yet another refinement of the
construction should take this action into account, and would produce a
tower of $\pi_{1}(F,*)^{\op{red}}$-equivariant pointed stacks. We did
not investigate this issue carefully but strongly believe that such a
refinement should exists. In particular, the schematic homotopy type
$F$ should be reconstructed from its
$\pi_{1}(F,*)^{\op{red}}$-equivariant tower by the formula
\[
F\simeq \left[\left. \op{holim}_{i \in \mathbb{T}}\displaylimits
 W^{(i)}F^{0}\right/\pi_{1}(F,*)^{\op{red}}\right].
\]
\end{rem}

\subsection{Real structures}

For a field extension $L/K$, one has a
base change functor
$$-\otimes_{K}L : \op{Ho}(\sff{SPr}_{*}(K)) \longrightarrow
\op{Ho}(\sff{SPr}_{*}(L)),$$
which is right adjoint to the forgetful functor
$$\op{Ho}(\sff{SPr}_{*}(L)) \longrightarrow
\op{Ho}(\sff{SPr}_{*}(K)).$$

The functor $-\otimes_{K}L$ preserves pointed schematic homotopy
types, and so  
the universal property of the schematization implies that for any
pointed and connected simplicial set $(X,x)$ there exists a natural morphism
\[
(X\otimes L,x)^{\op{sch}} \longrightarrow (X\otimes
K,x)^{\op{sch}}\otimes_{K} L.
\] 
This morphism is in general not an
isomorphism, but it is known to be an isomorphism when $L/K$ is a
finite extension.

\begin{prop}\label{p20}
Let $L/K$ be a finite extension of fields and let $X$ be a pointed and
connected simplicial set. The natural morphism 
\[
(X\otimes L,x)^{\op{sch}} \longrightarrow (X\otimes
K,x)^{\op{sch}}\otimes_{K}L 
\]
is an isomorphism of pointed stacks.
\end{prop}
{\bf Proof:} For any discrete group $\Gamma$,  the
natural morphism $\Gamma^{\op{alg},L} \longrightarrow
\Gamma^{\op{alg},K}\otimes_{K}L$ is an isomorphism of affine group
schemes (here $\Gamma^{\op{alg},k}$ denotes the pro-algebraic
completion of $\Gamma$ over the field $k$). As the schematization of
$X$ over a field $k$ can be written (see \cite{small}) as
$BG_{*}^{\op{alg},k}$ for a simplicial group $G_{*}$
this implies the proposition. \hfill $\Box$

\

\medskip

\noindent
We will use roposition \ref{p20} to
endow the stack $(X\otimes \mathbb{C},x)^{\op{sch}}$ with a natural
real structure 
\[
(X\otimes \mathbb{C},x)^{\op{sch}} \longrightarrow (X\otimes
\mathbb{R},x)^{\op{sch}}\otimes_{\mathbb{R}} 
\mathbb{C}. 
\]
In the case where $X$ is a smooth compact manifold, this real structure can be
directly seen at the level of the equivariant  algebra
of differential forms
$(A^{\bullet}_{DR}(\mathcal{O}(G_{X})),\nabla)$. Indeed, the $\opi{Ind}$-flat
bundle $\mathcal{O}(G_{X})$ has 
a natural real form given as the Tannakian dual
of the category of real flat bundles over $X$. This induces a natural
real structure on the de Rham complex
$(A^{\bullet}_{DR}(\mathcal{O}(G_{X})),\nabla)$ 
which is the real structure of $(X\otimes \mathbb{C},x)^{\op{sch}}$
discussed above. 

Of course, for a pointed schematic homotopy type $F$ over
$\mathbb{R}$, one can construct a real weight tower $W^{(*)}F^{0} \in
\op{Ho}(\sff{SPr}_{*}(\mathbb{R})^{\mathbb{T}})$, compatible with the base
change from $\mathbb{R}$ to $\mathbb{C}$.  This can be seen as
follows.

Let $F$ be a pointed schematic homotopy type over $\mathbb{R}$, and
let $F^{0}$ be the homotopy fiber of the projection $F \longrightarrow
K(\pi_{1}(F,*)^{\op{red}},1)$, where $\pi_{1}(F,*)^{\op{red}}$ is the
real maximal reductive quotient of $\pi_{1}(F,*)$. Clearly, one has
$\pi_{1}(F,*)^{\op{red}}\otimes_{\mathbb{R}}\mathbb{C}\simeq
(\pi_{1}(F,*)\otimes_{\mathbb{R}}\mathbb{C})^{\op{red}}$, and
therefore
\[
F^{0}\otimes_{\mathbb{R}}\mathbb{C} \simeq
(F\otimes_{\mathbb{R}}\mathbb{C})^{0}. 
\]
Hence $F^{0}\simeq BG_{*}$, for a simplicial real affine group scheme
$G_{*}$ (see \cite[Theorem 3.5]{small}).  Consider now the functor
$H_{*} \mapsto H_{*}(\mathbb{C})$, from the category of simplicial
real affine group schemes to the category $\mathbb{Z}/2-\sff{SSet}$ of
simplicial sets with an action of $\mathbb{Z}/2$. This functor is
exact, conservative and has a left adjoint. It can therefore be used
to construct a resolution $LG_{*} \longrightarrow G_{*}$, whichin
turns gives rise to the standard free resolution of
$G_{*}\otimes_{\mathbb{R}}\mathbb{C}$ used in the definition of the
weight tower of $F\otimes_{\mathbb{R}}\otimes \mathbb{C}$. Finally, as
the lower central serie and the construction $H \mapsto H^{\op{uni}}$ are
compatible with base change from $\mathbb{R}$ to $\mathbb{C}$, this
shows that the weight tower
$W^{(*)}(F\otimes_{\mathbb{R}}\mathbb{C}^{0})$ has a natural real
structure $W^{(*)}(F^{0})\in
\op{Ho}(\sff{SPr}_{*}(\mathbb{R})^{\mathbb{T}})$ given by
\[
W^{(i)}(F^{0}):=B(LH_{*}/LH_{*}^{(i)}) \qquad H_{*}=LG_{*}^{\op{uni}}.
\]
This construction gives a functor
\[
W^{(*)} : \op{Ho}(\sff{SHT}_{*}(\mathbb{R}))^{\op{red}}
\longrightarrow \op{Ho}(\sff{SPr}_{*}(\mathbb{R})^{\mathbb{T}})
\]
from the homotopy category of pointed schematic homotopy types over
$\mathbb{R}$ and reductive morphisms, to the homotopy category of
towers of real pointed stacks. The compatibility with the base change
from $\mathbb{R}$ to $\mathbb{C}$ is an isomorphism of towers of
pointed stacks (functorial with respect to reductive morphisms)
\[
(W^{(*)}F^{0})\otimes_{\mathbb{R}}\mathbb{C} \simeq
W^{(*)}((F\otimes_{\mathbb{R}}\mathbb{C})^{0}).
\]
\

\subsection{The spectral sequence of a tower}

Let $\{F_{i}\} \in
\op{Ho}(\sff{SPr}_{*}(\mathbb{C})^{\mathbb{T}})$ be a tower of pointed
stacks: 
\[
\xymatrix{ \dots \ar[r] & F_{i} \ar[r] & F_{i-1} \ar[r] & \dots \ar[r]
  & F_{1} \ar[r] & F_{0}=\bullet,}
\]
and let $W_{i}$ denote the homotopy fiber of the morphism $F_{i+1}
\longrightarrow F_{i}$. The associated 
long exact sequences of  homotopy groups give rise to an exact couple
of sheaves 
\[
\xymatrix{
\dots \ar[r] & \pi_{*}(F_{i}) \ar[r] & \ar[ld] \pi_{*}(F_{i-1}) \ar[r]
& \ar[ld] \dots \ar[r] & \pi_{*}(F_{1}) \\ 
& \pi_{*}(W_{i-1}) \ar[u] & \pi_{*}(W_{i-2}) \ar[u] & \dots &
\pi_{*}(W_{0}) \ar[u],}
\]
and hence defines a spectra sequence $\{E_{r}^{p,q}(F_{*})\}$ of
sheaves of groups 
(abelian when $q-p>1$). To explicate, consider first the groups
\[
\begin{split}
Z^{p,q}_{r} & = \op{Ker} \left(
\pi_{q-p}(W_{p-1}) \longrightarrow
\frac{\pi_{q-p}(F_{p})}{\op{Im} \left(\pi_{q-p}(F_{p+r-1}) \rightarrow
  \pi_{q-p}(F_{p})\right)}\right) \\ \vspace{0.2in}
B^{p,q}_{r} & = \op{Ker}\left(
\pi_{q-p+1}(F_{p-1}) \longrightarrow \pi_{q-p+1}(F_{p-r})\right).
\end{split}
\]
The boundary operator $\pi_{q-p+1}(F_{p-1}) \longrightarrow
\pi_{q-p}(W_{p-1})$ 
induces a morphism 
\[
\partial : B^{p,q}_{r}\longrightarrow
Z^{p,q}_{r},
\] 
and we set
\[
E_{r}^{p,q}=\frac{Z^{p,q}_{r}}{\partial B^{p,q}_{r}}.
\]
The differential
\[
d_{r} : E_{r}^{p,q} \longrightarrow E_{r}^{p+r,q+r-1}
\]
is given by the composition
\[
E_{r}^{p,q} \longrightarrow \op{Im} \left(\pi_{q-p}(F_{p+r-1})
\rightarrow \pi_{q-p}(F_{p})\right) 
\longrightarrow E_{r}^{p+r,q+r-1}.
\]
We refer to \cite[VI \S 2]{gj} for more details on this spectral
sequence.  

The convergence of this spectral sequence can be quite subtle in
general. However it is known that under certain conditions it
converges and its limit computes the groups $\pi_{*}(\op{holim}
F_{i})$. For our purposes the following simple case of the
\textit{complete convergence lemma} of Bousfield and Kan will suffice.

\begin{lem}\label{l32}
Let $\{F_{i}\} \in \op{Ho}(\sff{SPr}_{*}(\mathbb{C})^{\mathbb{T}})$ be a
tower of pointed  
and connected affine stacks, and let $F:= \op{holim}  F_{i}$.
Let $\{E_{r}^{p,q}(F_{*})\}$ be the associated
spectral sequence in homotopy, and assume that there is an integer
$N$ such that $d_{r}=0$ for any $r\geq N$. Then the following two
conditions are satisfied.
\begin{enumerate}
\item The limiting tableu of $\{E_{r}^{p,q}(F_{*})\}$ is given by
\[
E_{\infty}^{p,q}\simeq \frac{\op{Ker}\left(\pi_{q-p}(F)\rightarrow
  \pi_{q-p}(F_{p-1})\right)} 
{\op{Ker}\left(\pi_{q-p}(F)\rightarrow \pi_{q-p}(F_{p})\right)}
\]
\item The natural morphism
\[
\pi_{*}(F) \longrightarrow \lim \pi_{*}(F_{i})
\]
is an isomorphism.
\end{enumerate}
\end{lem}
{\bf Proof:} Consider the global section functor $\mathbb{R}\Gamma :
\op{Ho}(\sff{SPr}_{*}(\mathbb{C}))\longrightarrow
\op{Ho}(\sff{SSet}_{*})$.  The fact that $H^{i}(\op{Spec}\, A,H)=0$ for any
$i>0$ and any affine unipotent group scheme $H$ implies that
\linebreak 
$\mathbb{R}\Gamma(K(H,n))\simeq K(H(\mathbb{C}),n)$, for any unipotent
affine group scheme $H$. Postnikov induction combined with
\cite[Proposition 1.2.2]{t1} and Theorem \ref{t0} implies that for any
pointed and connected affine stack $F$ the natural morphism
\[
\pi_{i}(\mathbb{R}\Gamma(F),*) \longrightarrow
\pi_{i}(F,*)(\mathbb{C})
\]
is always an isomorphism. In particular, $\mathbb{R}\Gamma$ is a
conservative functor which commutes with taking homotopy
groups. Recall also that the functor $H \mapsto H(\mathbb{C})$ is an
exact and conservative functor from the category of sheaves of groups
represented by affine group schemes to the category of groups.

Now, since each stack $F_{i}$ is a pointed and connected affine stack,
the spectral sequence $\{E_{r}^{p,q}(F_{*})\}$ consists of affine
unipotent group schemes. Therefore, the spectral sequence
$\{E_{r}^{p,q}(F_{*})\}(\mathbb{C})$ is the spectral sequence of the
tower of pointed simplicial sets $\{F_{i}(\mathbb{C})\}$ described in
\cite[VI \S 2]{gj}.  The proposition now follows from
\cite[VI Corollary 2.22]{gj}. \hfill $\Box$

\

\medskip

By construction the spectral sequence $\{E_{r}^{p,q}(F_{*})\}$ is
functorial in the tower $\{F_{i}\}$. Therefore, if a discrete group
acts on the tower of pointed stacks $\{F_{i}\}$, then it also acts on
the spectral sequence $\{E_{r}^{p,q}(F_{*})\}$. Similarly, if 
the tower $\{F_{i}\}$ is defined over $\mathbb{R}$, then so is the
spectral sequence $\{E_{r}^{p,q}(F_{*})\}$.

\subsection{Purity and degeneration of the weight spectral sequence} 

Let $F$ be a pointed schematic homotopy type.  The weight tower
$W^{(*)}F^{0}$ gives rise to the \textit{weight spectral sequence}
$\{E_{r}^{p,q}(W^{(*)}F^{0})\}$ of $F$. By construction this spectral
sequence is functorial in $F$ with respect to reductive morphisms.

\begin{thm}\label{tpurity}
Let $(X,x)$ be a pointed complex projective  manifold and
$F=(X^{\op{top}}\otimes\mathbb{C},x)^{\op{sch}}$ be its
schematization. Then the 
weight spectral sequence
$\{E_{r}^{p,q}(W^{(*)}F^{0})\}$ of $F$ degenerates at $E_{2}$
(i.e. $d_{r}=0$ for all $r\geq 2$).
\end{thm}
{\bf Proof:} The proof is based on the standard purity argument. 

Consider a category $\mathcal{C}$ defined as follows. The objects of
$\mathcal{C}$ are linearly compact $\mathbb{R}$-vector spaces $V$
(i.e. abelian unipotent affine group schemes over $\op{Spec}\,
\mathbb{R}$) together with an action of the discrete group
$\mathbb{C}^{\times\delta}$ on the linearly compact
$\mathbb{C}$-vector space $V\otimes_{\mathbb{R}}\mathbb{C}$. The
morphisms in $\mathcal{C}$ are morphisms of linearly compact
$\mathbb{R}$-vector spaces whose base extension to $\mathbb{C}$
commutes with the action of $\mathbb{C}^{\times}$.  The category
$\mathbb{C}$ is clearly an abelian category. For an object $V$ in
$\mathcal{C}$, the action of an element $\lambda \in
\mathbb{C}^{\times}$ on $v\in V\otimes_{\mathbb{R}} \mathbb{C}$ will
be denoted by 
$\lambda(v)$, while $\lambda\cdot v$ will denote the scaling action of
$\mathbb{C}^{\times}$ coming from the $\mathbb{C}$-vector space
structure.

\begin{df}\label{d33}
An object $V \in \mathcal{C}$ is pure of weight $n \in \mathbb{Z}$ if
for any $\lambda \in U(1) \subset \mathbb{C}^{\times}$ and any
$v \in V\otimes_{\mathbb{R}}\mathcal{C} $ one has
\[
\lambda^{n}\cdot \overline{\lambda(v)}=\lambda(\overline{v}),
\]
where $x\mapsto \overline{x}$ denotes complex conjugation.
\end{df}

Clearly, the full subcategory of objects in $\mathcal{C}$ which are
pure of weight $n$ is stable under passing to extensions, sub-objects
and cokernels.  Furthermore, it is clear that there are no non-zero
morphisms between two objects in $\mathcal{C}$ which are pure of
different weights. Finally, the (completed) tensor products of
linearly compact vector spaces, turns $\mathcal{C}$ into a symmetric
monoidal category, and the tensor product of two pure objects of
weights $p$ and $q$ gives a pure object of weight $p+q$.

Going back to the proof of theorem \ref{tpurity}, note that the Hodge
decomposition on $(X^{\op{top}}\otimes\mathbb{C},x)^{\op{sch}}$,
induces an action of $\mathbb{C}^{\times\delta}$ on the tower
$W^{(*)}F^{0}$ and hence on the weight spectral sequence
$\{E_{r}^{p,q}(W^{(*)}F^{0})\}$. Furthermore, the natural real
structure $(X^{\op{top}}\otimes\mathbb{R},x)^{\op{sch}}$ on the
schematic homotopy type $(X^{\op{top}}\otimes\mathbb{C},x)^{\op{sch}}$
gives rise to a natural real structure on the weight spectral sequence
$\{E_{r}^{p,q}(W^{(*)}F^{0})\}$ of $F$. This shows that the spectral
sequence $\{E_{r}^{p,q}(W^{(*)}F^{0})\}$ can be viewed as a spectral
sequence in the abelian category $\mathcal{C}$ described above.

Now, in order to check that $d_{r} : E_{r}^{p,q} \longrightarrow
E_{r}^{p+r,q+r-1}$ vanish for all $r\geq 2$, it suffices to check that
$E_{r}^{p,q}$, as an object of $\mathcal{C}$, is pure of weight
$-q$. This, in turn, will follow if we can show that the object
$E_{1}^{p,q} \in \mathcal{C}$ is pure of weight $-q$. By
Proposition~\ref{p21}(5) it is enough to check that
$H_{q}(F^{0},\mathcal{O})$ is pure of weight $-q$ (recall that
$L_{p,q}$ can be identified with a sub-space of
$H_{d_{1}}(F^{0},\mathcal{O})\otimes \dots \otimes
H_{d_{q}}(F^{0},\mathcal{O})$ with $\sum d_{j}=p$).

Finally, $H_{q}(F^{0},\mathcal{O})$ is the dual of the vector space
$H^{q}(X,\mathcal{O}(G_{X}))$, where $G_{X}$ is the pro-reductive completion
of the group $\pi_{1}(X,x)$. As an $\opi{Ind}$-local system on
$X$, $\mathcal{O}(G_{X})$ is isomorphic to
$\oplus_{L}L^{\dim\, L}$, where the sum is over the set of isomorphism
classes of simple $\mathbb{C}$-linear representations of $\pi_{1}(X,x)$.
Using the non-abelian Hodge correspondence this $\opi{Ind}$-local system
correspond to the $\opi{Ind}$-Higgs bundle
$\oplus_{L}V^{\dim\, V}$ where $V$ runs trough the set of isomorphism
classes of stable Higgs bundles of degree $0$. The space
$H^{q}(X,\mathcal{O}(G_{X}))$ can therefore 
be described as
\[
H^{q}(X,\mathcal{O}(G_{X}))\simeq \bigoplus_{(V,\theta)}
H^{q}(A^{\bullet}_{Dol}(V^{dim\, V},D'')),
\]
where the right hand side is the Dolbeault cohomology as recalled in
Section \ref{ss-review}. Now, if
$(V,D'')=(V,\overline{\partial}+\theta)$ is a stable Higgs 
bundle of degree $0$ corresponding to a local system $L$ on $X$, the
Higgs bundle $(\overline{V},\overline{\partial}-\theta)$ corresponds to the
local system $\overline{L}$ complex conjugate of $L$ (see
\cite{s2}). Therefore, 
the complex conjugation acts on $H^{q}(X,\mathcal{O}(G_{X}))$ by sending
a differential form $v \in A^{i,j}_{Dol}(V,\overline{\partial}+\theta)$ to
the form $\overline{v} \in
A^{j,i}_{Dol}(\overline{V},\overline{\partial}-\theta)$. 
On the other hand, $\mathbb{C}^{\times}$ acts on
$H^{q}(X,\mathcal{O}(G_{X}))$ by sending a differential form $v \in
A^{i,j}_{Dol}(V,\overline{\partial}+\theta)$ to 
the form $\lambda^{j}\cdot v \in
A^{i,j}_{Dol}(V,\overline{\partial}+\lambda\cdot \theta)$. 
Therefore, for $\lambda\in \mathbb{C}^{\times}$, one has
\[
\lambda^{q}\cdot \overline{\lambda(v)}=\lambda(\overline{v}) \in
A^{j,i}_{Dol}(\overline{V},\overline{\partial}-\lambda\cdot \theta).
\]
This shows that for any $v\in H^{q}(X,\mathcal{O}(G_{X}))$ and any
$\lambda \in \mathbb{C}^{\times}$, one has
$\lambda^{q}\cdot
\overline{\lambda(v)}=\lambda(\overline{v})$. Dualizing, one gets
that $H_{q}(F^{0},\mathcal{O})\simeq
H^{q}(X,\mathcal{O}(G_{X}))^{\vee}$ is pure 
of weight $-q$. The theorem is proven. \ \hfill $\Box$ 

\

\bigskip

\noindent
Theorem \ref{tpurity} and lemma \ref{l32} have the 
following important corollary.

\begin{cor}\label{c30}
Let $(X,x)$ be a pointed   complex projective manifold.
Let $F=(X^{\op{top}}\otimes\mathbb{C},x)^{\op{sch}}$ be its schematization,
$W^{(*)}F^{0}$ be the corresponding weight tower and set
\[
F^{(p)}_{W}\pi_{q}(F^{0},*):= \op{Ker} \left(\pi_{q}(F^{0},*)
\longrightarrow \pi_{q}(W^{(p)}F^{0},*)\right). 
\]
Then, one has
\begin{enumerate}
\item[(1)] $\pi_{1}(F^{0},*)= \op{Ker} \left( \pi_{1}(F,*) \longrightarrow
\pi_{1}(F,*)^{\op{red}} \right).$
\item[(2)] $\pi_{q}(F^{0},*)\simeq \pi_{q}(F,*), \qquad \forall \; q>1.$
\item[(3)] $\pi_{q}(F^{0},*)\simeq \lim\,
  \pi_{q}(F^{0},*)/F^{(p)}_{W}\pi_{q}(F^{0},*).$ 
\item[(4)]
  $F^{(p-1)}_{W}\pi_{q}(F^{0},*)/F^{(p)}_{W}\pi_{q}(F^{0},*)\simeq 
  E_{\infty}^{p,q+p} 
\simeq E_{2}^{p,q+p}(W^{(*)}F^{0}).$ 
\end{enumerate}
\end{cor}

\

\medskip

\noindent
Corollary \ref{c30} gives a concrete way of computing the homotopy
groups of the schematization $(X^{\op{top}}\otimes \mathbb{C},
x)^{\op{sch}}$ of a projective manifold, while very little is
known of these groups for a general topological space $X$. Indeed, by
Proposition \ref{p21} (5) the
complex $(E_{1}^{p,q},d_{1})$ can be  described explictly. For
instance, the groups  
$E_{1}^{p,q}$ are given by
\[
E_{1}^{p,q}\simeq L_{q-1,p} \subset
\bigoplus_{d_{1}+\dots+d_{p}=q-1}H_{d_{1}}(F^{0},
\mathcal{O})\widehat{\otimes}   
\dots \widehat{\otimes}H_{d_{p}}(F^{0},\mathcal{O}). 
\]
The differential $d_{1}$ is given by the 
(co-)products listed in Proposition \ref{p21} (5), whereas the
vanishing of $d_{2}$ can be interpreted as 
the vanishing of Massey (co)-products (see \cite{arch}). 

\section{Restrictions on homotopy types}

In this section we will give some applications of the existence of a
mixed Hodge structure on the schematization of a projective manifold.

\subsection{An example}

In this section we will use the existence of the Hodge decomposition
on $\sch{X^{\op{top}}}{\mathbb{C}}$ in order to give examples of
homotopy types which are not realizable as homotopy types of smooth
projective varieties. These examples are obtained after defining new
homotopy invariants of a space $X$ in terms of the stack $(X\otimes
\mathbb{C})^{\op{sch}}$ and the action of $\pi_{1}(X,x)$ on the spaces
$\pi_{i}(\sch{X}{\mathbb{C}},x)$. The existence of the Hodge
decomposition implies strong restrictions on these invariants, and it
is relatively easy to find explicit examples of homotopy types
violating these restrictions.

We should note also that our invariants are trivial as soon as one
restricts to the case when the action of $\pi_{1}(X,x)$ on
$\pi_{i}((X\otimes \mathbb{C})^{sch},x)$ is nilpotent. Therefore, it
appears that our examples can not be ruled out by using Hodge theory
on rational homotopy types in the way it is done for example in
\cite{dgms,mo}.

Let us start with a pointed schematic homotopy type $F$ such that
$\pi_{1}(F,*)$ is an affine group scheme. According to
\cite[Proposition $3.2.9$]{t1} the group scheme $\pi_{i}(F,*)$ is an
abelian unipotent affine group scheme for $i>1$. Hence $\pi_{i}(F,*)$
is isomorphic to a (possibly infinite) product of $\mathbb{G}_{a}$'s.

Consider now the
maximal reductive quotient of the affine group scheme $\pi_{1}(F,*)$
\[
\pi_{1}(F,*) \longrightarrow \pi_{1}(F,*)^{\op{red}}.
\]
Using  the Levy decomposition, let us choose a section of this
morphism $s : \pi_{1}(F,*)^{\op{red}} \longrightarrow
\pi_{1}(F,*)$. The morphism $s$ allows us to consider the induced
action of $\pi_{1}(F,*)^{\op{red}}$ on $\pi_{i}(F,*)$.
But, since  $\pi_{i}(F,*)$ is a linearly compact vector space, and
$\pi_{1}(F,*)$ is a reductive affine group
scheme acting on it, there exists a decomposition of $\pi_{i}(F,*)$ as
a (possibly infinite) product
\[
\pi_{i}(F,*)\simeq \prod_{\rho \in
  R(\pi_{1}(F,*))}\pi_{i}(F,*)^{\rho},
\]
where $R(\pi_{1}(F,*))$ is the set of isomorphism classes of finite
dimensional simple linear representations of
$\pi_{1}(F,*)$, and $\pi_{i}(F,*)^{\rho}$ is a product (possibly
infinite) of representations in  the class $\rho$.
Using the fact that the Levy decomposition is unique
up to an inner automorphism one can check that the
set $\{\rho \in R(\pi_{1}(F,*)) \; | \; \pi_{i}(F,*)^{\rho}\neq 0\}$
is independent of the choice of the section $s$.  With this notation
we have the
following:

\begin{df}\label{dsupp}
Let $F$ be a pointed schematic homotopy type.
The subset
\[
\op{Supp}(\pi_{i}(F,*)) = \{\rho \in R(\pi_{1}(F,*)) \; | \;
\pi_{i}(F,*)^{\rho} \neq
0 \}\subset
R(\pi_{1}(F,*))
\]
is called {\em the  support of $\pi_{i}(F,*)$} for every $i > 1$.
\end{df}
Note that for a pointed and connected simplicial set in
$\mathbb{U}$ the supports $\op{Supp}(\pi_{i}((X\otimes
\mathbb{C})^{\op{sch}},x))$ of  $X$ are  homotopy
invariants of $X$.

Suppose now that $X^{\op{top}}$ is the underlying pointed space of a smooth
projective algebraic variety over
$\mathbb{C}$. Then, the stack $\sch{X^{\op{top}}}{\mathbb{C}}$
comes equipped with the Hodge decomposition
defined in the previous  section.  The naturality of the construction
of the Hodge decomposition implies  the following:
\begin{lem}
For any $i>1$, the subset $\op{Supp}(\pi_{i}((X^{\op{top}}\otimes
\mathbb{C})^{\op{sch}},x))$
is invariant under the $\mathbb{C}^{\times }$-action on
$R(\pi_{1}(\sch{X^{\op{top}}}{\mathbb{C}},x))\simeq
R(\pi_{1}(X^{\op{top}},x))$.
\end{lem}

As a  consequence of this lemma we get.

\begin{cor}\label{vhs}
Let $X$ be a pointed smooth projective complex algebraic variety.
\begin{enumerate}
\item If $\rho \in \op{Supp}(\pi_{i}((X^{\op{top}}\otimes
\mathbb{C})^{\op{sch}},x)$ is an isolated point
(for the induced topology of $R(\pi_{1}(X,x))$), then its
corresponding local system on $X$ underlies
a polarizable complex variation of Hodge structure.
\item If $\pi_{i}((X^{\op{top}}\otimes \mathbb{C})^{\op{sch}},x)$ is
  an affine group 
scheme of finite type, then each
simple factor of the semi-simplification of the representation
of $\pi_{1}(X,x)$ to the vector space $\pi_{i}((X^{\op{top}}\otimes
\mathbb{C})^{\op{sch}},x)$ underlies a
polarizable complex variation of Hodge structure on $X$.
\item Suppose that $\pi_{1}(X,x)$ is abelian. Then each isolated
character $\chi \in \op{Supp}(\pi_{i}((X^{\op{top}}\otimes
\mathbb{C})^{\op{sch}},x))$ is unitary.
\end{enumerate}
\end{cor}
{\bf Proof:} Clearly   Parts $(2)$ and $(3)$ of the corollary are direct
consequences of $(1)$.  For the proof of part $(1)$ note that the fact that
the action of $\mathbb{C}^{\times }$ on the
space $R(\pi_{1}(X,x))$ is
continuous yields that $\rho \in
\op{Supp}(\pi_{i}((X^{\op{top}}\otimes \mathbb{C})^{\op{sch}},x))$ is
fixed by the action of $\mathbb{C}^{\times }$. Part $(1)$ of the corollary
now follows from
\cite[Corollary
$4.2$]{s2}. \hfill $\Box$

\

\bigskip

\noindent
The previous corollary shows that the existence of the Hodge
decomposition imposes severe restrictions on the supports
$\op{Supp}(\pi_{i}((X^{\op{top}}\otimes \mathbb{C})^{\op{sch}},x))$,
and therefore on the homotopy type of $X^{\op{top}}$.  The next
theorem will provide a family of examples of homotopy types which are
not realizable as homotopy typer of smooth projective varieties.

Before stating the theorem let us recall that the notion of a good
group (see \cite[$\S 3.4$]{t1}). Recall
that for a group $\Gamma$ we denote by $\Gamma^{\op{alg}}$ its
pro-algebraic completion.

\begin{df}
A group $\Gamma$ (in $\mathbb{U}$) is called algebraically good
(relative to $\mathbb{C}$) if the natural
morphism of pointed stacks
\[
(K(\Gamma,1)\otimes \mathbb{C})^{\op{sch}} \longrightarrow
K(\Gamma^{\op{alg}},1)
\]
is an isomorphism.
\end{df}
\

\noindent
\textit{Remark:} An immediate but important remark is that a group
$\Gamma$ is algebraically good if and only if for any linear
representation of finite dimension $V$ of $\Gamma^{\op{alg}}$ the
natural morphism $\Gamma \longrightarrow \Gamma^{\op{alg}}$ induces an
isomorphism
\[
H^{\bullet}_{H}(\Gamma^{\op{alg}},V) \simeq H^{\bullet}(\Gamma,V),
\]
where $H^{\bullet}_{H}$ denotes Hochschild cohomology of the affine
group scheme $\Gamma^{\op{alg}}$.

\

\medskip

\noindent
Besides finite groups, there are three known classes of
examples of algebraically good 
groups (see \cite[\S 4.5]{small} for proofs and details):
\begin{itemize}
\item Any finitely generated abelian group is an algebraically good
group.
\item Any fundamental group of a compact Riemann surface is
an algebraically good group.
\item A successive extension of a finitely presented free groups is an
algebraically good group.
\end{itemize}
\

\medskip

\noindent
In addition to the goodness property we will often need the following
finiteness condition: 

\begin{df} \label{def-finite}
Let $\Gamma$ be an abstract group. We will say that the
group $\Gamma$ is of type $(\boldsymbol{F})$ (over ${\mathbb C}$) if: 
\begin{itemize}
\item[(a)] For every $n$ and every finite dimensional complex
representation $V$ of $\Gamma$ the group $H^{n}(\Gamma,V)$ is finite
dimensional; 
\item[(b)] $H^{\bullet}(\Gamma,-)$ commutes with inductive limits of
finite dimensional complex representations of $\Gamma$.
\end{itemize}
\end{df}

\begin{rem} \label{rem-finite} {\bf (i)} The condition $(\boldsymbol{F})$ is a
slight relaxation of a standard finiteness condition in combinatorial
group theory. Recall \cite[Chapter~VIII]{brown} that a group
$\Gamma$ is defined to be of type $\fp_{\infty}$ if the trivial
${\mathbb Z}\Gamma$-module ${\mathbb Z}$ admits a resolution by free
${\mathbb Z}\Gamma$-modules of finite type.  It is well known,
\cite[Section~VIII.4]{brown} that if $\Gamma$ is of type
$\fp_{\infty}$, then all cohomologies of $\Gamma$ with finite type
coefficients are also of finite type and that the cohomology of
$\Gamma$ commutes with direct limits. In particular such a $\Gamma$
will be of type $(\boldsymbol{F})$.

\

\noindent
{\bf (ii)} As usual, it is convenient to try and study the cohomology of
a group $\Gamma$ trough their topological incarnation as the
cohomology of the classifying space of $\Gamma$. In particular, if
$K(\Gamma,1)$ admits a realization as a CW complex having only
finitely many cells in each dimension, it is clear that both
$\fp_{\infty}$ and $(\boldsymbol{F})$ hold for $\Gamma$

\

\noindent 
{\bf (iii)} If $\Gamma$ is an algebraically good group of type
$(\boldsymbol{F})$, then for any  
linear representation $V$ of $\Gamma^{\op{alg}}$, maybe of infinite
dimension, one has 
\[
H^{*}(\Gamma,V)\simeq H^{*}(\Gamma^{\op{alg}},V).
\]
Indeed, every linear representation of $\Gamma^{\op{alg}}$ is the
inductive limit of its  
finite dimensional sub-representations, and the Hochschild cohomology
of $\Gamma^{\op{alg}}$ always 
commutes with inductive limits.
\end{rem}

\

In the next theorem we use the notion of 
a group of Hodge type, which can be found in e.g.  \cite[$\S4$]{s2}.

\begin{thm}\label{cex1}
Let $n>1$ be an integer. Let $Y$ be a pointed and connected simplicial
set in $\mathbb{U}$ such that:
$\pi_{1}(Y,y)=\Gamma$ is an algebraically good group of type
$(\boldsymbol{F})$;
$\pi_{i}(Y,y)$ is of finite type for any $1<i<n$, and $\pi_{i}(Y,y)=0$
for $i\geq n$.

Let $\rho : \Gamma \longrightarrow Gl_{m}(\mathbb{Z})$ be
an integral representation and let
$\rho_{\mathbb{C}} : \Gamma \longrightarrow Gl_{m}(\mathbb{C})$ be the
induced complex linear representation .
Denote by $\rho_{1}, \dots, \rho_{r}$ the simple factors of the
semi-simplification of $\rho_{\mathbb{C}}$.

Let $X$ be the homotopy type defined by the following homotopy cartesian
diagram
$$\xymatrix{
Z \ar[r] \ar[d] & Y \ar[d] \\
K(\Gamma,\mathbb{Z}^{m},n) \ar[r] & K(\Gamma,1).}$$
Suppose that there exists a smooth and projective complex algebraic
variety $X$,
such that the $n$-truncated homotopy
types $\tau_{\leq n}X^{\op{top}}$ and $\tau_{\leq n}Z$ are equivalent, then
the real Zariski closure of the image of each $\rho_{j}$ is a group of
Hodge type.
\end{thm}
{\bf Proof:} The theorem is based of the following lemma, describing
the homotopy groups 
of $(Z\otimes \mathbb{C})^{\op{sch}}$.

\begin{lem}{(\cite[Prop. 4.14]{small})}
For any $i>1$ there are isomorphisms of affine group schemes
\[
\pi_{i}((Z\otimes \mathbb{C})^{\op{sch}},x)\simeq \pi_{i}(Z,x)\otimes
\mathbb{G}_{a}.
\]
\end{lem}

\smallskip

\noindent
Theorem \ref{cex1} now follows from the  previous lemma,
Corollary \ref{vhs} and \cite[Lemma $4.4$]{s2}.  \hfill $\Box$

\

\medskip

\noindent
As an immediate consequence we get:

\begin{cor}\label{ccex1}
If in theorem \ref{cex1}, the real Zariski closure of the image of one
of the representations
$\rho_{j}$ is not a group of Hodge type, then
$Z$ is not the $n$-truncation of the homotopy type of a smooth and
projective algebraic variety defined over
$\mathbb{C}$.
\end{cor}

A list of examples of representations $\rho$ satisfying the hypothesis
of the previous corollary can be obtained
using the list of \cite[$\S 4$]{s2}. Here  are two explicit examples:

\begin{itemize}
\item[A.] Let $\Gamma=\mathbb{Z}^{2g}$, and $\rho$ be any reductive integral
representation such that its
complexification $\rho_{\mathbb{C}}$ is non-unitary. Then, one of the
characters $\rho_{j}$ is not unitary, which
implies by \cite[$4.4.3$]{s2} that the real Zariski closure of its
image is not of Hodge type.
\item[B.] Let $\Gamma$ be the fundamental group of compact Riemann
surface of genus $g>2$, and let $m>2$. Consider any surjective morphism
$\rho : \Gamma \longrightarrow Sl_{m}(\mathbb{Z})\subset
Gl_{m}(\mathbb{Z})$.
Then the real Zariski closure of $\rho$ is $Sl_{n}(\mathbb{R})$ which
is not of Hodge type (see \cite[$\S 4$]{s2}).
\end{itemize}

\subsection{The formality theorem}

In this paragraph we will show that for a pointed
smooth and projective manifold $X$ the pointed schematic homotopy type
$\sch{X^{\op{top}}}{\mathbb{C}}$ is formal. This result is a generalization
of the formality result of \cite{dgms}, and could possibly be used to
obtain other restrictions on homotopy types of projective manifolds.
We leave this question for the future, and restrict ourselves to present
a proof of the formality theorem.

To state the main definition of this paragraph, let us recall that
for any affine group scheme $G$, and any
$G$-equivariant commutative differential graded algebra $A$, or any
$G$-equivariant co-simplicial algebra $A$,
one has the cohomology algebra $H^{\bullet}(A)$. This algebra is in a
natural way
a graded algebra with an action of $G$, and therefore can be considered as
a $G$-equivariant commutative differential graded algebra with trivial
differential.

\begin{df}\label{dfo}
Let $G$ be an affine group scheme.
\begin{itemize}
\item A $G$-equivariant commutative differential graded algebra $A$ is
$G$-formal if it is isomorphic to $H^{\bullet}(A)$ in
$\op{Ho}(\opi{G-\sff{CDGA}})$.
\item A $G$-equivariant co-simplicial algebra $A$ is $G$-formal if
it is isomorphic in $\op{Ho}(\opi{G-\sff{Alg}}^{\Delta})$ to the
denormalization 
$D(H^{\bullet}(A))$.
\item A pointed schematic homotopy type $F$ is formal, if it is
isomorphic
in $\op{Ho}(\sff{SPr}(\mathbb{C}))$ to an object of the form
$[\mathbb{R}\op{Spec}_{G}\, A/G]$, whith $A$ a $G$-formal $G$-equivariant
co-simplicial algebra,
and $G$ an affine reductive group scheme.
\end{itemize}
\end{df}

Let $X$ be a  pointed and connected simplicial set in $\mathbb{U}$, and
$G_{X}:=\pi_{1}(X,x)^{\op{red}}$ be the
pro-reductive completion of its fundamental group. One can
consider the $G_{X}$-equivariant
commutative differential graded algebra
$H^{\bullet}(C^{\bullet}(X,\mathcal{O}(G_{X})))$, of cohomology
of $X$ with coefficients in the local system of
algebras $\mathcal{O}(G_{X})$. As this
cohomology algebra is such that
$H^{0}(C^{\bullet}(X,\mathcal{O}(G_{X})))\simeq \mathbb{C}$,
Proposition \ref{p5} implies that
$[\mathbb{R}\op{Spec}_{G_{X}}\, H^{\bullet}(X,\mathcal{O}(G_{X}))/G_{X}]$
is a pointed schematic homotopy type.

\begin{df}
The formal schematization of the pointed and connected simplicial set
$X$ is the pointed schematic
homotopy type
\[
(X\otimes \mathbb{C})^{\op{for}}:=[\mathbb{R}\op{Spec}_{G_{X}}\,
H^{\bullet}(X,\mathcal{O}(G_{X}))/G_{X}]
\in \op{Ho}(\sff{SPr}_{*}(\mathbb{C})).
\]
\end{df}

The main theorem of this section is the following formality
statement, which in particular answers Problem~2 in \cite[\S 7]{ght}.

\begin{thm}\label{tfo}
Let $X$ be a pointed smooth and projective complex manifold and
$X^{\op{top}}$ its underlying
topological space. Then, there exist
an isomorphism in $\op{Ho}(\sff{SPr}_{*}(\mathbb{C}))$, functorial in $X$
\[
\sch{X^{\op{top}}}{\mathbb{C}} \simeq (X^{\op{top}}\otimes
\mathbb{C})^{\op{for}}.
\]
\end{thm}
{\bf Proof:} By using Proposition \ref{pdiff},
it is enough to produce a functorial isomorphism in \linebreak 
$\op{Ho}(\opi{G_{X}-\sff{CDGA}})$
\[
(A^{\bullet}_{DR}(\mathcal{O}(G_{X}),\nabla) \simeq
H^{\bullet}_{DR}(X,\mathcal{O}(G_{X})),
\]
where $H^{\bullet}_{DR}(X,\mathcal{O}(G_{X})) \in
\op{Ho}(\opi{G_{X}-\sff{CDGA}})$ is the cohomology algebra of
$(A^{\bullet}_{DR}(\mathcal{O}(G_{X}),\nabla)$.

But, writing $\mathcal{O}(G_{X})$ as an $\opi{Ind}$-object $(V,\nabla) \in
T_{DR}(X)$,
corresponding to $(V,\nabla,D'') \in T_{D'}(X)$,
and applying
Corollary \ref{cfo} one obtains a diagram of quasi-isomorphisms
\[
\xymatrix{
(A^{\bullet}_{DR}(\mathcal{O}(G_{X})),D) & \ar[l]
(A^{\bullet}_{D'}(V),D'') \ar[r] &
(H^{\bullet}_{DR}(X,\mathcal{O}(G_{X})),0)).}
\]
By the compatibility of the functor $(V,\nabla) \mapsto
(A^{\bullet}_{D'}(V),D'')$ with the tensor products,
this diagram is actually a diagram of $G_{X}$-equivariant
commutative differential graded algebras. Passing to the homotopy
category, one obtains
well defined and functorial isomorphisms in $\op{Ho}(\opi{G_{X}-\sff{CDGA}})$
\[
(A^{\bullet}_{DR}(\mathcal{O}(G_{X})),\nabla) \simeq
H^{\bullet}_{DR}(X,\mathcal{O}(G_{X})).
\]
The theorem is proven. \  \hfill $\Box$

\

\bigskip

\begin{rem} \label{rem-formality} The above theorem implies in
particular that for a smooth projective variety $X$ the
algebraic-geometric invariants we have considered in this paper are
captured not only by the cochain dg algebra by also by the
corresponding graded algebra of cohomology. More precisely, if $G$
denotes the pro-reductive completion of $\pi_{1}(X)$, and ${\mathcal
O}(G)$ denotes the algebra of functions on $G$ viewed as a local
system on $X$, then the schematization $\sch{X}{\mathbb{C}}$ and in
particular the support invariants of $X$ are all determined by the
$G$-equivariant graded algebra of cohomology $H^{\bullet}(X,{\mathcal
O}(G))$. We do not know if, as for the rational homotopy type, the
$\mathbb{C}^{\times}$ action on $\sch{X}{\mathbb{C}}$ is also
determined by the its action on the group $G$ and the $G$-equivariant
graded algebra $H^{\bullet}(X,\mathcal{O}(G))$.
\end{rem}

The functoriality assumption in theorem \ref{tfo} has the following
striking consequence.

\begin{cor}\label{ccfo}
Let $f,g : (X,x) \longrightarrow (Y,y)$ two morphisms between pointed
smooth and projective complex manifolds. Suppose that
\begin{enumerate}
\item The induced morphisms $f,g : \pi_{1}(X,x)^{\op{red}}
  \longrightarrow \pi_{1}(Y,y)^{\op{red}}$ 
are equal.
\item For any simple local system $L$ on $Y$ the induces morphisms
\[
f,g : H^{*}(Y,L) \longrightarrow H^{*}(X,f^{*}(L)=g^{*}(L))
\]
 are
equal.
\end{enumerate}
Then the two morphisms of pointed stacks $f,g : (X^{\op{top}}\otimes
\mathbb{C},x)^{\op{sch}} \longrightarrow (Y^{\op{top}}\otimes
\mathbb{C},y)^{\op{sch}}$ are equal as morphisms in
$\op{Ho}(\sff{Spr}_{*}(\mathbb{C}))$. 
\end{cor}

\subsection{Hurewitz maps} 

In this section we use the Hodge decomposition on the schematic
homotopy type of a smooth projective $X$ in order to show that under
certain conditions the image of the Hurewitz map 
\[
\pi_{n}(X) \to
H_{n}(X,{\mathbb Z})
\]
 is a sub Hodge structure. For $n =2$ this
answers a question of P.Eyssidieux.

\begin{thm} \label{philippe} Suppose  $X$ is smooth and projective over
  ${\mathbb C}$ satisfying
\begin{enumerate}
\item[(a)] $\pi_{1}(X)$ is good,
\item[(b)] $\pi_{i}(X)$ is finitely generated for $1 < i < n$.
\end{enumerate}
Then $\op{Im}[\pi_{n}(X) \to
H_{n}(X,{\mathbb Z})]$ is a Hodge substructure.
\end{thm}
{\bf Proof.} The natural map $\bs : X \to \sch{X}{\mathbb{C}}$
induces a commutative diagram 
\begin{equation} \label{eq:hurewitz}
\xymatrix{
\pi_{n}(X)\otimes {\mathbb C} \ar[r]^-{\bs} \ar[rd]^-{\psi} \ar[d]_-{\op{Hu}} &
   \pi_{n}(\sch{X}{\mathbb{C}}) \ar[d]^-{\op{Hu}^{sch}} \\
H_{n}(X,{\mathbb C}) \ar[r]^-{\cong} & H_{n}((X\otimes {\mathbb
   C})^{\op{sch}},{\mathcal O}) 
}
\end{equation}
Since the schematic Hurewitz map is functorial for morphisms in the
  homotopy category of simplicial presheaves and since $\sch{X}{\mathbb{C}}$ can be considered as ${\mathbb
  C}^{\times}$-equivariant object in this category,
  it follows that the image of $\op{Hu}^{\op{sch}}$ is preserved under the
  ${\mathbb C}^{\times}$  action on $H_{n}(X,{\mathbb C})$. Therefore
  it suffices to show that the maps $\op{Hu}^{\op{sch}}$ and $\psi$ have
  the same image.

By the commutativity of \eqref{eq:hurewitz} we know that
$\op{im}(\psi) \subset \op{im}(\op{Hu}^{\op{sch}})$ and so we will be done
if we can show the opposite inclusion $\op{im}(\psi) \supset
\op{im}(\op{Hu}^{\op{sch}})$. For this we will need the following 

\begin{lem} Let $X$ be a smooth complex projective variety satisfying
  conditions (a) and (b). For every complex finite dimensional
  representation $V$ of 
  $\pi_{1}(X)$, the morphism $\bs$ induces a bijection
\[
\bs^{*} : \op{Hom}_{\pi_{1}(X)}^{\op{cont}}(\pi_{n}(\sch{X}{\mathbb{C}}),V) \widetilde{\to} \op{Hom}_{\pi_{1}(X)}(\pi_{n}(X),V).
\]
where $\op{Hom}_{\pi_{1}(X)}(\pi_{n}(X),V)$ denotes the space of
$\pi_{1}(X)$-equivariant morphisms between the abelian groups
$\pi_{n}(X)$ and $V$, and $\op{Hom}_{\pi_{1}(X)}^{\op{cont}}(\pi_{n}(\sch{X}{\mathbb{C}}),V)$  denotes the space of
  $\pi_{1}(X)$-equivariant cotinuous morphisms between the linearly
  compact vector spaces $\pi_{n}(\sch{X}{\mathbb{C}})$ and $V$. 
\end{lem}
{\bf Proof.} Let $\tau_{\leq n-1} : X \to X_{\leq n-1}$ be the
Postnikov truncation. Consider the commutative diagram
\[
\xymatrix{
X \ar[r]^-{\tau_{\leq n-1}} \ar[d]_-{\bs} & X_{\leq n-1} \ar[d]^-{\bs}
\\
\sch{X}{\mathbb{C}} \ar[r]^-{\tau_{\leq n-1}} & (\sch{X}{\mathbb{C}})_{\leq n-1}
}
\]
The vertical maps in this diagram induce a morphism between the five term
long exact sequences coming from to the Leray spectral sequence of the
truncation maps:
\[
\xymatrix@C-1pc{
H^{n}((\sch{X}{\mathbb{C}})_{\leq n-1},V) \ar[d]^-{(1)} \ar[r] & H^{n}(\sch{X}{\mathbb{C}},V) \ar[d]^-{\cong} \ar[r] &
  \op{Hom}_{\pi_{1}(X)}^{\op{cont}}(\pi_{n}((X\otimes 
  {\mathbb C})^{\op{sch}}),V) \ar[d]^-{\bs^{*}}  \ar[r] & \\ 
H^{n}(X_{\leq n-1},V) \ar[r] & H^{n}(X,V)  \ar[r] &
  \op{Hom}_{\pi_{1}(X)}(\pi_{n}(X,V)  \ar[r] &  
\\  
\ar[r] & H^{n+1}(\sch{X}{\mathbb{C}})_{\leq n-1},V)
\ar[d]^-{(2)} \ar[r] 
\ar[r] & H^{n+1}((X\otimes 
  {\mathbb C})^{\op{sch}},V) \ar[d]^-{\cong}  \\
\ar[r] & H^{n+1}(X_{\leq n-1},V) \ar[r] & H^{n+1}(X,V)  
}
\]
To show that $\bs^{*}$ is an isomorphism, it suffices to show that (1)
and (2) are isomorphisms. This will follow if we know that the natural
map $X_{\leq n-1} \to  (\sch{X}{\mathbb{C}})_{\leq n-1}$ induces an isomorphism between the
  schematic homotopy types $(X_{\leq n-1}\otimes
  {\mathbb C})^{\op{sch}}$ and $(\sch{X}{\mathbb{C}})_{\leq n-1}$. This is equivalent to showing that
  the stack $(X_{\leq n-1}\otimes
  {\mathbb C})^{\op{sch}}$ is $(n-1)$-truncated, which in turn follows from
  the hypothesis (a) and (b) and \cite[Proposition~4.21]{small} applied to the
  moprphism $X_{\leq n-1} \to X_{\leq 1}$. The lemma is proven. \
  \hfill $\Box$

\

\smallskip

\noindent
To finish the proof of the theorem take $V = H_{n}(X,{\mathbb
  C})/\op{im}(\psi)$. Now the composition 
\[
\xymatrix@1{
\pi_{n}(X)\otimes {\mathbb C} \ar[r]^-{\bs} & \pi_{n}((X\otimes
   {\mathbb C})^{\op{sch}}) \ar[r]^-{\op{Hu}^{\op{sch}}} &
H_{n}(X,{\mathbb 
  C}) \ar[r] &  V
}
\]
is zero by construction, the lemma implies that  the composition 
 \[
\xymatrix@1{
\pi_{n}((X\otimes
   {\mathbb C})^{\op{sch}}) \ar[r]^-{\op{Hu}^{\op{sch}}} &
H_{n}(X,{\mathbb 
  C}) \ar[r] &  V
}
\]
is zero as well. In particular $\op{im}(\op{Hu}^{\op{sch}}) \subset
\op{im}(\psi)$ which completes the proof of the theorem. \ \hfill
$\Box$

\

\bigskip

\begin{rem} It is not hard to construct interesting examples of
  varieties $X$ satisfying the hypothesis of Theorem
  \ref{philippe}. For instance we can start with a variety $Z$ which
  is  a $K(\pi,1)$ with 
  good fundamental group, e.g. $Z$ can be a product of curves or and
  abelian variety, and take $X$ to be a smooth hyperplane section of
  dimension $n$. 
\end{rem}

\

\bigskip

\begin{rem} Note that the lemma implies that
  $\pi_{n}(\sch{X}{\mathbb{C}})$ is the $\pi_{1}(X)$-equivariant
  unipotent 
  completion of $\pi_{n}(X)$. In particular, the image of $\pi_{n}(X)$
  in $\pi_{n}(\sch{X}{\mathbb{C}})$ is Zariski dense.
\end{rem}

\

\bigskip

\noindent
\textsf{L. Katzarkov: UCI, Irvine, CA 92612, e-mail: lkatzark@math.uci.edu}

\

\medskip

\noindent
\textsf{T. Pantev: University of Pennsylvania, Philadelphia, PA 19104,
e-mail: tpantev@math.upenn.edu}

\

\medskip

\noindent
\textsf{B. Toen: Laboratoire Emile Picard, Universit\'e Paul Sabatier,
  118 route 
de Narbonne, 31062 Toulouse, France, e-mail: toen@picard.ups-tlse.fr}

\end{document}